\documentclass [twoside,10pt]{article}
\usepackage{graphicx}
\usepackage{amssymb}
\usepackage{amsmath}
\usepackage{color}
\usepackage{url}
\usepackage{hyperref}

\newfont{\bbb}{msbm10 scaled\magstep 1} 

\makeatletter
\g@addto@macro\@floatboxreset\centering
\makeatother

\newcommand{\diam}{{\rm diam}}

\newcommand{\bN}{{\mathbb N}}
\newcommand{\bR}{{\mathbb R}}

\newcommand{\gO}{\mathsf{O}}

\newcommand{\gGL}{\mathsf{GL}}

\newcommand{\gSL}{\mathsf{SL}}

\newcommand{\ve}{\varepsilon} 
\newcommand{\Tt}{T} 

\DeclareMathOperator{\conv}{\mathrm{conv}}
\DeclareMathOperator{\vol}{vol}
\DeclareMathOperator{\Tr}{Tr}

\newtheorem{thm}{Theorem}
\newtheorem{cor}{Corollary}

\newtheorem{ddd}{Definition}
\newtheorem{pro}{Proposition}

\newtheorem{prob}{Problem}

\title{\bf Banach-Mazur distance and related methods of comparing convex bodies}

\author{Marek Lassak and Stanis{\l}aw J. Szarek}

\begin{document}
\date{}
\maketitle
\centerline
{ }
\pagestyle{myheadings} \markboth{ }{ }

\baselineskip 16pt

\bigskip 
\noindent
{\bf Abstract} \ 
The Banach-Mazur distance is one of the standard tools used in the study of mathematical objects such as convex bodies or normed spaces, and it becomes indispensable when investigating properties that are invariant under linear isomorphisms. It appears in convex geometry and convex analysis as a way to quantify the approximation of convex bodies. In functional analysis, it plays a basic role in the study of local properties of Banach spaces. 
The goal of this survey is to present a cohesive picture of these two perspectives. We supply an extensive collection of known results for reference and list open problems of varying difficulty, which could stimulate further research. We include in our discussion several distinct variations on the concept of the Banach-Mazur distance, some of them exploratory, that haven’t received much attention to date in the literature. We follow the trend -- which emerged in recent decades -- of considering arbitrary convex bodies, and not just those that are centrally symmetric.  Finally, we include a comprehensive list of references on the Banach-Mazur distance.



\newpage
\tableofcontents

\baselineskip 19.5pt

\newpage
\section{Introduction} 
\label{sec:intro} 

Let $E^n$ be the $n$-dimensional Euclidean space, which we will often identify with the {\em vector} space $\bR^n$.  
By ${\mathcal C}^n$ we will always denote the family of convex bodies  of $E^n$ (i.e., compact convex subsets with nonempty interior). 
Similarly,   ${\cal M}^n$ will represent the family of centrally symmetric convex bodies of $E^n$. 
The letter $B$ will usually stand for a Euclidean ball,  $P$ for a parallelotope, and $\Tt$ for a triangle, tetrahedron or, more generally, for a simplex. If we want to emphasize the dimension, we may write (for example) $B^n$. 
If $C\in {\cal C}^n$, then by ${\rm int}(C)$ we will denote the interior of $C$ and by $\partial C$ its boundary. 
 
\vskip 0.09cm
We say that convex bodies $C_1$ and $C_2$ are {\it affinely equivalent} if there exists an affine transformation $a$ such that $a(C_1) =C_2$. 

\vskip 0.09cm
\begin{ddd} \label{def:dBM}
The Banach-Mazur distance (or the BM-distance, for brevity) between convex bodies $C$ and $D$ of $E^n$ is the quantity 
\begin{equation} \label{eq:dBM}
d_{BM} (C,D) = \inf_{a, h_{\lambda}} \{ \lambda : \, D \subset a(C) \subset h_{\lambda}(D) \}, 
\end{equation} where $a$ is an affine transformation and $h_{\lambda}$ is a homothety with a {\sl positive}  ratio $\lambda$. 
\end{ddd}

Since the dimension is finite, a simple compactness argument shows that the infimum is attained and it is  readily verified that $d_{BM} (C,D)  = d_{BM} (D,C)$. 
While not assumed explicitly, the condition in \eqref{eq:dBM} forces $\lambda \geq 1$ and the center of the homothety $h_{\lambda}$ to be in $D$. An equivalent -- and sometimes convenient -- form of the condition from \eqref{eq:dBM} is as follows: {\sl $a(C)$ contains  $D$ and is contained in a translate of $\lambda D$}. 
(It then follows that $D$ and the translate of $\lambda D$ are related by a homothety with ratio $\lambda$ and center in $D$.) 
Another equivalent form of the condition is $a(D) \subset C \subset h_{\lambda}(a(D))$ (for some other $a, h_\lambda$). 
 
\vskip 0.09cm
The above concept was defined in the Banach's monograph \cite{Ban32} in the essentially equivalent language of normed spaces 
(Remarques, Ch. XI, \S 6, p. 242 in the 1932 French language edition, or p. 151 in the 1987 English language edition; it is stated there that the notions were introduced ``in collaboration with S. Mazur''). 
If $X,Y$ are two isomorphic normed spaces, one defines their distance as 
\begin{equation} \label{def:dBMns}
d(X,Y) = \inf \|u\| \cdot \|u^{-1}\|, 
\end{equation} 
where $u$ ranges over all linear isomorphisms from $X$ to $Y$.  (In the infinite-dimensional case, which is beyond the scope of this survey, the infimum may or may not be attained.) 
Denoting by $B_X$ and $B_Y$ the (closed) unit balls of $X$ and $Y$, we have then
\begin{equation} \label{d-d_BM}
 d(X,Y) = d_{BM} (B_X,B_Y).    
\end{equation}
Indeed, since  $\|u^{-1}\|^{-1} B_Y \subset uB_X \subset \|u\| B_Y$ by the definition of the norm of an operator, the only point that needs to be addressed is that in \eqref{eq:dBM}  one considers all {\sl affine} maps/homotheties, while  
in \eqref{def:dBMns} only {\sl linear} operators are allowed.  
However, it is a simple exercise to verify that if $C,D$ are $0$-symmetric, then the infimum in \eqref{eq:dBM} is attained on linear maps. 

Definitions in the spirit of \eqref{d-d_BM}, applying to not-necessarily-symmetric convex bodies, have been used in an {\em ad hoc} fashion at least since the middle of the 20th century. 
In 1963, Gr\"unbaum formally introduced in \cite{[Grun63]} a slightly different concept, allowing homotheties with negative ratios (see Definition \ref{grunbaum} in Section \ref{sec:simplex}). A definition equivalent to \eqref{d-d_BM} was formulated for example in \cite{[L 92a]} in 1992. 
 
\vskip 0.09cm
One checks easily that the multiplicative version 
\begin{equation} \label{triangle_ineq}
 d_{BM} (C, D) \leq d_{BM} (C, A) \,\cdot \,d_{BM} (A, D)  
\end{equation}
of the triangle inequality holds  for any convex bodies $C, D, A \subset E^n$.  
To have true metric spaces, one would need to work with $\log d_{BM}$ as the distance function and to identify convex bodies that are affinely equivalent. 
We will denote the resulting set of equivalence classes of $n$-dimensional convex bodies by $\mathfrak{C}^n$  
(and, respectively, $\mathfrak{M}^n$ for centrally symmetric bodies), but -- unless the distinction is meaningful 
-- we will follow the usual practice of glossing over these and related subtleties.  

An important family of convex bodies is induced by the classical $\ell_p$-norms. Given  $x= (x_j) \in \bR^n$. we set $\|x\|_p := \left(\sum_{j=1}^n |x_j|^p\right)^{1/p}$ for $p\in [1,\infty)$ and $\|x\|_\infty := \max_j |x_j|$. The unit ball with respect to $\ell_p$-norm is often denoted by  $B_p^n$ and the associated normed space by $\ell_p^n$.  Obviously, $B_\infty^n$ is an $n$-dimensional cube, an instance of a parallelotope. The $B_1^n$, or any body affinely equivalent to it, is called a {\em cross-polytope}.

\vskip 0.09cm
A key role in the study of these concepts is played by duality considerations,  particularly in the case  of 
centrally symmetric bodies or normed spaces, where we have the identity
\begin{equation} \label{eq:duality}
d(X,Y) =  d(X^*,Y^*),     
\end{equation} 
which is readily verified via the definition \eqref{def:dBMns}.  
In the more general setting,   
one defines the {\em polar} of  $A \subset E^n$ (identified with the vector space $\bR^n$) by 
\begin{equation} \label{def:polar}
A^\circ = \{x \in \bR^n :  \langle x,y\rangle \leq 1 \hbox{ for all } y \in A\}.   
\end{equation} 
The archetypal example is  that of $A = B_X$, the unit ball with respect to some norm, in which case $A^\circ = B_{X^*}$, the unit ball with respect to the dual norm.  In particular, 
\begin{equation} \label{ellp-ellq}
(B_p^n)^\circ = B_q^n, \hbox{ where } \tfrac 1p + \tfrac 1q =1.
\end{equation}
Specifying to $p=\infty$, we conclude that the polar of a ($0$-symmetric) parallelotope is a cross-polytope.

An important elementary result (based on a baby version of the Hahn-Banach theorem) is the {\sl bipolar theorem}:  {\sl if $A$ is closed, convex and $0\in A$, then $(A^\circ)^\circ = A$.} 

\vskip 0.09cm
In the above notation, \eqref{eq:duality} can be restated as follows:  {\sl if  $C, D \in {\cal M}^n$, then } 
\begin{equation} \label{eq:polar}
d_{BM} (C, D) = d_{BM} (C^\circ, D^\circ).  
\end{equation}
 
It is tempting to ask whether  \eqref{eq:polar}  holds also for $C, D \in {\cal C}^n$, but a moment of reflection shows that there is a major problem due to the fact that the center of the homothety implicit in the definition of  $d_{BM}$ is not necessarily $0$. 
What follows directly is that if  $a, h_\lambda$ are such that the infimum in \eqref{eq:dBM} is attained, $p$ is the center of the homothety $h_\lambda$, and  $q= a^{-1}(p)$, then  $d_{BM} \big((C-p)^\circ, (D-q)^\circ\big) \leq d_{BM} (C, D)$. 
 
\vskip 0.09cm
 Clearly, the equality can not hold without {\em some} additional assumptions. For example, if $C \in {\cal M}^n$ and $D=C$ (hence 
 $d_{BM} (C, D) = 1$), it is possible that one of $(C-p)^\circ, (D-q)^\circ$ is centrally symmetric and the other is not (hence $d_{BM} \big((C-p)^\circ, (D-q)^\circ\big) >1$).  
\begin{prob}
It would be interesting to find reasonably general conditions under which the equality {\em does} hold. 
\end{prob}

The concept of the Banach-Mazur distance and its properties were extensively studied in the monograph of Tomczak-Jaegermann \cite{[Tom 89]}.  
Since 1989, when it was published, many new important results have been obtained. 
Other surveys and books addressing the subject are (in the chronological order) \cite{Pelc84,  Grub, LiMi, [Thompson], [JL], [BMP],  [Szar06], Ver11, ArtGiaMi-1, To15, [Gi-EoM], ArtGiaMi-2, BMT}. 
However, most of these works concentrate on the setting of normed spaces or, equivalently, 
on that of centrally symmetric convex bodies, and often on results that are meaningful only when the dimension is large.  
Here we consider all settings, including a selection of questions in low dimensions and those involving 
not-necessarily-centrally-symmetric bodies. 

There is a broad spectrum of mathematical tools that are used to calculate or estimate BM-distances, from elementary (but often ingenious) geometric/calculus arguments for some low-dimensional questions (for example, Theorem \ref{M2}), through subtle variational arguments (for example, Theorems \ref{john}, \ref{decomposition} and \ref{thm:GLMP}) to links to areas such as combinatorics or number theory (for example, Section \ref{para-cross}). Next, one employs a variety of invariants that allow to differentiate between bodies and lower-bound the distances, see  Sections \ref{sec:projection} and  \ref{sec:weak} for representative examples. Somewhat surprisingly, some of the deepest results were obtained using random methods, both for constructing examples/counterexamples (e.g., Theorem \ref{gluskin}) and for finding affine transformations figuring in \eqref{eq:dBM}  (e.g., Theorem \ref{bizeul-klartag}).

\vskip0.4cm

\noindent
\section{John's theorem, its consequences, and variants : approximating convex bodies by ellipsoids}   \label{sec:john} 

The arguably most famous result that can be phrased in terms of the Banach-Mazur distance is due to John \cite{[John]}.

\begin{thm}[John's theorem] \label{john}
If $K$ is a convex body of $E^n$, then there exists a $0$-symmetric ellipsoid $E$ and $z\in K$ such that $E \subset K-z \subset n E$.  Consequently, $d_{BM}(K,B) \leq n$. Moreover, if $K$ is centrally symmetric, then $d_{BM}  (K,B) \leq \sqrt{n}$.
\end{thm}
 
The fact that the first statements of Theorem \ref{john} implies the second is immediate, it follows from $d_{BM}(E,B) = 1$ and the triangle inequality \eqref{triangle_ineq}. 
(Note that the two statements are not formally equivalent, see the discussion of the estimate  \eqref{john-ctr} in Section \ref{sec:one-symmetric}.) 
Both bounds are optimal.  First, $d_{BM} (\Tt,B) =n$, where $\Tt$ is the $n$-dimensional simplex. 
Additionally, Leichtweiss \cite{Lei} and Palmon \cite{[Pal 92]} proved that the equality $d_{BM} (A,B) =n$ holds only if $A$ is a simplex.  (Some authors omit the first -- and much earlier -- reference, perhaps because Leichtweiss's paper was in German; 
see Section \ref{sec:simplex} for these and more results concerning the functional $d_{BM} (\Tt,\cdot)$.) 
Next, in the centrally symmetric case, we have $d_{BM} (M,B) = \sqrt n$ when $M$ is a parallelotope or a cross-polytope. A more detailed analysis of the bodies whose BM-distance to is extremal or near-extremal can be found in Section \ref{sec:extreme}.  
 
\vskip 0.09cm
Before proceeding, let us comment on the proof of John's theorem.  
Let $E_J$ be the ellipsoid of maximal volume contained in $K$. It can be shown that such ellipsoid is unique;  
it is usually called the John ellipsoid of $K$. 
It then follows that if $z$ is the center of  $E_J$, then  $K-z \subset n (E_J-z)$.   If $K$ is $0$-symmetric, it similarly follows that $z=0$ and that  $K \subset \sqrt n \,E_J$. 
An alternative ``dual'' construction is as follows:  if $E_L$ is the ellipsoid of minimal volume containing $K$ (called L\"owner or John-L\"owner ellipsoid, and again unique; see \cite{Henk} for historical background and more references) and $z$ is its center, then $n^{-1} (E_L-z)\subset K-z$, with an analogous improvement if  $K$ is centrally symmetric. While our presentation is focused on analyzing the properties of $E_J$, all the facts we state can be {\em mutatis mutandis} reformulated for $E_L$. 
 
\vskip 0.09cm
The properties of John and L\"owner ellipsoids are usually established via variational arguments of various level of sophistication 
(see \cite{Ball92, Ball97, Bar02} or Section 4.2 in \cite{AuSz17}). 
For the benefit of a non-expert reader, here is a high-level outline of such an argument. 
 
\vskip 0.09cm
First, the existence of John's ellipsoid. Since an ellipsoid is an affine image of the Euclidean ball $B$, we are looking for $E = z+ uB \subset K$, $u\in \gGL(n)$ such that $|\det u\hskip0.4mm|$ is maximal. 
Note that the invertibility requirement is nonessential (if $u$ is not invertible, then $\det u =0$), and if we 
only require that $u$ be a linear operator, then $\mathcal{E} := \{(z,u) : z+ uB \subset K \}$ is a compact convex set, on which the supremum of $|\det u\hskip0.4mm|$ is attained. 
 
\vskip 0.09cm
For uniqueness, we observe that, without loss of generality, we may assume that $u$ is positive semi-definite ($u \in PSD$). (If not, consider the polar decomposition $u=u'V$ such that $u' \in PSD$ and $V$ is an isometry, and replace $u$ by $u'$.) We next use the fact that the function $\phi(u) = \log (\det u)$ is {\em strictly} concave on the set of positive definite operators, which implies that the maximum of $\phi$ on can be attained for only one $u=u_0$. It remains to exclude the possibility that the maximum of $\vol E$ is attained at $(z_1,u_0)$ and $(z_2,u_0)$ with $z_1\neq z_2$, which is easy :  in that case, $\conv\big((z_1+ u_0B)\cup (z_2+ u_0B\big) \subset K$ clearly contains an ellipsoid whose volume is larger than that of $u_0B$ (start by considering the case $u_0=I$ and then generalize). 
 
\vskip 0.09cm
The key ingredient of the argument is the following fact, arguably more fundamental than Theorem \ref{john}, which is often referred to as ``John's decomposition of the identity''. 
To state it, it is convenient to introduce the following concepts. 
First, for a convex body $L \subset K$, we say that $y$ is their {\em contact point} if $y \in L \cap \partial K$ (then, of course, necessarily $y \in \partial L \cap \partial K$). 
Next, we say that $K$ is in the {\em John position} provided its John ellipsoid $E_J$ is the unit Euclidean ball $B$. (This can be always achieved by replacing $K$ by an appropriate affine image.) 
We then have 

\begin{thm} [John's decomposition] \label{decomposition}
Let $K \subset E^n$ be in the John position. 
Then there exists a  set $\{y_1, y_2,\ldots, y_N\}$ of contact points of $B$ and $K$ and a sequence of positive scalars $(\lambda_1, \lambda_2, \ldots , \lambda_N)$ such that $\sum_{k=1}^N \lambda_k y_k = 0$ and, for every $x\in E^n$,
\begin{equation}  \label{eq:john}
 x = \sum_{k=1}^N \lambda_k  \langle y_k, x \rangle y_k .   
\end{equation}
Moreover, we can assume that $N \leq \frac{n(n+3)}2$ (or $N \leq \frac{n(n+1)}2$ if $K$ is centrally symmetric). \\
Conversely, if $B\subset K$, the existence of positive scalars $(\lambda_k)$ and contact points $(y_k)$ verifying \eqref{eq:john} guarantees that 
$K$ is in the John position.
\end{thm}

The equation \eqref{eq:john} can be rewritten as 
\begin{equation} \label{eq:john2}
I_n= \sum_{k=1}^N \lambda_k \, y_k y_k^t,
\end{equation}
where $I_n$ is the $n\times n$ identity matrix and $y^t$ denotes the transpose of $y$. This justifies the ``decomposition of identity'' language. 
A reasonably detailed proofs of Theorem \ref{decomposition} can be found in Barvinok's book \cite{Bar02} or Ball's survey \cite{Ball97}, but here is an outline of the main idea. 
 
\vskip 0.09cm
Since necessarily $\sum_{k=1}^N\lambda_k =n$, \eqref{eq:john2} is equivalent to  
\begin{equation} \label{eq:john3}
I_n/n \in \conv \{yy^t : y \hbox{ is a contact point of } B \hbox{ and } K \} =: C_K, 
\end{equation} 
with the bounds on $N$ following from Carath\'eodory's theorem. 

We now argue by contradiction. If $I_n/n$ doesn't belong to $V_K$, it can be (strictly) separated from it be a linear functional. Given that both $I_n/n$ and $V_K$ lie in the affine space of trace one self-adjoint matrices, this means that there exist $\ve >0$ and a self-adjoint matrix $v$ with $\Tr v =0$ and $\Tr (vyy^t) = \langle vy,y\rangle \geq \ve$ for all contact points $y \in \partial B \cap \partial K$. 
It is then elementary to verify that if $t>0$ is small enough, the ellipsoid defined by $E := \{x : \langle (I_n +tv)x,x\rangle \leq 1\}$ is contained in $K$ while $\vol E > \vol B$, contradicting the hypothesis $B=E_J$. The condition $\sum_{k=1}^N \lambda_k y_k = 0$ can be obtained by making the above argument slightly more complicated (working in the space of affine maps rather than in that of linear maps, and separating $(0,I_n/n)$ from $\conv\{(y,yy^t)\}$). However in the most important case when $K$ is $0$-symmetric, it can be achieved for free (or almost for free, at the cost of doubling $N$): just replace each $(\lambda_k,y_k)$ with a pair $(\lambda_k/2,y_k)$, $(\lambda_k/2,-y_k)$. A trick allowing to formally deduce the general case from the centrally symmetric case is described in \cite{AuSz17} (see Propositions 4.4, 4.6, and 4.7 there). 
 
\vskip 0.09cm
We will now outline a derivation of Theorem \ref{john} from Theorem \ref{decomposition}. Assume, as we may, that $K$ is in the John position. 
The key observation is that if $y$ is a contact point, then $\langle x,y\rangle \leq 1$ for all $x \in K$. If $K$ is symmetric, this implies that $\big|\langle x,y\rangle\big|\leq 1$ and so, by \eqref{eq:john},

\begin{equation} \label{eq:john4}
| x |^2 = \langle x,x \rangle = \Big\langle x,\sum_{k=1}^N \lambda_k  \langle y_k, x \rangle y_k \Big\rangle =  \sum_{k=1}^N \lambda_k \langle y_k, x \rangle^2 \leq \sum_{k=1}^N \lambda_k =n,
\end{equation}
where $ | \cdot | =\|\cdot \|_2$ is the Euclidean norm. 
So, if $x\in K$, then $| x |$ does not exceed $n^{1/2}$. In other words, $K \subset n^{1/2} B$, as asserted. 
 
\vskip 0.09cm
The argument in the not-necessarily-symmetric case is a little more involved due to the fact that the terms in the sum in \eqref{eq:john4} have to be handled differently depending on their sign. 
Given $x\in K$, denote $S = \{k : \langle y_k, x \rangle \geq 0\}$ and $s= \sum_{k\in S} \lambda_k$. Note that $ \sum_{k\not\in S} \lambda_k \big|\langle y_k, x \rangle\big| =  \sum_{k\in S} \lambda_k \langle y_k, x \rangle \leq s$ (the equality because $\sum_{k=1}^N \lambda_k y_k = 0$ and the inequality because $0\leq \langle y_k, x \rangle \leq 1$ for $k\in S$). 
Accordingly (cf. \eqref{eq:john4}), 

$$| x |^2 = \sum_{k=1}^N \lambda_k \langle y_k, x \rangle^2 = \sum_{k\in S} \lambda_k \langle y_k, x \rangle^2 +\sum_{k\not\in S} \lambda_k \langle y_k, x \rangle^2 \leq s + s|x|,$$
which resolves to $s\geq |x|^2/(|x|+1)$.  A similar but simpler bound, based just on $\sum_{k\not\in S} \lambda_k = n-s$ and $\big\langle y_k, x \rangle\big| \leq |x|$, leads to $| x |^2 \leq s + (n-s)|x|^2$, which -- assuming $|x|>1$ resolves to $s \leq (n-1)|x|^2/(|x|^2-1)$.  For these two bounds on $s$ to be compatible, we must have $|x|^2/(|x|+1) \leq (n-1)|x|^2/(|x|^2-1)$, which is equivalent to $|x|\leq n$. Similarly as before, this means that $K \subset n B$, as asserted. 
 
\vskip 0.09cm
An immediate consequence of John's Theorem \ref{john} and the multiplicative triangle inequality \eqref{triangle_ineq} is the following. 
\begin{cor} \label{john2} We have \\
{\rm (i)} $d_{BM} (C,D) \leq n^2$ for arbitrary $C, D \in {\cal C}^n$\\
{\rm (ii)}  $d_{BM} (M,N) \leq n$ for arbitrary $M, N \in {\cal M}^n$\\ {\rm (iii)} $d_{BM} (M, C) \leq n^{3/2}$ if $M\in {\cal M}^n$ and $C \in {\cal C}^n$. 
\end{cor}
We will address the sharpness of these bounds below and in subsequent sections. 
 
\vskip 0.09cm
Another -- slightly less immediate -- consequence of John's theorem is the fact that  $\mathfrak{C}^n$  
and $\mathfrak{M}^n$ are compact spaces with respect to $d_{BM}$.  Indeed, if $K$ is a convex body containing $0$, it's {\em gauge} $\|\cdot \|_K$ is defined by

\begin{equation}\label{def:gauge}
\|x\|_K := \inf\{ t\geq 0 : x \in tK\}.
\end{equation}
(If $K$ is $0$-symmetric, this is the norm for which $K$ is the unit ball.) In the language of gauges, 
John's theorem asserts that there exists $L$ affinely equivalent to $K$ such that 

\begin{equation}\label{John:gauge}
n^{-1}  |\cdot|  \leq\|\cdot\|_L\leq |\cdot| .
\end{equation}

In particular, the set of gauges verifying \eqref{John:gauge} is uniformly Lipschitz and uniformly bounded (say, on $B$), and it remains to appeal to the Arzel\`a-Ascoli theorem. 
This justifies referring to  $\mathfrak{C}^n$ or $\mathfrak{M}^n$ as the {\em Banach-Mazur compactum} or, more precisely, the  $n$th {\em Banach-Mazur compactum}.  
(Some authors employ the term {\em Minkowski compactum}, particularly for $\mathfrak{M}^n$.)  

Theorems \ref{john} and \ref{decomposition} focus on properties of the John ellipsoid and the John position, with estimates for the Banach-Mazur distance being an afterthought. However, an earlier version of Theorem \ref{decomposition}, which deals directly with the quantity $d_{BM} (K, B)$ and the corresponding ``distance position,'' can be found in \cite{Ader}. While that paper addresses the case $n=3$, it has been recently generalized to an arbitrary dimension in \cite{GruKob25} (the centrally symmetric case) and \cite{GruKob26} (the general case). We state here the result from \cite{GruKob25}. 

\begin{thm} \label{distance-position} 
Let $R\geq 1$ and let $K \subset E^n$ be a $0$-symmetric convex body such that $B\subset K \subset RB$. 
Then the following are equivalent. \\
$1^\circ$ $d_{BM} (K, B) = R$. \\
$2^\circ$ For every symmetric $n\times n$ matrix $A$, there exist an outer contact point $y \in \partial K \cap \partial (RB) $ and an inner contact point $z \in \partial K \cap \partial B $ such that $\langle y, Ay \rangle \geq R \langle z, Az \rangle$. \\
$3^\circ$ There exist integers $l,m\geq 1$, outer contact points $y_1,y_2, \ldots , y_l $, inner contact points $z_1, z_2,\ldots, z_m$, 
 and sequences of positive scalars $(\lambda_1, \lambda_2, \ldots , \lambda_l)$ and $(\mu_1, \mu_2, \ldots , \mu_m)$  such that
\begin{equation} \label{eq:GruKob}
\sum_{j=1}^l \lambda_j   y_j  y_j^t = \sum_{k=1}^m \mu_k   z_k z_k^t .   
\end{equation}
In the last condition, 
 there exists a choice with $l+m \leq \frac{n(n+1)}2 +1$.
\end{thm}
The characterization of the distance position in the not-necessarily-symmetric case is in the same vein, but more complicated. Additionally, \cite{GruKob26} contains {\em necessary} conditions for two arbitrary convex bodies to be in the distance position that are of similar nature. We believe that all these results will be very useful. At the very least, they should allow to produce conceptually simpler arguments for many known facts, whose original proofs were highly computational. 

To complete this section, let us emphasize the crucial influence of the circle of ideas related to the Banach-Mazur distance in general, and the John Theorems \ref{john} and \ref{decomposition} in particular, on various areas of computational geometry and convex analysis. This influence spans decades, from providing theoretical basis for the ellipsoid method in linear programming in 1970's to facilitating progress on frame theory and the Kadison-Singer problem in 2000's and 2010's.  Sample references include 
 \cite{YN76}, \cite{Le79}, \cite{Kha79}, \cite{Kha80},  \cite{NY83}, \cite{Pelc90}, \cite{Ver01}, \cite{LutYanZha05},  \cite{CEKP07}, \cite{Sri10}, \cite{SS12}, \cite{Youssef14}, \cite{Lasse15}, \cite{BG25}.

\vskip0.4cm

\section{The radii and the diameters of the Banach-Mazur compacta}  \label{sec:radius}

The results from the preceding section can be rephrased as simple and compelling statements about the (global) geometry of the Banach-Mazur compacta. 
First, the radius of  $\mathfrak{C}^n$ with respect to $d_{BM}$ and with $B$ as a center is precisely $n$.   (See Theorem \ref{thm:GLMP} in Section \ref{sec:one-symmetric} for a far reaching generalization.) 
Similarly, the radius of  $\mathfrak{M}^n$ with $B$ as a center is  $\sqrt n$. 
Finally, the diameter  of $\mathfrak{C}^n$ is upper-bounded by $n^2$ and that of $\mathfrak{M}^n$ by $n$. 
 
\vskip 0.09cm
This point of view leads to many natural and apparently difficult questions, and to a number of highly nontrivial results 
concerning the sharpness of these bounds.
We shall review the questions  -- and the results -- starting with the case $n=2$, and then proceed to the discussion of the general dimension. Later, in subsequent sections, we will explore other questions of similar flavor, typically under some natural constraints on one or both of the bodies under consideration. 
For example, if $M\in \mathfrak{M}^n$, then the radius of  $\mathfrak{C}^n$ with $M$ as a center is (as we already noticed) 
upper-bounded by $n^{3/2}$; we will see in Section \ref{sec:one-symmetric} that this bound if far from being sharp. We will also discuss the radii of $\mathfrak{C}^n$ and/or $\mathfrak{M}^n$ with specific centers such as the parallelotope $P$ or the cross-polytope (Section \ref{sec:parallelotope}) or the simplex  $T$ (Section \ref{sec:simplex}).

\subsection{The radii and the diameters of $\mathfrak{C}^2$ and $\mathfrak{M}^2$}  \label{radius:n=2}

The basic geometric parameters of $\mathfrak{M}^2$ are well understood. We have

\begin{thm}  [\cite{Aspl, [St 81]}]  \label{M2}
There exists $M_0 \in {\cal M}^2$ such that the radius of ${\cal M}^2$ with $M_0$ as a center is $\sqrt{3/2}$.  That is, $d_{BM}(M_0, M) \leq \sqrt{3/2}$ for all $M \in {\cal M}^2$ and, moreover, the equality holds if and only if $M$ is an affine-regular hexagon $H^{ex}$ or a parallelogram $P$.
\end{thm}
To clarify, the {\em affine-regular hexagon} is an affine image of the regular hexagon. Affine-regular pentagons and, more generally, affine-regular $n$-gons are defined analogously. (Of course, affine-regular quadrilaterals are just parallelograms, and all triangles are affinely equivalent.)  

The center $M_0$ from Theorem \ref{M2} was cleverly constructed by Stromquist in \cite{[St 81]} to be the body bounded by the straight lines $y= \pm 1$, by the circle $x^2+y^2=2$, and the ellipse $\frac{1}{2}x^2 +y^2 = \frac{4}{3}$. Such center is not unique as clearly the polar $M_0^\circ$ is just as good. Moreover, it is stated in \cite{[St 81]} that there are infinitely many small perturbations of $M_0$ which likewise work. 
 
\vskip 0.09cm
A formal consequence of Theorem \ref{M2} is the following. 
\begin{cor} \label{M2a} The diameter of $\mathfrak{M}^2$ is precisely $\frac{3}{2}$.   Moreover, the inequality $d_{BM} (M, N)\leq \frac{3}{2}$, valid for every $M, N \in {\cal M}^2$, is saturated if and only if  $\{M,N\}=\{H^{ex},P\}$. 
\end{cor}
The last equality is (of course) meant ``up to affine quivalence.'' 
The fact that $d_{BM} (H^{ex}, P) = \frac{3}{2}$ is stated in \cite{Aspl} and \cite{[St 81]}) without a proof. 
A detailed argument can be found in \cite{[L 21]} (see also \cite{Born,VK,AMP26}).  
 
\vskip 0.09cm
Theorem \ref{M2} and Corollary \ref{M2a} completely determine the values of  the radius and the diameter of  $\mathfrak{M}^2$ with the somewhat surprising conclusion that the Euclidean disc $B=B^2$ is {\em not} the optimal center. 

\vskip 0.09cm
The situation is far less clear for not-necessarily-symmetric bodies, but there are numerous partial results and reasonable conjectures, some of which will be elaborated on in subsequent sections.   
First, the ``trivial'' upper bound of $4$ on the diameter of $\mathfrak{C}^2$ (following from Theorem \ref{john}) has been improved to $3$ in \cite{[L 07]} and subsequently to $\frac{19-\sqrt{73}}4 \approx 2.614 $ in \cite{BPP18}.  
For a lower bound, we have $\rho := 1+ \frac{1}{2}\sqrt{5} = d_{BM} (\Tt,F) \approx 2.118$, where $\Tt$ is a triangle and $F$ is an affine-regular pentagon  (\cite{[L 92],FMRWY92}). 
It is apparently believed that this pair is extremal (for instance, see \cite{BPP18}), though we couldn't identify the original reference stating this question or conjecture.

\begin{prob} \label{prob:diamC2}
Is the diameter of $\mathfrak{C}^2$ equal to $d_{BM} (\Tt,F) = 1+ \frac{1}{2}\sqrt{5}$?
\end{prob}

Concerning the radius of $\mathfrak{C}^2$, it is shown in \cite{BPP18} that, for an appropriately chosen heptagon  $H^{ep}$, 
$\max \{d_{BM} (H^{ep},K) : K \in \mathfrak{C}^2\} \leq \frac{117}{70} < 1.672$.  
While there is no convincing evidence that this is the best choice, it is clear that, again, $B^2$ is not the optimal center. 
The best known lower bound for the radius seems to be $\rho^{1/2} \approx 1.455$.  
This suggests the next possible step, motivated by the strategy from \cite{[St 81]}: 

\begin{prob}
Identify bodies that are ``midpoints'' between the triangle $\Tt$ and the affine-regular pentagon $F$,  and try to calculate the radii of $\mathfrak{C}^2$ with those bodies as centers. 
\end{prob}
  
\subsection{The radii and the diameters of $\mathfrak{C}^n$ and $\mathfrak{M}^n$ for general $n$, and the Bachmann-Landau asymptotic notation}  \label{radius:ngeneral}

The results cited in Section \ref{radius:n=2} suggest that, for general $n\in \bN$, calculating the exact values of the diameters and/or identifying the optimal centers and the corresponding radii  will be very difficult, if not hopeless. 
Accordingly, one settles for the next best objective, determining the asymptotic order (as $n \to \infty$) of the quantities in question.  
In this context, it will be convenient to use the Bachmann-Landau notation, which is standard (particularly) in theoretical computer science.  

\begin{ddd}[Bachmann-Landau notation] Let  $f=(f_n)$ and $\phi=(\phi_n)$  be nonnegative sequences. We then say that \\
${\rm (i)}$ \ $f = O(\phi)$ if there exists a constant $c > 0$ such that $f_n\leq c\phi_n$  for all $n$\\ 
${\rm (ii)}$ \ $f = \Omega(\phi)$ if there exists a constant $c'> 0$ such that $f_n\geq c'\phi_n$  for all $n$\\ 
${\rm (iii)}$ \ $f = \Theta(\phi)$ if there exist constants $c,c' > 0$ such that $c' \phi_n \leq f_n\leq c   \phi_n$  for all $n$\\ 
${\rm (iv)}$ \ $f \sim\phi$ if $f_n/\phi_n \to 1$  as $n\to \infty$ \\
${\rm (v)}$ \ $f = o(\phi)$ if $f_n/\phi_n \to 0$  as $n\to \infty$.
\end{ddd}
Clearly, the above definition can be generalized to functions or quantities depending on variables (or parameters) other than the integers. 
However, in our discussion the parameter will be (most often) the dimension, so the setting of sequences will be generally sufficient. It should also be noted that while the Bachmann-Landau notation is usually described as ``asymptotic'', only the last two conventions are  
{\em truly} of that nature. The first three require that a specific bound holds for any fixed $n$ and are thus {\em de facto} ``nonasymptotic.''   Finally, sometimes it is more convenient to use an alternative notation for $f = o(\phi)$, namely $f \ll \phi$ or $\phi \gg f$.
 
\vskip 0.09cm
In this framework, the question of the diameter of $\mathfrak{M}^n$ was fully resolved in the seminal work of
Gluskin \cite{[Gl 81],[Gl 87]}, who proved that the upper bound of $n$, implied by John's theorem, is asymptotically optimal. 

\begin{thm} \label{gluskin} There exists a constant $c>0$ such that, for every $n$, there exist $M_1, M_2 \in {\cal M}^n$ such that $d_{BM} (M_1,M_2) \geq cn$. Consequently, $\diam \, \mathfrak{M}^n = \Theta(n)$. 
\end{thm}

As satisfying as Theorem \ref{gluskin} is, it still leaves out a few ``loose ends.''

\begin{prob}
    What is the diameter of $\mathfrak{M}^n$ for small $n$, say $n=3,4$?
\end{prob}

\begin{prob}
    Is the $\diam \, \mathfrak{M}^n \leq \theta n$, where $\theta \in (0,1)$ is a universal constant?
\end{prob}

Here is an outline of (a variant of) Gluskin's construction, on which Theorem \ref{gluskin} is based.  
Given $n\in \bN$ and $N\geq n$, consider $N$  points $x_1,x_2, \ldots , x_{N}$ chosen at random (uniformly and independently) from the sphere $S^{n-1}$,  and let $M$ be the convex hull of $\{\pm \, x_i : i=1,2,\ldots, N\}$. 
Then, with probability $1$, $M$ is a $0$-symmetric convex body in $E^n$. Next, let $M_1, M_2$ be two independent copies of $M$ and choose  $N=3n$. 
Then one can show that there is a constant $c_1>0$ (independent of $n$), such that, with probability close to $1$, we have $vM_1 \not\subset c_1\sqrt{n} M_2$ for every $v \in \gSL(n)$ (i.e., $v$ a linear operator with $\det v =1$). 
Since, by symmetry, the same is true about  $v^{-1}M_2 \not\subset c_1\sqrt{n} M_1$, it follows that for most choices of $M_1, M_2$ we have  $d_{BM}(M_1, M_2) > c_1^2 \hskip0.2mm n$, as needed.  

In turn, the failure -- with large probability -- of the inclusion $vM_1 \subset c_1\sqrt{n} M_2$  for {\em all} $v \in \gSL(n)$ is established in two steps. First, one estimates the probability of the inclusion holding for a {\em fixed} $v \in \gSL(n)$, and then one uses discretization and the union bound. Various tweaks are needed and some estimates require finesse, but the independence -- which is inherent in the scheme -- makes the calculation workable. 
 
\vskip 0.09cm
A remarkable feature of the above construction is that it {\em does not} allow to pinpoint any bodies, for which 
 $d_{BM} (M_1,M_2) \geq cn$.  Indeed, to the best of our knowledge, 
 the answer to the following is unknown. 
 
 \begin{prob} \label{explicit}
 Give an {\em explicit} sequence of (pairs of) bodies $K_n, L_n \in {\cal M}^n$,  for which the sequence  of ratios $d_{BM}(K_n, L_n)/\sqrt{n}$ is unbounded.  Alternatively, supply an efficient algorithm generating such pairs.  
 \end{prob}

 \vskip 0.09cm 
 We point out that there are several ways of generating ensembles of random convex bodies that appear in the literature. Modulo some simplifications,\footnote{In the original Gluskin's paper \cite{[Gl 81]} and in some other works that followed, the random body $M$ is the absolute convex hull of $\{e_1,\ldots, e_n, x_1,\ldots,x_N\}$, where $(e_i)$ is the standard basis of $\bR^n$ and $x_j$'s are selected independently at random.} the most common schemes can be outlined as follows (see \cite{LMOT}, \cite{Gi-Hm} for even more general schemes). 

\noindent (i) The Gluskin's construction sketched above, with the random body $M$ being the absolutely convex hull of $N$ points chosen independently from the sphere $S^{n-1}$.\\
(ii) A similar scheme, except that the $N$ points are drawn form the standard $n$-variate Gaussian distribution (cf. \cite{Sz83}). \\
(iii) The random body $M$ being the unit ball of the quotient space $\ell_1^N/F$, where $F$ is a random $n$-codimensional subspace of $\ell_1^N$ drawn according to the invariant measure on the corresponding Grassmann manifold ${\bf Gr}_{N,N-n}$(cf. \cite{FKP}).  

While it has been known to the experts that the three ensembles above are equivalent (or essentially equivalent; cf. \cite{ManSza,BarVit}), this interchangeability hasn't been well documented in the literature. Here is a sketch of the argument. 

First, the ensembles (ii) and (iii) are precisely equivalent in the sense that they induce the same distribution on $\mathfrak{M}^n$. Indeed, the random body $M$ generated by (ii) is the unit ball of the quotient of $\ell_1^N$, where the quotient map is an $n \times N$ matrix $A$ whose entries are independent standard Gaussian variables. Since the multivariate Gaussian measure is rotationally invariant, it follows that the distribution of the nullspace of $A$ is likewise invariant under rotations, which implies that it coincides with the invariant measure on the Grassmann manifold ${\bf Gr}_{N,N-n}$. (Note that the nullspace of $A$ can be of dimension greater than $N-n$, but that's an event of zero probability.) 

The ensemble (i) is not precisely equivalent to the other two, but it is ``nearly'' equivalent in the following sense. Let $y_1,y_2, \ldots , y_{N}$ be a sample from the $n$-variate Gaussian distribution and let $M$ be the corresponding convex body of type (ii). Set $x_j = y_j/|y_j|$, where $ |\cdot | $ is the Euclidean norm, and let $M_1$ be the convex hull of $\{\pm \, x_i : i=1,2,\ldots, N\}$. Then $x_j$'s are independent and -- again by the rotational invariance of the Gaussian measure -- uniformly distributed on $S^{n-1}$. In other words, $M_1$ is identical with the ensemble (i).  It remains to note that, with probability close to $1$, 
$d_{BM} (M,M_1)$ is small. This is a consequence of the well-known fact that, for each $j$, $|y_j|^2$ is a $\chi^2(n)$ random variable, which is strongly concentrated around the value $n$ (see, e.g., \cite{AuSz17}, Section 5.2.2). In other words -- at least if $N$ isn't too large -- $M_1$ is essentially $M/\sqrt{n}$. More precisely, if $\log N = o(n)$, then 
$d_{BM} (M,M_1)) \leq 1+O\left(\sqrt{\frac{\log N}{n}}\right)$. On the other hand, if  $\log N = \Omega(n)$, then -- with probability close to $1$ -- we have $d_{BM} (M,B^n)=O(1)$ and $d_{BM} (M_1,B^n)=O(1)$, hence also $d_{BM} (M,M_1)=O(1)$. 
 
\vskip 0.09cm 
By the triangle inequality, Theorem \ref{gluskin} implies that the radius of $\mathfrak{M}^n$ is at least 
$\sqrt{cn}$. 
Consequently, it is precisely of order $\sqrt{n}$ and the Euclidean ball $B = B^n$ is an asymptotically optimal center. However, we know from Theorem \ref{M2} that it is not {\sl precisely} optimal for $n=2$. A natural question to ask is whether this phenomenon can be reproduced for $n>2$.
\begin{prob} \label{rad-sym} For $n \geq 3$, is the radius of $\mathfrak{M}^n$ \\
(i) strictly smaller than $\sqrt{n}$ ? \\
(ii) smaller than $\theta\sqrt{n}$, where $\theta \in (0,1)$ is a universal constant?
\end{prob} 

Another natural question (that was asked, e.g., in \cite{Pelc84}) is whether $B = B^n$ is an asymptotically unique center. That is, if $K_n \in {\cal M}^n$ is such that the radius of $\mathfrak{M}^n$ with $K_n$ as a center is $O(\sqrt n)$, does it follow that the sequence $d_{BM} (K_n,B^n))$ is bounded? A resounding negative answer to this question was supplied in \cite{[BoSz88]} (see also \cite{Szar91}). We have

\begin{thm} \label{center-non-unique}
If $K_n$ is the Cartesian product of an $m$-dimensional parallelotope and $B^{n-m}$ with $m\leq n/2$, then the radius of $\mathfrak{M}^n$ with $K_n$ as a center is $O(\sqrt n)$, while $d_{BM} (K_n,B^n)) \geq \sqrt{m}$.
\end{thm}

The geometry of $\mathfrak{C}^n$ is less well understood, but it is known that the overall picture is quite different than in the symmetric setting. 
While a naive expectation could be that the Gluskin's construction that led to Theorem \ref{gluskin} can be modified to prove that the upper bound of $n^2$ on the diameter of $\mathfrak{C}^n$ that follows from John's Theorem \ref{john} can be asymptotically saturated, a simple argument shows that this is not possible. 

To see this, let us observe first that the constructions (i)-(iii) of the random set described earlier has to be fine-tuned to work well in the not-necessarily-symmetric setting. Indeed, if we simple-mindedly define $M$ be the convex hull of $N$ points chosen independently from the sphere $S^{n-1}$ (or according to the Gaussian measure) and if $N < 2n$, then  $0 \not\in M$ with probability greater than $\frac 12$ (see \cite{wendel}). To avoid this apparent anomaly we tweak the constructions a little, starting with the scheme (iii). 

Given $N>n$, let $T \subset \bR^{N-1}$ be an $(N-1)$-dimensional regular simplex with centroid at the origin.  
Next, let $E$ be a random $n$-dimensional subspace of $\bR^{N-1}$ drawn according to the invariant measure on the corresponding Grassmann manifold ${\bf Gr}_{N-1,n}$. Finally, let $M=P_E(T)$, where $P_E$ is the orthogonal projection on $E$. 
It can be shown that if $n$ is large and 
$N \geq \kappa n$ for some $\kappa >1$,  then -- with probability close to $1$ -- $d_{BM} (M,B^n) \leq C \sqrt{n}$, where $C$ is a constant depending only on $\kappa$. Hence, if $M_1, M_2$ are independent copies of $M$, then typically  $d_{BM} (M_1,M_2) \leq  C^2 n$. This is a far cry from the hoped for asymptotic saturation of the $O(n^2)$ upper bound on the diameter of $\mathfrak{C}^n$. 

For completeness, let us sketch the proof of the claim about random projections of the simplex that was stated above.  We will focus on the crucial case when $\kappa n \leq N \leq \lambda n$ for some fixed $\kappa, \lambda >1$. Note that this is the scenario that was at the heart of Gluskin's construction; the requirement that $\kappa >1$ assures that our example is ``nontrivial,'' while the upper bound on $N$ allows to simplify the arguments. 

We will exhibit $T$ in the following manner. 
First, let us identify $\bR^{N-1}$ with the hyperplane $H := \left\{a = (a_i)_{i=1}^N  : \sum_{i=1}^N a_i = 0\right\} \subset \bR^N$. Next, let $B_1^N$ be the unit ball in the $\ell_1$-norm and let $T_1$ be the positive face of $B_1^N$, i.e., $T_1 = \left\{a : \sum_{i=1}^N a_i = 1, a_i \geq 0\right\}$. Finally, we let $T=P_{H}T_1$, the orthogonal projection of $T_1$ on $H$. 

Now let $K=\conv(-T \cup T)= P_{H}B_1^N$. Clearly, the inradius of $K$ is at least $N^{-1/2}$ (because $B_1^N \supset N^{-1/2} B$). At the same time,  $\vol(K)\leq \vol(B_1^{N-1})$ (see \cite{Ball95}) and so its volume radius is $O(N^{-1/2})$ (see \cite{Sz78,SzTom}). (The volume radius of a $d$-dimensional set $T$ is defined as $\big(\vol(T)/\vol(B)\big)^{1/d}$, i.e., the radius of a Euclidean ball whose volume is the same as that of $T$. 
Actually, both the inradius and the volume radius of $K$ can be {\em exactly} evaluated, but such precision is not needed for our purposes.) 
Accordingly (\cite{Sz78,Sz80,SzTom}, cf. Theorem \ref{vr-kashin} in Section \ref{sec:sections}), for any $\theta \in (0,1)$ there is a constant $C(\theta)$ such that if $k \leq \theta N$, then $K \cap E \subset C(\theta)N^{-1/2} B$ \, for most of $k$-dimensional subspaces of $E \subset H$. {\em A fortiori}, for such $E$, $T \cap E \subset C(\theta)N^{-1/2} B$, and this conclusion holds for $k=n$ if we choose $\theta = \kappa^{-1}$.  
By duality (inside the $(N-1)$-dimensional space $H$, with $(\frac 1N,\ldots, \frac 1N)$ as the center), $P_E (T^\circ) = (T\cap E)^\circ \supset C(\kappa^{-1})^{-1}N^{1/2} B$. It remains to notice that $T^\circ = -\sqrt{N(N-1)} T$ to obtain $M=P_ET \supset C(\kappa^{-1})^{-1}(N-1)^{-1/2}B \supset C(\kappa^{-1})^{-1}\lambda^{-1}n^{-1/2}B$, which shows the claim (again, under the additional assumption that $N\leq \lambda n$). 

The schemes (i) and (ii) can be similarly modified. For example, one can consider $M=AT$, where $A$ is an $n \times (N-1)$ Gaussian matrix with independent $N(0,1)$ entries. Alternatively, one can select the vertices $x_1,x_2,\ldots,x_N$ of $M$ independently (from $S^{n-1}$, or according to the Gaussian measure) and then shift it so that the average $(x_1+x_2+\ldots +x_N)/N$ is at the origin. 

 While this observation didn't fully rule out examples with 
$d_{BM} (M_1,M_2) \gg  n$ and $N-n = o(n)$ (cf. Problem \ref{diam-nonsym} below), it strongly suggested that the diameter of $\mathfrak{C}^n$ was much smaller than $n^2$. 
The above heuristics have been subsequently confirmed by rigorous upper bounds. 
Rudelson \cite{[Rud 00]} proved that the diameter in question is $O\big(n^{4/3}\log^9(1+ n)\big)$.  
(A weaker but non-trivial  $O(n^{3/2})$ bound follows from the results of \cite{BLPS99}.) Finally,  very recent eprints \cite{BK,Bizeul} by Bizeul and Klartag and \cite{Bizeul} by Bizeul show the following ``nearly optimal'' result.  

\begin{thm} \label{bizeul-klartag}
The diameter of \, $\mathfrak{C}^n$ is $O\big(n\log (1+n)\big)$.  
\end{thm}
Without going into details (which are rather involved and based on deep recent breakthroughs in high-dimensional convex geometry, the area often referred to as {\em asymptotic geometric analysis}), the proof depends on showing that, given $K,L \in {\mathcal C}^n$, there exist affine images $K'=a(K), L'=b(L)$ such that, for most (in the sense of the corresponding Haar measure) rotations $V, W \in \gO(n)$,  
\begin{equation} \label{eq:bizeul-klartag}
 V(K') \subset \beta L' \ \hbox{ and } \ W(K') \subset \beta L',   
\end{equation}
where $\beta = O\big(\sqrt{n \log(1+n)}\big)$. The conclusion then follows by selecting any $(V,W)$ satisfying \eqref{eq:bizeul-klartag} and $W=V^{-1}$. 

While the upper bound supplied by Theorem \ref{bizeul-klartag} matches the lower bound of $n$ up to a polylogarithmic factor, the {\em precise} asymptotic behavior of the diameter of $\mathfrak{C}^n$ remains an open problem.  
In particular, one could ask the following ``intermediate'' questions.

\begin{prob} \label{diam-nonsym}
For $n\geq 3$, do there exist $K_n,L_n \in \mathfrak{C}^n$ with \\(i) $d_{BM} (K_n,L_n) >  n$? \\
(ii) $d_{BM} (K_n,L_n)/n \to \infty$? 
\end{prob}
Problem \ref{diam-nonsym}(i) asks in essence whether the example of the pair triangle-pentagon -- for which the distance is {\sl strictly} larger than the dimension -- can be generalized to an arbitrary $n\geq 2$. Regarding part (ii), there is no clear path to either proving or disproving the statement (but see the comments in the paragraph preceding Theorem \ref{bizeul-klartag}). It should be noted that the arguments from \cite{BK, Bizeul,[Rud 00],BLPS99} -- ultimately going back to \cite{BeGo} -- are based on the auxiliary bounds (involving the so-called $MM^*$-estimate, or $K$-convexity constant) that do {\em necessarily}  involve polylogarithmic factors, see \cite{Bou84,BLPS99}.  

In contrast, the question about the magnitude of the radius of $\mathfrak{C}^n$ is wide open, with the upper and lower bounds -- for general $n$  -- being $n$ and $\sqrt{n}$. (See Section \ref{radius:n=2} for the discussion of the case $n=2$.) Answers to the following ``intermediate'' questions would be of some interest. 

\begin{prob} \label{rad-nonsym} Is the radius of $\mathfrak{C}^n$ \\
(i) strictly smaller than $n$ for $n \geq 3$? Is it $\ll n$ (as $n\to \infty$)?\\
(ii) strictly larger than $\sqrt{n}$ for $n \geq 3$? Is it $\gg \sqrt{n}$ (as $n\to \infty$)?
\end{prob}

\vskip0.4cm

\section{Comparing convex bodies with a simplex}   \label{sec:simplex}

In the preceding section we discussed the problem of approximating $K \in  {\cal C}^n$ by (an affine image of) the Euclidean ball $B^n$, the $n$-dimensional body with most symmetries. 
Here we will go in the opposite direction and consider approximating by a simplex $\Tt=\Tt^n$. 
As we shall see, the scenery in the latter case is much more blurred than in the former, with many fundamental questions lacking a definite answer. 
 
\vskip 0.09cm
Before proceeding, let us comment here that the notion of {\em symmetry} is used in convex geometry and geometric functional analysis in several different ways. 
First, it may refer to the {\em size} of the group of affine maps that preserve $K$, the largest such group being isomorphic to $\gO(n)$ (in which case $K$ is an ellipsoid) and the smallest the trivial group consisting just of the identity map $I_n$.  
Next, there is the concept of {\em symmetric  normed spaces} (see Section \ref{sec:classes}). 
Finally, one may just refer to $K$ being centrally symmetric (or $0$-symmetric); this leads to the classical {\em quantitative} invariant, which was introduced by Minkowski in \cite{Mink} and can be defined as follows. 
\begin{ddd}
Let $C\in {\cal C}^n$. Given  $x\in {\rm int}(C)$ and a lines $L$ containing $x$,  let $r(L,x) \geq 1$ be the ratio in which $x$ divides the segment  $L\cap C$. We now define the {\em Minkowski measure of symmetry} of $C$ as 
\begin{equation} \label{eq:measure-dual}
{\rm as} (C) := \inf_{x\in {\rm int}(C)}\sup_{L\ni x} \; r(L,x).
\end{equation} 
\end{ddd} 
We note that both $\inf$ and $\sup$ in \eqref{eq:measure-dual}  are attained. This is immediate (by compactness) for the $\sup$. Concerning the $\inf$, it is enough to observe that $\phi(x)=\sup_{L\ni x} \; r(L,x)$ tends to infinity when $x$ approaches $\partial C$. 

An alternative ``dual'' definition (showing equivalence is a simple exercise in convex analysis) involves hyperplanes $H$ containing $x$. If $H_1, H_2$ are the two support hyperplanes of $C$ parallel to $H$, denote by $r'(H,x) \geq 1$ the ratio in which $H$ divides the distance between $H_1$ and $H_2$. 
We then have 
\begin{equation} \label{eq:measure}
{\rm as} (C)= \inf_{x\in {\rm int}(C)}\sup_{H\ni x} \; r'(H,x)
\end{equation}
Another convenient reformulation of \eqref{eq:measure-dual} is as follows (cf. the comments in the paragraph following Definition \ref{def:dBM}) 
 \begin{equation} \label{eq:measure-translates}
{\rm as} (C) = \inf \, \{\lambda >0 : -C \hbox{ is contained in a translate of } \lambda C\}.
\end{equation} 
The invariant ${\rm as} (C)$ and other measures of symmetry are discussed at length in Gr\"unbaum's foundational paper \cite{[Grun63]} and in the monograph \cite{To15} by Toth.  
See also \cite{GluLit} for symmetry properties of polytopes with ``few'' vertices. 

\vskip 0.09cm
The functional $C \to {\rm as} (C)$ achieves its minimum value $1$ if and only if $C$ is a centrally symmetric convex body; this is immediate from \eqref{eq:measure-dual}.  
Next, it is well known that ${\rm as} (C)$ attains its maximum value of $n$ if and only if $C$ is a simplex, which is thus -- in this sense -- the most asymmetrical body. 
A stable version of the last statement is also true, with better and better estimates obtained successively in  \cite{[Gu 05]}, \cite{Sch09},  \cite{BY26}, and\cite{BGG}. We have 
\begin{thm}[\cite{BGG}] \label{thm:BGG} 
If $0\leq \ve < 1$, $C \in  {\cal C}^n$, and ${\rm as}(C)\geq n-\ve$,  then 
$$d_{BM} (C,\Tt)\le \frac{1-\ve^2/2}{1-\ve}.$$ 
\end{thm}

A related  result that characterizes the simplex as the ``most asymmetrical body'' is the Rogers-Shephard inequality \cite{RoSh1,RoSh2}. 
\begin{thm} \label{RS} 
If $K\in \mathcal{C}^n$, then 
\begin{equation} \label{eq:RS}
 2^n \vol(K) \leq    \vol(K-K) \leq \binom{2n}{n}\vol(K). 
\end{equation}
 The first inequality is saturated iff $K\in \mathcal{M}^n$ and the second one iff $K$ is a simplex.    
\end{thm}
Similarly as in Theorem \ref{thm:BGG}, we have a stable version of \eqref{eq:RS} (\cite{Bor05,BGG}) : {\em There is a constant $\gamma =\gamma(n) >0$ such that if $0< \ve < \gamma^{-1}$ and $\vol(K-K) \geq (1- \ve) \binom{2n}{n}\vol(K)$, 
then $d_{BM}(K, T) \leq 1 + \gamma \ve$.} 
  
\vskip 0.09cm
An important (and obvious) feature of the Minkowski measure of symmetry is its affine invariance, i.e., ${\rm as} \big(a(C)\big) =  {\rm as} (C)$  if $a$ is an affine map.   
Moreover, we have the following. 

\begin{pro} \label{lipschitz}  If $C, D \in  {\cal C}^n$, then 
\begin{equation} \label{eq:lipschitz}
{\rm as} (C) \leq {\rm as} (D) \, d_{BM} (C, D). 
\end{equation}
\end{pro} 
That is, in the multiplicative sense, the functional $C \to {\rm as} (C)$ is $1$-Lipschitz with respect to  $d_{BM}$.  
Functionals having property \eqref{eq:lipschitz} are useful for establishing {\em lower} bounds for $d_{BM}$; for example, it follows immediately from what we said above that 
 \begin{equation} \label{eq:measure-simplex}
 d_{BM} (M, \Tt) \geq n  \hbox{ if } M \in {\cal M}^n .  
\end{equation}

Another consequence of \eqref{eq:lipschitz} is one more   characterization of the  Minkowski measure of symmetry
\begin{equation} 
{\rm as} (C) = \min_{\in {\cal M}^n}d_{BM} (M, C).
\end{equation}
Indeed, the inequality $``\leq \hskip0.2mm "$ follows immediately from \eqref{eq:lipschitz}, 
and it is easy to see that it is saturated for $M=C \cap (x-C)$, where $x\in {\rm int}(C)$ is a point at which the infimum in \eqref{eq:measure-dual} is attained. 

While we believe that Proposition \ref{lipschitz} was folklore since the ideas of Gr\"unbaum's paper \cite{[Grun63]} were assimilated by the convex geometry community, reasonably explicit statements that we could locate in the literature were much more recent (see Brandenburg and K\"onig \cite{BK15} and Huang \cite{Hua18}). Accordingly, it seems useful to supply a sketch of the argument, which is based on the characterization \eqref{eq:measure-translates}. 

Denote  $\lambda = {\rm as} (D)$ and  $\delta =d_{BM} (C, D)$, and let $C'=a(C)$  be such that $D\subset C'$, while  $C'$ is contained in a translate of $\delta D$ and, consequently, $-C'$ is contained in a translate of $\delta (-D)$. Since $-D$ is contained in a translate of $\lambda D \subset \lambda C'$, it follows that $-C'$ is contained in a translate of $\delta \lambda C'$.  
In other words, ${\rm as} (C)= {\rm as} (C') \leq \delta \lambda = d_{BM} (C, D){\rm as} (D)$, as needed. 

In view of \eqref{eq:measure-simplex} (and Corollary \ref{cor:dBMmod-sym} below), it would be tempting to conjecture that, for all $C \in  {\cal C}^n$,  $d_{BM} (C, \Tt) \leq n$ 
with equality if and only if $C \in {\cal M}^n$ (the if part was shown by Gr\"unbaum in \cite{[Grun63]} on p. 259, see also our Section 5). 
This is, however, false already in the planar case: if $F$ is a regular pentagon, then $d_{BM} (F,\Tt)= 1+ \frac{1}{2}\sqrt{5}=: \rho  \approx 2.118$, see \cite{[L 92],FMRWY92}. 
 Both 
 these papers conjecture that this is the extremal configuration. We thus have (cf. Problem \ref{prob:diamC2})
 \begin{prob} Is \ $\max \{d_{BM} (C,\Tt) : C \in {\cal C}^2\} = \rho$?  \end{prob} 
 What seems to be rigorously established (by Fleischer, Mehlhorn, Rote, Welzl and Yap in \cite{FMRWY92}) is only that $d_{BM} (C,\Tt) \leq 9/4$ for every $C \in {\cal C}^2$.   
(A very simple and short proof in \cite{[L 92]} gives a weaker bound $d_{BM}(C,\Tt) \leq 5/2$.) 
Let us also mention that it follows by continuity that, for some $\theta\in (0,1)$, 
the non-centrally symmetric body $K= (1-\theta)T + \theta F$ verifies $d_{BM} (K,\Tt) =2$.
 
\vskip 0.09cm 
However, the naive guess stated above is not far from the truth. 
Namely,  \cite{[L 92a]} establishes the following (see \cite{[L 11]} and \cite{[Gu 12]} for alternate proofs).
\begin{thm} \label{d(C,T)}
If $n \geq 2$, then  
\begin{equation}\label{n+2} 
d_{BM}(C,\Tt) \leq n+2 
\end{equation}
for every $C\in {\cal C}^n$. 
\end{thm} 

\begin{figure}[htbp]

\includegraphics[width=0.38\textwidth]{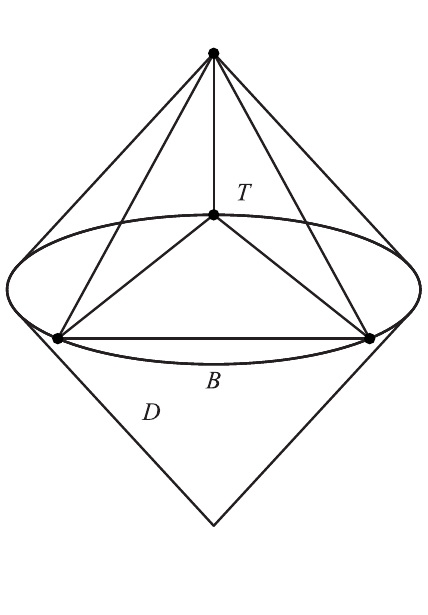} 
\vskip-5mm
\caption{A simplex $\Tt$ of the largest volume in a double cone $D \subset E^3$   \label{fig:double-cone} } 

\end{figure}

Similarly as is in the case of John's theorem, Theorem \ref{d(C,T)} can be proved by considering any simplex $\Tt_{max}$ of maximal volume contained in $C$. 
While -- as opposed to ellipsoid -- such simplex is in general far from being unique, a simple variational argument shows  (see \cite{[L 92a], [L 11]}) that the homothetic image of $\Tt_{max}$ with ratio $n+2$ and the homothety center in the centroid of $\Tt_{max}$ contains $C$.  
 With these additional constraints, the estimate \eqref{n+2} cannot be improved, which follows from the example of a double cone $D$ in $E^n$ whose base is an $(n-1)$-dimensional ball $B$. 
The reason for this is that a simplex $\Tt$ contained in $D$ has the maximum possible volume if  (and only if, when $n\geq 3$) one vertex of $\Tt$ is at an apex of $D$ and the remaining vertices of $\Tt$ form a regular $(n-1)$-dimensional simplex inscribed in $B$. 
This was shown in \cite{[L 11]} and is illustrated in Figure \ref{fig:double-cone} for $n=3$.

\vskip 0.3cm 

Other interesting properties of the simplex $\Tt_{max}$ of maximal volume contained in the body $C$ is that $\Tt_{max}$ is inscribed in $C$ (i.e., every vertex of $\Tt_{max}$ belongs to $\partial C$), and for every vertex of $\Tt_{max}$ the hyperplane containing this vertex and parallel to the opposite face  supports $C$ (see \cite{[Grun63]}, p. 259, or the Proposition of \cite {[L 06]} and Lemma 3 of \cite {[L 10]} for details). 
Equivalently, the homothetic image of $\Tt_{max}$ with {\sl negative} ratio $-n$ and the homothety center in the centroid of $\Tt_{max}$ contains $C$. 
This motivates a modification of the definition \eqref{eq:dBM} of the Banach-Mazur distance by allowing homotheties with negative ratios.
\begin{ddd} [Gr\"unbaum distance] \label{grunbaum}
    Given $C, D \in \mathcal{C}^n$, denote
\begin{equation} \label{def:dBMmod}
d_{BM}^{\pm} (C,D) := \inf_{a, h_{\lambda}} \{ |\lambda| : \, a(D) \subset C \subset h_{\lambda}(a(D)) \}. 
\end{equation} 
\end{ddd}
%
In this notation, the preceding observation can be restated as    
\begin{thm} \label{thm:dBMmod-bound} If $C \in  {\cal C}^n$, then 
\begin{equation} \label{eq:dBMmod-bound}
d_{BM}^{\pm}  (C, \Tt) \leq n . 
\end{equation}
\end{thm}
The same argument shows that $d_{BM}  (M, \Tt) \leq n$ if $M \in {\cal M}^n$, and combining this with \eqref{eq:measure-simplex} we can deduce 
\begin{cor} \label{cor:dBMmod-sym} If $M \in {\cal M}^n$, then
\begin{equation} \label{eq:dBMmod-sym}
d_{BM}^{\pm} (M,\Tt) = d_{BM}(M,\Tt)  = n.  
\end{equation}
\end{cor}
We will see in Section \ref{sec:one-symmetric} that \eqref{eq:dBMmod-bound} and \eqref{eq:dBMmod-sym} are not specific to the simplex; $T$ can be replaced by any convex body. 

The quantity $d_{BM}^{\pm} (\cdot,\cdot)$ was introduced and studied in \cite{[Grun63]}, where it was denoted by $\delta_2(\cdot,\cdot)$, and where
 the relations \eqref{eq:dBMmod-bound} and \eqref{eq:dBMmod-sym} were shown there on p. 259 by the same method as above.   
Presumably for this reason, it is often referred to as the {\em Gr\"unbaum distance.}

\vskip 0.09cm
While \eqref{n+2}, \eqref{eq:dBMmod-bound}  and \eqref{eq:dBMmod-sym} are quite elegant, the concepts that are behind them do not tell the entire story,  even in the planar case. 
To illustrate this, we will discuss in detail the pair  $\Tt, Q$, where $Q$ is a square; this will also reveal various intricacies of the not-necessarily-symmetric setting and suggest some questions that may be worthwhile asking, also in more general settings.  
 
\begin{figure}[htbp]

 \includegraphics[width=0.85\textwidth]{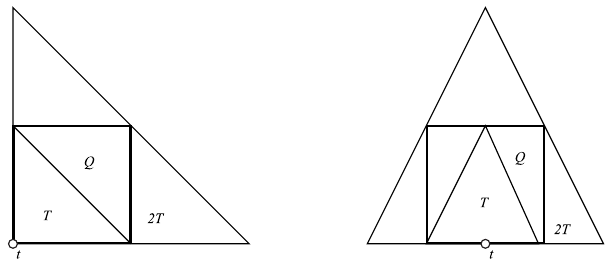} \\

\caption{Two positions of a triangle $\Tt$ of the largest area in a square $Q \subset E^2$ \label{fig:la25100a}.} 

\end{figure}

\vskip 0.09cm
We first consider two triangles $T=\Tt_{max}$ of largest area contained in $Q$, see Figure \ref{fig:la25100a}.   
The centers of the homotheties $h$ which transform $T$ into $h(T)\supset Q$ are denoted -- in both cases -- by $t$ and, to reduce clutter, we write $2T$ instead of $h(T)$ or $h_2(T)$. 
(The first configuration is a $2$-dimensional version of the example shown in Figure \ref{fig:double-cone}.)  
Moreover, we can consider a family of positions of $T$ whose ``top" vertex moves along the ``top" side of $Q$ between the two 
positions; a simple geometric argument shows that -- up to symmetries -- these are all possible positions of $\Tt_{max}$. In all such cases, the ratio of the homothety $h$ required to assure that $h(T) \supset Q$ is $2$. 

 \vskip 0.25cm 

It is a little less obvious that the same family of configurations yields -- again, modulo symmetries -- all possible positions  certifying the equality $d_{BM}(Q,\Tt)=2$. This must be known, but since we couldn't find a convenient reference, we sketch the argument, which depends on ``reverse engineering'' the  proof of \eqref{eq:lipschitz}. (See also Proposition \ref{prop:dTM} and Problem \ref{prob:dTP} in Section \ref{non-sym}.)  

Consider $T \subset Q \subset 2T$ and let $T'$ be the symmetric image of $T$ with respect to the center of $Q$. Then $T' \subset Q \subset 2T$ and $T'$ is a translate of $-T$. Since the only way a translate of $-T$ can fit inside $2T$ is when its vertices are midpoints of the sides of $2T$, it follows that all these midpoints are boundary points of $Q$. If all of them are vertices of $Q$, we have the first configuration in Figure \ref{fig:la25100a}. If one of them is an interior point of a side of $Q$, the other two must be vertices of $Q$, and we have a configuration that is intermediate between the two pictured in Figure \ref{fig:la25100a}. 

A curious consequence of the above considerations is the following fact. 
\begin{pro} \label{boundary-only} 
    In all the configurations $T\subset P \subset h(T)$ certifying the equality $d_{BM}(P,T)=2$ according to Definition \ref{def:dBM}, the center of the homothety $h$ (with ratio $2$) is a boundary point of $T$ and of $P$.
\end{pro}

The feature pointed out in Proposition \ref{boundary-only} may be a drawback for some applications. For example, while Banach-Mazur distance estimates can be conveniently restated in the language of gauges (cf. \eqref{def:gauge}, \eqref{John:gauge}) and hence encode Lipschitz properties of the functionals in question, the gauge of $K-t$ takes infinite values if $t \in \partial K$; this makes quantitative analysis not quite feasible. Accordingly, it makes sense to consider constraints that ensure that the center of the homothety is, in some appropriate sense, a ``central point'' of the bodies that are being compared.

\begin{figure}[htbp]        
\includegraphics[width=0.5\textwidth]{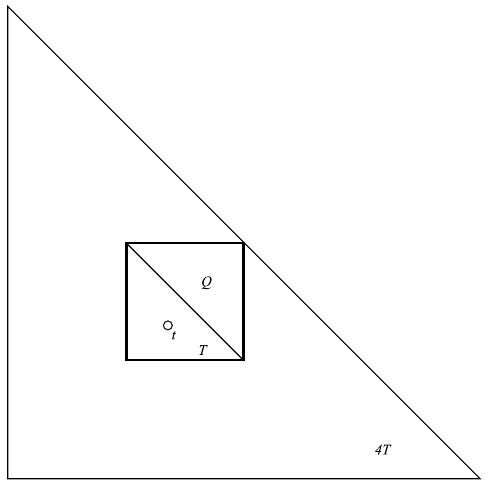} 
 \caption{If we require that the center of the homothety coincides with the centroid of  $\Tt_{max}$, then the smallest ratio $\lambda$ such that $h_\lambda(\Tt_{max})\supset Q$ equals to $4$. 
 This is a $2$-dimensional version of the example shown in Figure \ref{fig:double-cone}.\label{fig:la25100c}} 
\end{figure}

\vskip0.2cm
For the pair $(\Tt, Q)$ (or $(T,P)$), one possible approach is to allow exclusively homotheties whose center $t$ is the centroid of $\Tt=\Tt_{max}$. In that case, the smallest ratio required to enlarge $\Tt$ so that the image contains $Q$ is $4$, see Figure \ref{fig:la25100c}.  
Again, this illustrates the fact that, under this restriction, the bound \eqref{n+2} can be saturated (cf. the comments following Theorem \ref{d(C,T)}).

\begin{figure}[htbp]        
\includegraphics[width=0.95\textwidth]{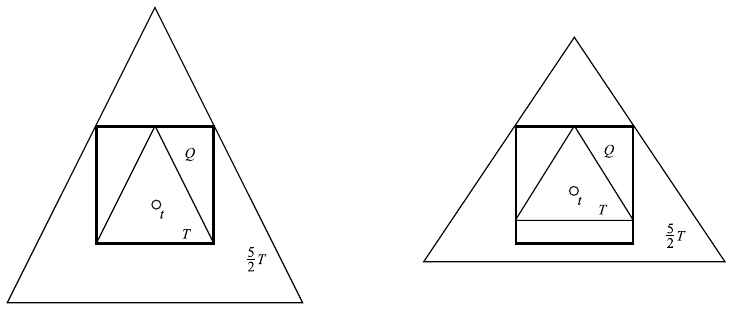} 
\caption{Two configurations of $\Tt \subset Q \subset h_{5/2}(T)$, where the center of the homothety is the centroid  of $T$.    \label{fig:la25100b}} 
\end{figure}

\vskip0.3cm
If we {\em do not} insist that $\Tt=\Tt_{max}$ -- while still requiring that the centroid of $\Tt$ is the center of homothety $h_\lambda$ such that $T \subset Q \subset h_\lambda(\Tt)$ -- an improvement is possible. 
Namely, we can have $\lambda =5/2$, with two configurations that realize this value  illustrated in Figure \ref{fig:la25100b}.  Note that in the second configuration the centroid of $\Tt$ coincides with the center of $Q$. 
The value $5/2$ is optimal as stated at the bottom of p. 259 of \cite{[Grun63]} (see \cite{[L 24a]} for a detailed argument). 
Moreover, this is not a special property of the square; we have the following fact, which was conjectured in  \cite{[Grun63]} and confirmed in \cite{[L 92]}. 
\begin{thm} \label{thm:T-center} If $M \in {\cal M}^2$ with center at $o$, then there exists a triangle $\Tt$ with centroid at $o$ such that $\Tt \subset M \subset \tfrac 52(\Tt)$. 
\end{thm}

Finally, let us note that the two configurations in Figure \ref{fig:la25100a} -- when looking at $Q$ and $2T$ -- represent parallelograms of maximal volume contained in a triangle. Together with the intermediate configurations, and up to affine equivalence, these are all positions of such parallelograms. It makes sense to look for homotheties $h_\lambda$ with center in the center of $Q$, for which $h_\lambda(Q) \supset 2T$, and it is easy to see that the smallest value of $\lambda$ that works is $3$. (A far-reaching generalization of this observation is given below in Theorem \ref{maxvol}.)

More systematic attempts directed at modifying the definition of $d_{BM}$ by imposing various ``centrality'' constraints on the center of the homothety $h_\lambda$ are considered in Sections  \ref{sec:one-symmetric} and \ref{sectioncentroid}.

\vskip 0.09cm
The parallelogram-triangle pair appears also in the following result  from \cite{[L 14]}.

\begin{thm} \label{quadrilaterals}
The Banach-Mazur distance between any two convex quadrangles is at most $2$ with equality only for the parallelogram and the triangle.
\end{thm}
Triangle is treated here as a degenerate quadrilateral.
 
\vskip 0.09cm
A survey of results on the approximation of planar convex bodies by triangles is presented in \cite{[BL]}, including 
in particular those concerning the BM-distance.  
The authors put forward the following conjecture, which earlier appeared in \cite{[L 92]}.

\begin{prob}
Given $C \in {\cal C}^2$,  can one inscribe a triangle $\Tt$ in $C$ such that $C\subset h_\lambda(\Tt)$, where $\lambda = 1+\frac{3}{5}\sqrt 5$ and the center of the homothety $h_\lambda$ coincides with the centroid of $\Tt$?  
\end{prob}
A better estimate can not hold even for centrally symmetric bodies as follows from the example in Figure 2 of \cite{[L 92]}. 
Note that it is {\em not} assumed here that the these centroids coincide with the centroid of $C$, cf. Section  \ref{sectioncentroid}.

Applying similar ideas and evaluations in three-dimensional space, Novotny \cite {[Nov]} proved that for every $C\in {\cal C}^3$ and for the tetrahedron 
$\Tt$ we have $d_{BM} (C,\Tt) \leq \frac{13}{3}$. 
\begin{prob}
Improve this estimate or show that it is optimal. 
More generally, improve the estimate \eqref{n+2} for bodies  $C \in {\cal C}^n$ for general $n\in \bN$. 
\end{prob} 

\vskip0.4cm

\section{Comparison of convex bodies, one of which is centrally symmetric }  \label{sec:one-symmetric}

As illustrated by our analysis of the pair   $T, Q \in {\cal C}^2$  in Section \ref{sec:simplex}, in some scenarios the definition  \eqref{eq:dBM} of the BM-distance has some drawbacks. This is in particular the case when our underlying setting is an $n$-dimensional normed space $X$ with the canonical $0$-symmetric convex set $M=B_X \in {\cal M}^n$, which we want to compare with another not-necessarily-symmetric convex set $C \in {\cal C}^n$, and to restate the resulting Banach-Mazur distance estimates in the language of gauges (cf. \eqref{def:gauge}, \eqref{John:gauge}).   As pointed out in  Section \ref{sec:simplex}, this objective is rather incompatible with the original definition of the BM-distance. First, the resulting gauge of $M$ is typically not necessarily the norm of $X$. Next, if -- as it can happen -- the center of the  homothety witnessing $d_{BM}(M,C)$ is a boundary point of one of the bodies, the corresponding gauge  takes infinite values, which makes any quantitative analysis difficult or simply infeasible.   

\begin{figure}[htbp]        
\includegraphics[width=0.7\textwidth]{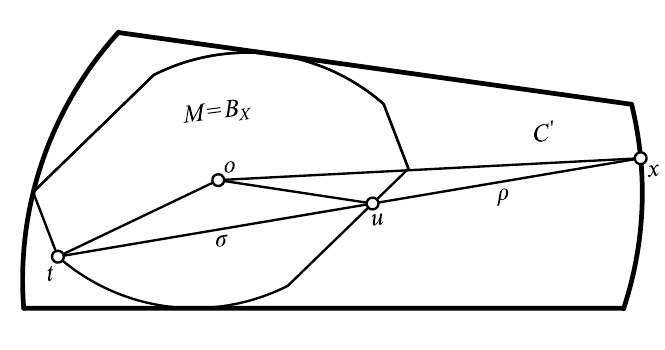} 
\vskip5mm
\caption{The illustration of the proof of Theorem \ref{thm:ctr}.    \label{2lambda-1}} 
\end{figure}

To alleviate these shortcomings, we propose the following ``centered'' variant of the Banach-Mazur distance.  
\begin{ddd} \label{def:center}
 Given $C \in \mathcal C^n$ and $M \in \mathcal M^n$, let us define 
\begin{equation} \label{def:dBM-center}
d_{BM}^{\rm ctr} (C,M) = \inf_{a, h_{\lambda}} \{ \lambda : \, M \subset a(C) \subset h_{\lambda}(M) \}, 
\end{equation}
where $a$ is an affine transformation and $h_{\lambda}$ is a homothety with center in the center of $M$ and a {\sl positive}  ratio $\lambda$.    
\end{ddd}
 If $M$ is $0$-symmetric, which we can and will assume through the remainder of this section, the condition in \eqref{def:dBM-center} simplifies to $M \subset a(C) \subset \lambda M$. It is worthwhile to point out that while in general $d_{BM}^{\rm ctr} \geq d_{BM}$, some of the facts concerning $d_{BM}$ stated in prior sections were in reality results about $d_{BM}^{\rm ctr}$. For example, the first assertion of John's Theorem \ref{john} can be rephrased as
 \begin{equation} \label{john-ctr}
     d_{BM}^{\rm ctr}(K,B) \leq n \quad \hbox{ for }  K \in \mathcal{C}^n.
 \end{equation}
Some other facts (for example, the bound \eqref{Gru-n2} in Theorem \ref{thm:GLMP} further in this section) are specific to $d_{BM}$ (or $d_{BM}^\pm$) and simply do not hold for $d_{BM}^{\rm ctr}$. 
Still, it is comforting to know that the two invariants, $d_{BM}$ and $d_{BM}^{\rm ctr}$ are not significantly different. We have 
\begin{thm} \label{thm:ctr} 
If $C \in \mathcal C^n$ and $M \in \mathcal M^n$, then 
\begin{equation} \label{eq:ctr} 
d_{BM} (C,M) \leq d_{BM}^{\rm ctr} (C,M) \leq  2d_{BM} (C,M) -1. 
\end{equation}
\end{thm}
The first inequality in \eqref{eq:ctr} is trivial. For the second, assume -- as we may -- that $M=B_X$, the unit ball for some norm $\|\cdot \|$ on $\bR^n$. Denote $\lambda = d_{BM} (C,M)$, and let $h_\lambda$ be a homothety with ratio $\lambda$ and center $t$ such that $M \subset C' \subset h_\lambda(M)$, where $C'=a(C)$, see Figure \ref{2lambda-1}. Recall that necessarily $t\in M$, i.e., $\|t\| \leq 1$. 

Let $x\in C$, our objective is to show that $\|x\|\leq 2\lambda-1$. 
If $\|x\|\leq 1$, there is nothing to show, so assume $\|x\| > 1$.  Consider the segment $[t,x]$; its intersection with $M$ is a closed segment $[t,u]$, where $u \in \partial M$, i.e., $\|u\|=1$. 
Denote $\rho = \|x-u\|$ and $\sigma=\|u-t\|$ (observe that $\sigma \leq 2$). Since $x \in h_\lambda(M)$, it follows that 
$\frac{\rho+\sigma}\sigma \leq \lambda$, which translates to $\rho \leq (\lambda-1)\sigma\leq 2(\lambda-1)$. We conclude by noting that $\|x\| \leq \|x-u\|+ \|u\| = \rho+1 \leq 2(\lambda-1)+1 =2\lambda-1$, as needed.

The following analogue of John's Theorem \ref{john} can be conveniently phrased in terms of the invariant $d_{BM}^{\rm ctr}$. In a nutshell, an ellipsoid is being replaced by an arbitrary centrally symmetric convex body.  
\begin{thm}[\cite{[L 98],[L 99a]}] \label{maxvol}
Let $C \in {\cal C}^n$ and let $M \in {\cal M}^n$. 
If $M'$ is an affine image of $M$ of maximal possible volume contained in $C$, then $C$ is contained in the homothetic copy of $M'$ with  ratio $2n-1$ and the homothety center in the center of $M'$. The ratio $2n-1$ cannot be lessened as shown by the example of a simplex and a parallelotope.
It follows that, for every $C \in {\cal C}^n$ and every $M \in {\cal M}^n$, we have 
\begin{equation} \label{2n-1} 
d_{BM}^{\rm ctr} (C,M) \leq 2n - 1.
\end{equation}
\end{thm} 

A generalization of 
Theorem \ref{maxvol} (which assumes symmetry of one of the convex bodies) was given by Guo and Kaijser \cite{[GK 03]}, who proved the following (see also \cite{[GK 99],[GK 02]} and the comments on p. 115 in the book \cite {To15} by Toth).

\begin{thm} \label{Guo+Kaijser}
For any two convex bodies $C, D$ in $E^n$ we have
$$d_{BM} (C, D)  \leq (n-1) \min\big({\rm as}(C), {\rm as}(D)\big) +n,$$
where ${\rm as}(\cdot)$ is the Minkowski measure of symmetry.
\end{thm}

Theorem \ref{maxvol}, and the subsequent papers \cite{[GianPT]} by Giannopoulos, Perissinaki and Tsolomitis and \cite{BasRom} by Bastero and Romance,  mobilized Gordon, Litvak, Meyer and Pajor  to prove the following. (Note that once \eqref{Gru-n2} is established, Theorem \ref{maxvol} becomes a special case of Theorem \ref{thm:ctr}.)    

\begin{thm} [\cite {[GLMP]}] \label{thm:GLMP}
Let $C, D \in {\cal C}^n$. We then have 
\begin{equation} \label{Gru-n1}
d_{BM}^\pm (C,D) \leq n \hskip0.3mm .
\end{equation}
Consequently, if  $C \in {\cal C}^n$ and $M \in {\cal M}^n$, then 
\begin{equation} \label{Gru-n2}
d_{BM} (C,M) = d_{BM}^\pm (C,M) \leq n,
\end{equation}
with equality if $C=T$. 
\end{thm}

The above theorem is a far-reaching generalization of Theorem \ref{thm:dBMmod-bound} and Corollary \ref{cor:dBMmod-sym}. Its proof follows similar lines and is based on the following.
\begin{pro}
Let $C,D \in {\cal C}^n$ and let $D'=a(D)$ be such that  
 $C \subset D'$ and, subject to this restriction, $D'$ is an affine image of $D$ of minimal volume. Then $D'$ is contained in a translate of $-nC$. Equivalently, $D'$ is contained in a homothetic image of $C$ with ratio $-n$.  
\end{pro}

Apart from the John's Theorem \ref{john}, Theorem \ref {thm:GLMP} is one of very few instances in this circle of ideas, where a precise answer to a question is known across full range of dimensions. It also suggests an intriguing problem of characterizing the equality cases.  

\begin{prob} If equality holds in \eqref{Gru-n1} (resp., \eqref{Gru-n2}), is it true  that one of $C,D$ is a simplex (resp., that $C$ is a simplex)? 
\end{prob} 

A partial solution to this problem was obtained by Jim\'enez, and Nasz\'odi \cite{JN11}. They show that if $K$ is a convex body, $L$ a strictly convex or a smooth convex body in $E^n$, and $d_{BM}^\pm(K, L) = n$, then $K$ is a simplex.  Further, 
Kobos \cite{Kobos20} proved that if $C$ is not a triangle and $M \in \mathcal M^2$, then $d_{BM}(C,M) < 2$.

\vskip 0.09cm
Next, let us note that while \eqref{Gru-n2} appears to be an improvement on \eqref{2n-1}, the former does not require that the centers of $a(M)$ and the homothety implicit in the definition of $d_{BM}(C,M)$ coincide, and thus it {\em does not} yield a matching bound on $d_{BM}^{\rm ctr}(C, M)$. 
Moreover (see Figure \ref{fig:la25100a} and its analysis in Section \ref{sec:simplex}), if $C =\Tt$ and $M=Q$, then the center of any homothety realizing $d_{BM}(C,M)=2$ is necessarily a contact point of all bodies in the corresponding configuration. So, at least in some cases, we have a strict inequality  $d_{BM}(C,M) < d_{BM}^{\rm ctr}(C, M)$. 
These complications notwithstanding, the question of the optimal bound for $d_{BM}^{\rm ctr}(C, M)$ is fully resolved in the planar case. 

\begin{thm} [\cite{[L 89]}] \label{sqrt 2 +1}
 For every $C \in {\cal C}^2$ and every $M \in {\cal M}^2$ there exists an affine transformation $a$ such that $M \subset a(C) \subset (\sqrt 2 +1)M$. Moreover, the factor $\sqrt 2 +1$ is, in general, the best possible. In other words, 
 \begin{equation} \label{eq:sqrt 2 +1}
    \max_{M\in {\cal M}^2,C \in {\cal C}^2} d_{BM}^{\rm ctr} (C,M) = \sqrt 2 +1.
 \end{equation}
 \end{thm}
We may consider Theorem \ref{sqrt 2 +1}  as a generalization of the John's Theorem \ref{john} to arbitrary two-dimensional real Banach spaces, with the unit ball of that space playing the role of the ellipsoid $E$. The optimality of the factor ${\sqrt 2} +1$ is shown by the example of $C = \Tt$ and $M=Q$. The only two optimal positions of $a(\Tt)$ with respect to $Q$ are shown in Figure 3 of \cite{[L 89]}.  
 For the general dimension, the problem remains open. 
\begin{prob} \label{center}
Given $n\geq 3$, find $\max\left\{d_{BM}^{\rm ctr}(C, M) : C \in {\cal C}^n, M \in {\cal M}^n\right\}$. At the very least, prove an upper estimate better than $2n-1$ or a lower estimate better than $n$.  
\end{prob}
Since the prime suspect for the extremal $C$ in the maximum in Problem \ref{center} is $C=\Tt$ (the simplex), let us mention the following observation.
\begin{thm} \label{ctr-polar}
   If $M \in {\cal M}^n$, then  $d_{BM}^{\rm ctr}(T, M)=d_{BM}^{\rm ctr}(T, M^\circ)$. 
\end{thm}
As was pointed out in Section \ref{sec:intro}, the duality considerations do not -- in general -- work seamlessly in the nonsymmetric case. However, in the present setting a minor miracle happens. Denote $\lambda=d_{BM}^{\rm ctr}(T, M)$ and let $T'=a(T)$ 
be such that $M \subset T' \subset \lambda M$. Since this is equivalent to $\lambda^{-1}M^\circ \subset (T')^\circ \subset M^\circ$ and since $(T')^\circ$ is again a simplex, it follows that  $d_{BM}^{\rm ctr}(T, M^\circ) \leq \lambda$, and we get the conclusion by symmetry. 
 
\vskip 0.09cm

We conclude this section by defining the following invariant, which was introduced in \cite{[L 89]} (at the top of p. 375). Given  $M \in \mathcal M^n$, denote 
$$\overline{\xi}_M =  \sup \{d_{BM}^{\rm ctr}(C, M) : C \in \mathcal C^n\}.$$
In \cite{[L 89]}, this quantity is called {\em capability}, meaning the capability of $M$ to approximate general convex bodies. 
It is shown there that, for any $M \in \mathcal M^2$, we have  $2 \leq \overline{\xi}_M \leq \sqrt 2 +1$.
Next, the class of all $M \in \mathcal M^2$ with the best capability $2$ is described.  Below we recall a question, stated as a conjecture 
on p. 375 in \cite{[L 89]}, which addresses the case $n=3$.

\begin{prob}
Is it true that $\overline{\xi}_M = 3$ for every $M \in \mathcal M^3$? 
\end{prob}

\vskip0.4cm

\section {Centroid Banach-Mazur distance} \label{sectioncentroid}

In this section we will introduce and study another modification of $d_{BM}$, the {\it centroid Banach-Mazur distance} of convex bodies. The literature studying this concept is somewhat limited, so our discussion will be concentrated on the planar case, with one major result and a few minor observations concerning the general case stated at the end of the section. 
\begin{ddd}
Given $C, D \in {\mathcal C}^n$, we define their centroid Banach-Mazur distance by
\begin{equation} \label{def:centroid}
d_{\rm BM}^{\rm cen} (C,D) = \inf_{a, h_\lambda} \{\lambda : \, C \subset a(D) \subset h_\lambda (C) \ {\rm and} \ {\rm cen}(C) =  {\rm cen}(a(D))=  {\rm cen}(h_\lambda(C)) \},
\end{equation}
where ${\rm cen}(K)$ denotes the centroid of $K$ and $h_\lambda$ is a homothety with positive ratio $\lambda$. 
\end{ddd}   
Note that while this is not explicitly required, the center of the homothety $h_\lambda$ coincides then necessarily with the centroid of $C$. Since the requirements for $a$ and $h_\lambda$ in \eqref{def:centroid} are more restrictive than those in the definition of $d_{\rm BM}^{\rm ctr}$, the following holds (cf. Theorem \ref{ctr-polar}). 
\begin{cor} \label{cor:cen}
Let $C\in \mathcal{C}^n$ and $M\in \mathcal{M}^n$. Then
$d_{\rm BM}^{\rm ctr}(C,M) \leq d_{\rm BM}^{\rm cen}(C,M)$ and 
$d_{\rm BM}^{\rm cen}(T,M) = d_{\rm BM}^{\rm cen}(T,M^\circ)$.
\end{cor}

The centroid distance was defined in \cite{[BL]} (where a different notation was used, see the last definition on p. 109) and mentioned in passing on p. 470 of the book \cite{[BMP]} by Brass, Moser and Pach. More recently, it has been investigated  in  \cite{[L 24a]}-\cite{[L 24c]}. 

Similar invariant was introduced by Gr\"unbaum at the bottom of p. 259 of \cite{[Grun63]}. While Gr\"unbaum's definition does not include the requirement concerning the center of the homothety, it appears from the context that this omission was an oversight. Namely, it was pointed out there, and illustrated in Figure 17, that the distance between the parallelogram $P$ and the triangle $T$ is at most $\frac{5}{2}$. However, 
this is correct only if the center of the homothety coincides with the centroids; in absence of such requirement, the optimal ratio appears to be $\frac 73$. Gr\"nbaum also asserts (without proof) that the distance between $P$ and $T$ is exactly $\frac{5}{2}$ and puts forward several conjectures/queries, some of which we restate here. 
\begin{prob} \label{5/2}
Show that $\max_{C\in {\mathcal C}^2} d_{\rm BM}^{\rm cen} (C,T) = \frac{5}{2}$.  Is $\max_{C,D\in {\mathcal C}^2} d_{\rm BM}^{\rm cen} (C,D) = \frac{5}{2}$?
\end{prob}

A partial answer to the first part of Problem \ref{5/2} was already stated as Theorem \ref{thm:T-center} in Section \ref{sec:simplex}. In the present notation, it can be fleshed out as follows.
\begin{thm} [\cite{[L 92],[L 24a]}]\label{thm:5/2}
We have $\max_{M \in {\cal M}^2} d_{\rm BM}^{\rm cen} (M,T) = \frac{5}{2}$. The maximum is attained when $M=P$, where $P$ is a parallelogram.
\end{thm}

The article \cite{[L 24a]} analyzes in detail the pair $(P,T)$. 
Besides giving a detailed argument for 
$d_{BM}^{\rm cen} (P, T) = \frac{5}{2}$, it shows that there are exactly two (up to affine equivalence) different positions of $T$ with respect to $P$ for which $a(T) \subset P \subset h_{5/2}a(T)$, where the origin $o$ is the common centroid of $P$ and $T$. This is illustrated in  Figure  \ref{fig:centroid-BM-dist(P,T)}  ($P$ being there a square). A more transparent illustration, with a regular triangle and parallelogram that are represented as a rectangle and a rhombus, is shown in Figure \ref{fig:centroid-BM-dist(T,P)}. 
\begin{figure}[htbp]
\includegraphics[width=0.55\textwidth]{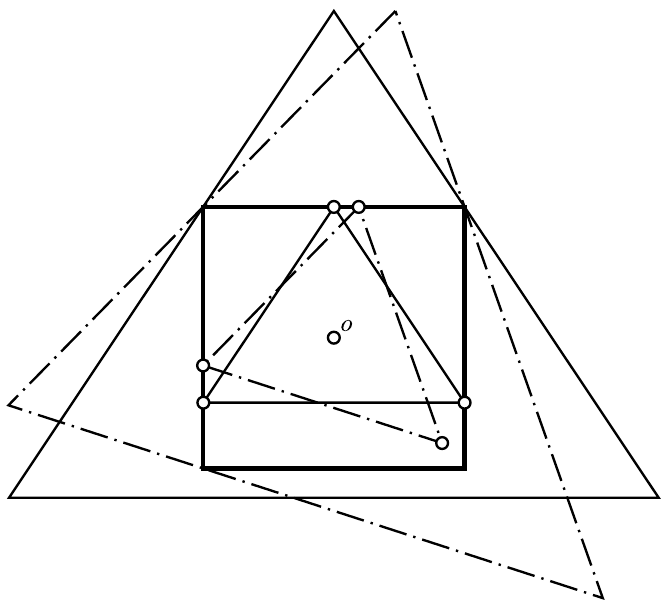}    
\caption{Two different positions of $T$ with respect to $P$ for which $a(T) \subset P \subset h_{5/2}a(T)$.  \label{fig:centroid-BM-dist(P,T)} 
See \cite{[L 21]} for details and clarifications.} 
\end{figure}

\begin{figure}[htbp]
\includegraphics[width=0.5\textwidth]{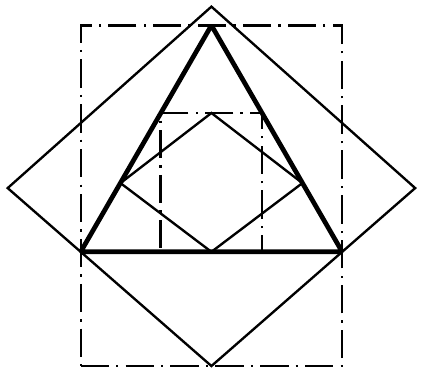} \\   
\caption{Two different positions of $P$ with respect to $T$ for which $a(P) \subset T \subset h_{5/2}a(P)$  \label{fig:centroid-BM-dist(T,P)}.  In this representation it is apparent that these two optimal positions are dual to each other (modulo rescaling and reflection).} 
\end{figure}

A partial result in the direction of the second part of Problem \ref{5/2} was obtained~in~\cite{[L 24c]}.
\begin{thm} [\cite{[L 24c]}] \label{69/17}
The centroid Banach-Mazur distance between arbitrary two planar convex bodies is at most $\frac{69}{17}$.\,\footnote{\, We have been informed by S. D. Kominers \cite{Kominers} that he proved that this upper bound is strictly smaller than $4$, which is slightly better than $\frac{69}{17}$, but still far away from the value of $\frac 52$ conjectured by Gr\"unbaum.}
\end{thm}

In the proof of Theorem \ref{69/17} the following two results are successively used.
The first is a well-known theorem of Besicovitch \cite{Bes} :  {\it For every planar convex body $A$ there exists at least one affine-regular hexagon $H^{ex}$ inscribed in $A$}.
The second is the fact shown in \cite{[L 24c]} and asserting the following : {\it If $H^{ex}$ is an affine-regular hexagon inscribed in a planar convex body $A$, then the centroid of $A$ belongs to the homothetic image of $H^{ex}$ with the ratio $\frac{4}{21}$ and the center in the center of $H^{ex}$}. (The ratio $\frac{4}{21}$ is optimal here.)

\begin{figure}[htbp]
\includegraphics[width=0.5\textwidth]{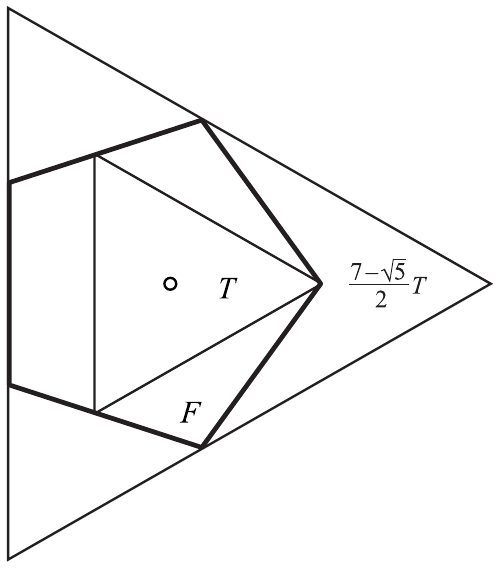} 
\caption{Illustration to the inequality $d_{\rm BM}^{\rm cen} (F, T) \leq \frac{7- \sqrt 5}{2}=:\eta$.  \label{fig:pentagon}}
\end{figure}

\vskip0.09cm
Another result from  \cite{[L 24c]} is that if $F$ is an affine-regular regular pentagon, then  $d_{\rm BM}^{\rm cen} (F, T) \leq \frac{7- \sqrt 5}{2}  =: \eta \approx 2.382$. 
This is illustrated in Figure \ref{fig:pentagon}, which shows the regular pentagon $F \subset \bR^2$ with vertices $(\cos \frac{2\pi}{5}k, \sin \frac{2\pi}{5}k)$, $k=0,1,  \ldots, 4$, the inscribed triangle $T$ with vertex at $(1,0)$ and the opposite side on the  line $x= -\frac{1}{2}$,  
and the homothetic image $\eta T$. 
The centroids of $T$, $F$ and $\eta T$ coincide, they are all at the origin $(0,0)$. We conjecture that this bound is optimal. 
\begin{prob}
Is $d_{\rm BM}^{\rm cen} (F, T) = \eta$?  
\end{prob}
Either way, given that $\eta < \frac 52$ and that the pair $(F,T)$ is  ``the usual suspect'' for extremal behavior among bodies from ${\mathcal C}^2$, the bound shown in \cite{[L 24c]} supplies additional circumstantial evidence that the answer to the second part of Problem \ref{5/2} is affirmative. 

Up to this point, our discussion of the centroid BM-distance was largely restricted to the planar setting. Various instances of $d_{\rm BM}^{\rm cen} (\cdot, \cdot)$ (mostly in higher dimensions) are analyzed below in Sections \ref{non-sym} and \ref{sec:vr} (see Theorem \ref{simplex-H} and Problems \ref{prob:simplex-H} and \ref{simplex-nonH}).
Let us conclude with the following observation of asymptotic nature, thus concerning general dimension, which is a corollary of (the proof of) Theorem \ref{bizeul-klartag}. 
\begin{cor} \label{bizeul-klartag-cor} We have \ 
$\max\{d_{\rm BM}^{\rm cen} (K,L) :  K,L \in \mathcal{C}^n \} =O\big(n\log (1+n)\big)$.  
\end{cor}
What makes the argument work is that the affine images $K'=a(K), L'=b(L)$ that appear in the comments following the statement of 
Theorem \ref{bizeul-klartag} and in \eqref{eq:bizeul-klartag} can be chosen so that they are in the so-called {\em isotropic position}, which implies in particular that their centroids are at $0$, and the latter property is preserved when applying a linear operator (such as a rotation from $\gO(n)$, or multiplication by a scalar, cf. \eqref{eq:bizeul-klartag}). 

While Corollary \ref{bizeul-klartag-cor} essentially resolves the problem of the diameter of $\mathcal{C}^n$ with respect to $d_{\rm BM}^{\rm cen}$ in the asymptotic sense, it does not address the low dimensions.  In the planar case, we do have a candidate for the extremal pair of bodies (see Problem \ref{5/2}), but for $n\geq 3$ the question is wide open. 
\begin{prob} 
For $n=3$ (and other small values of $n$), determine $\max_{C\in {\mathcal C}^2} d_{\rm BM}^{\rm cen} (C,T)$ and $\max_{C,D\in {\mathcal C}^2} d_{\rm BM}^{\rm cen} (C,D)$. At the very least, find instances for which these quantities are strictly larger than $n$.
\end{prob}

\vskip0.4cm


\section{Comparing convex bodies with the parallelotope and the cross-polytope}  
\label{sec:parallelotope} 

We discuss the parallelotope and the cross-polytope in the same section since (when $0$-symmetric) they are related by duality, each being the polar of the other (see \eqref{ellp-ellq}). 

\subsection{The planar case }
Let us start our considerations with the planar case, which is quite particular since $2$-dimensional cross-polytope and parallelogram are affinely equivalent, and both are equivalent to the square $Q$. In the language of normed spaces, the $2$-dimensional $\ell_1$ and $\ell_\infty$ spaces are isometric. Next, it follows from the material presented in the previous sections that answers to many natural questions are known. 

First, by Corollary \ref{M2a}, $d_{BM}(Q,M) \leq \frac 32$ for $M \in {\cal {M}}^2$, with equality iff $M$ is an affine-regular hexagon $H^{ex}$. For $C \in {\mathcal{C}}^2$, $d_{BM}(C,Q) \leq 2$, $d_{BM}^{\rm ctr} (C,Q) \leq 3$,  and $d_{BM}^{\rm cen} (C,Q) \leq \frac 52$, with all inequalities saturated when $C=T$ (by Theorems \ref{maxvol}, \ref{thm:GLMP}, \ref{thm:5/2}, and by \eqref{eq:measure-simplex}). 

\begin{figure}[htbp]
\centering
  \begin{minipage}[c]{0.35\textwidth}
    \includegraphics[width=\linewidth]{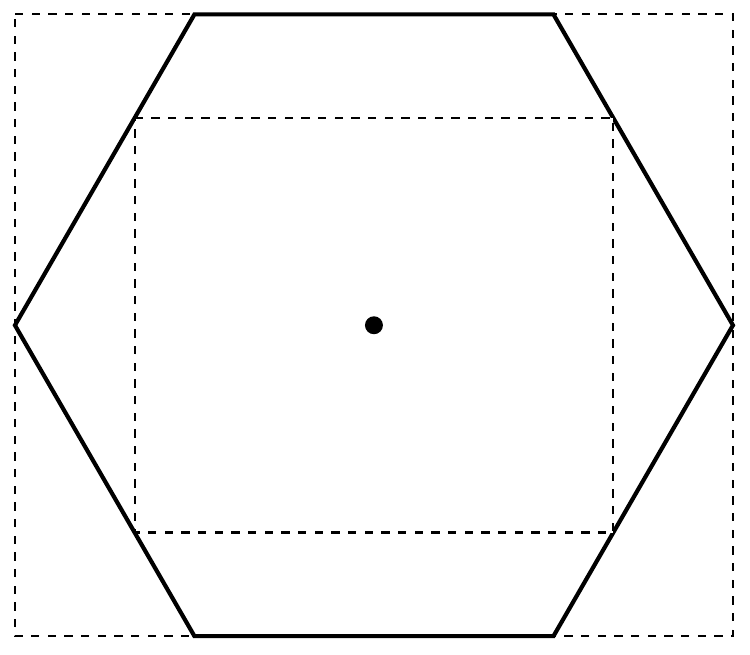}
  \end{minipage}
  \hfill
  \begin{minipage}[c]{0.45\textwidth}
    \includegraphics[width=\linewidth]{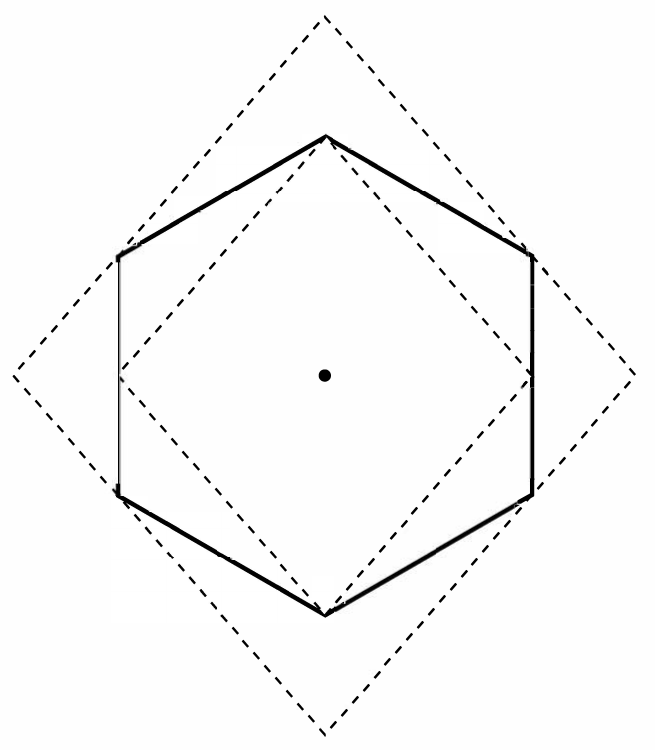}
  \end{minipage}

\vskip0.1cm
\caption{Two optimal pairs of homothetic affine images of $P$ with respect to $H^{ex}$. The second configuration is the dual of the first one : it consists of polars of bodies from the first one (up to scaling). \label{fig:optimumpositions2}} 

\end{figure}

When $M=H^{ex}$, more interesting information is available. As shown 
in \cite{[L 21]}, there are -- up to affine equivalence --  only two optimal positions of a parallelogram $P$ so that $P \subset H^{ex} \subset \frac{3}{2} P$, see Figure \ref{fig:optimumpositions2}. 
This is (modulo elementary calculus) a consequence of the following general fact. 
\begin{pro}[\cite{[L 21]}]
Let  $M \in \mathcal M^2$ and denote $\mu=d_{BM}(P,M)$.
If a particular parallelogram $P_o$ is such that $P_o \subset M \subset \mu P_o$, then $P_o$ is inscribed in $M$ and $\mu P_o$ is circumscribed about $M$. 
\end{pro}

The article \cite{[L 21]}  contains also the following result concerning the BM-distance from $P$ to the remaining affine-regular even-gons, which we include for reference purposes.
\begin{thm} \label{even-gons} 
For a parallelogram $P$ and every $j \in \{1,2, \ldots\}$ we have

\vskip0.1cm
{\rm(i)} $\:\;\;d_{BM} (P , P_{8j}) = \sqrt 2$,

\vskip0.1cm
{\rm (ii)} $\;d_{BM} (P , P_{8j+2}) \leq  \frac{1}{2} \sec \frac{2j}{8j+2}\pi  + \cos \frac{2j}{8j+2}\pi $,

\vskip0.1cm
{\rm (iii)} $d_{BM}(P , P_{8j+4}) = \sqrt 2 \cos\frac{1}{8j+4}\pi$,

\vskip0.1cm
{\rm (iv)}  $d_{BM} (P , P_{8j+6}) \leq  \sin \frac{2j+2}{8j+6}\pi \cdot \csc \frac{4j+2}{8j+6}\pi  + \cos \frac{2j+2}{8j+6}\pi $,

\vskip0.1cm
\noindent
where $P_m$ denotes the affine-regular $m$-gon.
\end{thm}

\noindent Moreover, we conjecture that the following improvement is possible.

\begin{prob}
Show that the inequalities in parts (ii) and (iv) of Theorem \ref{even-gons}  are actually equalities.
\end{prob}

\subsection{$n\geq 3$, the symmetric case}
For general $n$ and for $M \in {\cal M}^n$, it is convenient to keep in mind the language of normed spaces.  In this context, since $(\ell_\infty^n)^*=\ell_1^n$ (cf. \eqref{ellp-ellq}), any statement about $\ell_\infty^n$ will have a companion statement about $\ell_1^n$. Equivalently -- in the language of symmetric convex bodies -- any {\em affinely invariant} statement about $B_\infty^n$ (or the parallelotope $P^n$, or the cube $Q^n$) will have a companion statement about $B_1^n$ (or the cross-polytope). 

As an illustration of that principle, we notice that the radii of $\mathfrak{M}^n$ with parallelotope and with cross-polytope as centers are exactly the same; denote their common value by $\rho_n$. The immediate bounds on $\rho_n$ following from the results of Section \ref{sec:john} (more precisely, from John's Theorem \ref{john} and its sharpness when $K$ is a parallelotope or a cross-polytope) are   
\begin{equation} \label{cube-john}
    \sqrt{n} \leq \rho_n \leq n. 
\end{equation}
Apart from relatively minor improvements, this essentially remained ``the state of the art'' for several decades. For example, it was conjectured by Gr\"unbaum in \cite{[Grun63]} that $\rho_n \leq  \frac{1}{2}(n+1)$, but even the progress towards that modest goal has been elusive.

The first significant progress was achieved around 1990, when it was shown that both bounds are asymptotically imprecise, specifically, $\rho_n = O\big(n \exp(-c\sqrt {\log n})\big) = o(n)$ (for some universal constant $c>0$) and $\rho_n =\Omega\big(\sqrt{n} \log n \big) \gg \sqrt{n}$ (\cite{[BoSz88]} and \cite{Sza90}). In particular, it follows that the cross-polytope or the parallelotope cannot be an asymptotic center of the Banach-Mazur compactum (cf. Theorem \ref{center-non-unique}).

In the intervening years, both bounds have been successively sharpened. The upper bound for $\rho_n$ was shown to be $O(n^{7/8})$ in \cite{[SzTa1]} (see also \cite{[SzTa2]}), which was later improved to $O(n^{5/6})$ by Giannopoulos \cite{[Gi 95]}, to $(2n)^{5/6}$ by Youssef \cite{Youssef14}, and, very recently, to $O(n^{2/3})$ by Friedland, Naor and Youssef \cite{FNY}. For the lower bound, Tikhomirov \cite {Tik19} showed that $\rho_n = \Omega(n^{5/9}\log^{-\beta} n)$, where $\beta > 0$ is some universal constant. This was improved by Friedland  to $\Omega(n^{4/7}\log^{-\beta} n)$ in \cite{Fri1} and to $\Omega(n^{5/8}\log^{-1/4} n)$ in \cite{Fri2}. (Similar estimate, with a slighly larger power of the logarithm, was obtained independently and approximately at the same time by Hmadi \cite{Hmadi26}.) Finally, Friedland \cite{Fri3} supplied an AI-assisted proof giving an $\Omega(n^{2/3}\log^{-10} n)$ bound. Accordingly, the state of the art as of the time of this writing is given by the following.
\begin{thm} \label{thm:rad-cube}
 If $\rho_n:= \max_{M\in {\cal M}^n}d_{BM}(P,M) = \max_{\dim X = n}d(\ell_1^n)$, then 
 \begin{equation} \label{eq:rad-cube}
    c_1 \, n^{2/3}\log^{-10} (1+n) \leq \rho_n \leq c_2 \, n^{2/3},
 \end{equation}
 where $c_1,c_2 > 0$ are universal constants. 
\end{thm}
The asymptotic behavior of the sequence $(\rho_n)$ has been thus determined up to a polylogarithmic factor. To fully close the gap is likely a difficult problem. 
\begin{prob}
Determine the asymptotic order $\Theta(\rho_n)$. At the very least, improve the exponent $-10$ on the left hand side of  \eqref{eq:rad-cube}.
\end{prob}

\vskip 0.09cm We will now describe briefly the ideas  behind the proofs in \cite{[BoSz88]} and \cite{Sza90}, and in most of the works that followed. For the lower bound on $\rho_n$, bodies that are far from the cross-polytope are random polytopes obtained by one of the schemes (i)-(iii) from the discussion that followed Theorem \ref{gluskin}, with $N$ being linear or polynomial in the dimension $n$ (say, $N=2n$, $N=n^2$, or $N=n^3$). 

In functional analytic terms, the corresponding random $n$-dimensional normed space $X$ can be represented as a (random) quotient of $\ell_1^N$, with the quotient map being an $n\times N$ matrix whose entries are independent standard Gaussian variables. Then, we need to show that there is no operator $u\in \gGL(n)$ such that
\begin{equation}\label{cube-bound}
\|u : \ell_1^n \to X\| \leq 1 \quad \hbox{and} \quad \|u^{-1} : X \to \ell_1^n\|\leq \rho
\end{equation}
if $\rho < \rho_n$ (or whatever lower bound on $\rho_n$ we want to establish). As in the case of Theorem  \ref{gluskin}, the argument consists of two steps : a sharp probability estimate for a fixed operator (say, for the probability of the second inequality in \eqref{cube-bound}, conditioned on the first inequality), followed by discretization and the union bound. 

For the upper bound, the main technical step involved finding a ``nice'' subspace of an arbitrary $n$ dimensional $X$ of small codimension whose distance to an $\ell_1$-space is ``not too large.'' These requirements can be quantified, e.g.,  by the so-called {\em proportional Dvoretzky-Rogers factorization} (cf. \cite{DR50,[BoSz88]}) : given $\theta \in (0,1)$,  there exists $m\geq \theta n$ such that the formal identity map $i_{12} : \ell_1^m \to \ell_2^m$ can be represented as a composition $i_{12} = v \circ w$ with
\begin{equation} \label{DR}
    \|w : \ell_1^m \to X\| \times \|v : X \to \ell_2^m\| \leq c(\theta),
\end{equation}
where $c(\theta)$ depends only on $\theta$ (and not on $X$ or $n$). The bound \eqref{DR} encodes both the ``closeness'' of the $m$-dimensional subspace $E=w(\ell_1^m)\subset X$ to $\ell_1^m$ (specifically, $d(E,\ell_1^m) \leq c(\theta)\sqrt{m}$\,) and the ``niceness'' of $E$ as a subspace of $X$ (the fact that $v$ is defined on $X$ and not just on $E$). Then -- using \eqref{DR} or a similar property -- one proceeds to construct an operator $u$ verifying \eqref{cube-bound} for an appropriate $\rho$, with the magnitude of $\rho$ depending on the quality of the bound $c(\theta)$. (To the best of our knowledge, the best known dependence on $\theta$, $c(\theta) =O\big((1-\theta)^{-1}\big)$, was obtained in \cite{[Gi 96]}, with constants improved in  \cite{Youssef14}.) The most recent work \cite{FNY} dispenses with this ``piecemeal'' approach and selects $u$ verifying \eqref{cube-bound} ``in one step,'' which allows to improve the estimate.

While the asymptotic bounds such as \eqref{eq:rad-cube} are quite impressive, they are not necessarily useful for small $n$, particularly if they involve unspecified ``universal constants.'' Those constants are typically ``effectively computable,'' but more often than not they yield reasonable estimates only for very large values of $n$. Note that even the ``nice looking'' upper bound from \eqref{eq:rad-cube} is worse than the easy bound from  \eqref{cube-john} if $n<2^5=32$. Accordingly, it is of interest to come up with nontrivial bounds that work for all $n$ and not just asymptotically. 

In this direction, it was shown in Theorem 3 of  \cite{[L 91]} that, for any   $M \in {\mathcal M}^n$,
\begin{equation}\label{n^2-n+1}
 d_{BM} (M,P) =  d_{BM}^{\rm cen} (M,P) \leq \sqrt{n^2 -n +1} \end{equation}
and consequently 
\begin{equation}\label{rho_n^2-n+1}
\rho_n \leq \sqrt{n^2 -n +1} \,.
\end{equation}

\noindent 
Asymptotically, as $n \to \infty$, this estimate behaves only as $n-\frac 12$, only slightly better than the upper bound in \eqref{cube-john}. 
However, for small $n$ we get concrete estimates, e.g., $\rho_3 \leq \sqrt 7$, and the bound \eqref{n^2-n+1} was used by Kim and Nakprasit \cite{[KN]} in the proof of their main Theorem 2 on the chromatic number of the graph of the complement of the intersection of a family of translations of a convex body of $E^n$. 
 
\vskip 0.09cm
It is worthwhile to mention that in the proof of Theorem 3 of  \cite{[L 91]} the author constructed a convex body $W_n \in \mathcal M^n$ and conjectured  that $W_n$ is more central in $\mathcal M^n$ than $B^n$, i.e., that the radius of $\mathcal M^n$ with the center $W_n$ is strictly smaller than $\sqrt n$. If true, this would solve the first part of Problem \ref{rad-sym}.

Taschuk \cite{Taschuk} slightly improved the upper bound from \eqref{rho_n^2-n+1} for $n\geq 3$ by proving the following. 
\begin{thm} \label{Taschuk} We have 
\begin{equation} \label{eq:Taschuk}
\rho_n    \leq \sqrt{(n-1)^2+1 +2/(\sqrt {n+2}-1)}.  
\end{equation}
\end{thm}
As $n\to \infty$, the bound \eqref{eq:Taschuk} behaves as $n-1 \gg (2n)^{5/6}$, but it is smaller than  $(2n)^{5/6}$ for $n \leq 37$. (The bound from \eqref{rho_n^2-n+1} is smaller than  $(2n)^{5/6}$ for $n\leq 34$. For $n=2$, \eqref{eq:Taschuk} yields the trivial bound $\rho_2\leq 2$, but otherwise  both bounds are {\em always} strictly smaller than $n$.)

The paper \cite{FNY} doesn't state explicitly the value of the constant $c_2$ from \eqref{eq:rad-cube}. While an examination of the argument will likely lower the threshold of $n=37$ mentioned above, the bound \eqref{eq:Taschuk} will remain competitive for sufficiently small dimensions.

\subsubsection{The distance between the parallelotope and the cross-polytope} \label{para-cross}
An interesting question that involves both the parallelotope and the cross-polytope, and that leads to several deep problems in other areas of mathematics,  concerns the value of the Banach-Mazur distance between them. In functional analytic terms, we want to find $d(\ell_1^n,\ell_\infty^n)$.  While it has been known at least since 1960's (Gurari\u\i, Kadets and Matsaev \cite{GKM65}) that $d(\ell_1^n,\ell_\infty^n) = \Theta(n^{1/2})$, obtaining explicit bounds with good constants is harder. One result that can be found in the literature (Xue \cite{Xue 18}; see also his dissertation \cite{Xue 19} and \cite{GKM66}) is as follows.\footnote{Actually, \cite{Xue 18} contains a proof of a weaker form with $\alpha \approx 0.5832$ in place of $\frac 1{\sqrt{2}}$. The stronger inequality \eqref{eq:Xue} is formulated there as a guess, which would follow from an inequality stated as Conjecture 5.1. As was pointed out in \cite{[KV 24]} (and is implicit in \cite{GKM66}), that conjecture is an instance of Khinchine's inequality with the optimal constant, which was known since 1970's (\cite{Sz76,Haa82}).}
\begin{thm}  \label{thm:Xue} We have 
\begin{equation} \label{eq:Xue}
    \tfrac 1{\sqrt{2}} \sqrt{n} \leq  d(\ell_1^n,\ell_\infty^n) \leq (1+\sqrt{2}) \sqrt{n}.    
\end{equation}
\end{thm}
The constant in the first inequality in \eqref{eq:Xue} is clearly optimal, as seen by considering $n=2$. The constant in the second inequality is unlikely to be optimal. 
\begin{prob} \label{1-infty}
Improve the constant in the second inequality in \eqref{eq:Xue}. 
\end{prob}

The article \cite{Xue 18} by Xue contains also various numerical results, which in particular suggested that the BM-distance between the $3$-dimensional cube and cross-polytope is $\frac{9}{5}$. This was recently confirmed by Kobos and Varivoda \cite{[KV 24]}, who showed the following.
\begin{thm} 
\label{KobosVarivoda}
The BM-distance between the parallelotope and  cross-polytope is $\frac{9}{5}$ in dimension three and $2$ in dimension four. 
\end{thm}
Since $\frac 95 > \sqrt{3}$, we conclude that the constant in the second inequality in \eqref{thm:Xue} can not be $1$. However, it is conceivable that  $d(\ell_1^n,\ell_\infty^n) \leq \sqrt{n}$ for $n$ sufficiently large; the numerical results from \cite{Xue 18} show that this is the case when $n=6,7,8$, in fact with a strict inequality. (For $n=5$ the numerical results are not conclusive.) 

The authors of \cite{[KV 24]} also proved that if, for some  $u\in \gGL(3)$, we have $\|u : \ell_1^3 \to \ell_\infty^3\|\leq 1$ and $\|u^{-1} : \ell_\infty^3 \to \ell_1^3 \| \leq \frac{5}{9}$, then the matrix of $u$ must have the form 
\begin{equation} \label{special}
{\left[ \begin{matrix} 
1/3&-1&-1 \cr 
-1&1/3&-1 \cr 
-1&-1&1/3 \cr
\end{matrix}
\right]}
  \end{equation}
or any equivalent form obtained by permuting rows and/or columns, and multiplying some rows and/or columns by $-1$.
 
The article \cite{Xue 18} highlights the well-known connection of the question at hand to the concept of Hadamard matrices (\cite{Sylvester, Hadamard}). 
\begin{ddd} \label{def:H}
 A {\em Hadamard matrix} $H$ is an $n \times n$ matrix whose all entries are $\pm 1$ and  whose rows (or, equivalently, columns) are pairwise orthogonal. If such $H$ exists, we will call $n$ a {\em Hadamard number.}    
\end{ddd}
It is easy to see that if $H$ is a Hadamard matrix, then $H^{-1}=\frac 1n H^t$, $\|H : \ell_1^n \to \ell_\infty^n\| = 1$, and  $\|H^{-1} : \ell_\infty^n \to \ell_1^n\| \leq \sqrt{n}\; \|H^{-1} : \ell_\infty^n \to \ell_2^n\|= \sqrt{n}$. Accordingly, if $n$ is a Hadamard number, then $d(\ell_1^n,\ell_\infty^n) \leq \sqrt{n}$. Similarly, if an $m$ is a Hadamard number and $n=m+k$, then 
\begin{equation} \label{useH}
    d(\ell_1^n,\ell_\infty^n) \leq \sqrt{m} + d(\ell_1^k,\ell_\infty^k).  
\end{equation}
(If $k=0$, we interpret the second term on the right hand side as $0$.)  

The argument from \cite{Xue 18} was based on \eqref{useH} and the (easy) fact that every power of $2$ is a Hadamard number. A renown {\em Hadamard conjecture} states that every $k\in \bN$, $n=4k$ is a Hadamard number. (At the time of this writing, the smallest $n=4k$ for which the existence is unknown is $668$ (see \cite{CD-handbook,Dok10}). If the conjecture was true, it would follow from \eqref{useH} and Theorem \ref{KobosVarivoda} that, for $n\geq 3$,  
\begin{equation}
    d(\ell_1^n,\ell_\infty^n) \leq \sqrt{n-3}+\tfrac95 .
\end{equation}
What can be shown unconditionally is weaker, but still asymptotically better than the upper bound in \eqref{eq:Xue}.
\begin{thm} \label{thm:Paley} We have 
\begin{equation} \label{Paley}
    d(\ell_1^n,\ell_\infty^n) \leq \sqrt{n}+ O\big(n^{21/80}\big).
\end{equation}  
\end{thm}
This is an easy consequence of several known deep results. First, by \eqref{useH} and \eqref{eq:Xue}, it is enough to show that, given $n$, there is a Hadamard number $m \leq n$ such that $n-m = O(n^{21/40}$). Now, by a construction due to Paley \cite{Paley}, if $q$ is an odd prime power, then $2(q+1)$ is a Hadamard number. In particular, the gaps in the sequence of Hadamard numbers are at most twice the gaps in the sequence of prime numbers. At the time of this writing, the best upper bound on gaps for primes in the interval $[2,x]$ is $O(x^{21/40})$ (see \cite{BHP01}), which is exactly what we need. 

Since Theorem \ref{Paley} implies that $d(\ell_1^n,\ell_\infty^n) \leq (1+o(1))\sqrt{n}$, it is very likely that -- by chasing the constants implicit in the $O(\cdot)$ bounds and  some numerical work -- one can achieve some progress towards Problem \ref{1-infty}, and perhaps even find the optimal constant.  We will discuss improvements to the {\em lower }bound in \eqref{eq:Xue} in  subsection \ref{sec:projection}. 

\subsubsection{The projection constant} \label{sec:projection}
Since the BM-distance controls many invariants of normed spaces/convex bodies, those invariants can be used to obtain lower bounds on the distance between spaces/bodies. One such invariant was the Minkowski measure of symmetry ${\rm as}(\cdot)$ discussed in Section \ref{sec:simplex}. Here we will consider another invariant, which is more tailored to the setting of symmetric bodies/normed spaces. 
\begin{ddd}\label{def:projection}
 Given $n$-dimensional normed space $X$, denote by $I_X$ the identity map on $X$ and define the projection constant of $X$ as 
 \begin{equation} \label{eq:projection}
   \lambda(X) = \inf \|u\| \cdot \|v\|  
 \end{equation}
 over all $N\geq n$ and all pairs of operators $u : \ell_\infty^N \to X$, $v : X \to \ell_\infty^N$ such that $u \circ v = I_X$. 
\end{ddd}
This is not the usual definition of $\lambda(X)$, but it is easy to show that it is equivalent to the more common ones (see, e.g.,  \cite{[Tom 89]}, Propositions 13.3 and 32.1). It is clear that 
\begin{equation} \label{proj-dist}
    \lambda(X) \leq d_{BM}(X,\ell_\infty^n).  
\end{equation}
Indeed, $d_{BM}(X,\ell_\infty^n)$ is given by the same expression \eqref{eq:projection}, but restricted to $N=n$. The usefullness of $\lambda(\cdot)$ follows in part from the fact that, for many natural spaces $X$, $\lambda(X)$ can be evaluated precisely or at least approximately. For example, already in 1961 Gr\"unbaum calculated 
$\lambda(\ell_1^n)$ (see \cite{[Grun60],[KTS99]}) .We have, for $k \in \bN$, 
\begin{equation} 
    \lambda(\ell_1^{2k})=  \lambda(\ell_1^{2k-1})= 2 B(k, \tfrac 12)^{-1} =\tfrac 2{\sqrt{\pi}}\Gamma(k+\tfrac 12)/\Gamma(k),
\end{equation}
where $B(\cdot,\cdot)$ and $\Gamma(\cdot)$ are respectively the {\em Beta} and the {\em Gamma} functions. 
It follows that  $\lambda(\ell_1^n) \geq \frac 1{\sqrt{2}} \sqrt{n}$ (with equality iff $n=2$; this gives an alternate proof of the first inequality in \eqref{eq:Xue}) and that   $\frac{\lambda(\ell_1^n)} {\sqrt{n}} \to \sqrt{\frac{2}{\pi}}$, which implies that   
\begin{equation} \label{2/pi}
\liminf_n \tfrac{d_{BM}(\ell_1^n,\ell_\infty^n)}{\sqrt{n}} \geq \sqrt{\frac{2}{\pi}}.
\end{equation}
In the asymptotic sense, \eqref{2/pi} improves the first inequality from  \eqref{eq:Xue}. We can then ask
\begin{prob}
Is the constant $\sqrt{\frac{2}{\pi}}$ in \eqref{2/pi} optimal?
\end{prob}
It follows from \eqref{Paley} that the constant in question can not be strictly larger than $1$. 

One of the most important open problems in this circle of ideas is whether, in some sense, the inequality \eqref{proj-dist} can be reversed. This question is often referred to as ``the $P_\lambda$ problem'' and its precise statement is as follows.   
\begin{prob}
Does there exists a function $f : [1,\infty) \to [1,\infty)$ such that, if $n\in \bN$ and if $X$ is an $n$-dimensional normed space satisfying  $\lambda(X) \leq t$, then $d_{BM}(\ell_\infty^n, X) \leq f(t)$?
\end{prob}
Partial positive results have been obtained by Nachbin \cite{Nachbin} (for $t=1$, with $f(1)=1$), by Zippin \cite{Zip81, Zip84} (for $t\in [1,1+\ve_0]$, for some $\ve_0 > 0$), and Sch\"utt \cite{Sch}  (under the hypothesis that $X$ is a symmetric normed space, see \eqref{def:symmetric} for a definition). It follows from \cite{Fri3} that, as $t \to\infty$, $f(t)$ -- if it exists at all - must grow at least as fast as $t^{4/3}$, up to a polylogarithmic factor.

\subsection{$n\geq 3$, the not-necessarily-symmetric case} \label{non-sym}

Since the parallelotope and the cross-polytope are centrally symmetric, general results from Section \ref{sec:one-symmetric} apply directly in the present setting. In particular, it follows that if $M$ is either the $n$-dimensional parallelotope or the cross-polytope and $C \in {\mathcal C}^n$, then \\
$1^\circ$ $d_{BM} (C,M) = d_{BM}^\pm (C,M) \leq n$ with equality when $C=\Tt$ (by Theorem \ref{thm:GLMP}), \\
$2^\circ$ 
$d_{BM}^{\rm ctr} (C,M) \leq 2n - 1$  (by Theorem \ref{maxvol}). 

We know from Theorem \ref{sqrt 2 +1} that $2^\circ$ is not precise when $n=2$. The following is a special case of Problem \ref{center}, which addresses the issue at hand in higher dimensions. 
\begin{prob} \label{C-P-ctr}
If  $n\geq 3$ and $M$ is either the $n$-dimensional parallelotope or the $n$-dimensional cross-polytope, find \ $\max\left\{d_{BM}^{\rm ctr}(C, M) : C \in {\cal C}^n\right\}$.
\end{prob}
An interesting instance of Problem \ref{C-P-ctr} (also relevant to Problem \ref{center}) is that of $C=T^n$. We have the following curious observation. 
\begin{thm} \label{simplex-H}
If $n+1$ is a Hadamard number (in particular, if $n+1$ is a power of $2$) and $M$ is either the $n$-dimensional parallelotope or the $n$-dimensional cross-polytope, then 
   \begin{equation} \label{eq:simplex-H}
    d_{\rm BM}^{ctr} (T,M)= d_{\rm BM}^{\rm cen} (T,M) = n.
   \end{equation}
\end{thm}
It will be convenient to postpone the proof of Theorem \ref{simplex-H} until 
Section \ref{sec:vr} (see the paragraph following Problem \ref{prob:vrPT}). The following follow-up problem is probably not very difficult. 
\begin{prob} \label{prob:simplex-H}
Do equalities in \eqref{eq:simplex-H} (both of them and separately) characterize Hadamard numbers?    
\end{prob} 
A more ambitious problem is the following.
\begin{prob} \label{prob:dTP}
Given $n\geq 3$, characterize all configurations certifying $d_{BM}(T,P)=n$. More concretely, characterize all simplices $T \subset Q:=[-1,1]^n$ that contain a translate of $Q/n$. 
\end{prob}
For $n=2$, the problem was solved in Section \ref{sec:simplex} (see the comments following Figure \ref{fig:la25100a}). The argument presented there can be easily modified to show the following more general fact, which can be seen as a step towards answering Problem \ref{prob:dTP} which, however, is not explicit enough to have a bearing on questions such as Problem \ref{prob:simplex-H}. 
\begin{pro} \label{prop:dTM}
Let $M \in \mathcal{M}^n$. The configuration $M \subset T$ certifies $d_{BM}(T,M)=n$ (that is, $M$ contains a translate of \hskip0.3mm $T/n$) if and only if the centroids of the faces of $T$ belong to $M$. 
\end{pro}
We note, however, that (for example) for $n=3$ there are many non-equivalent configurations satisfying the requirements from Problem \ref{prob:dTP}, so the answer to that problem -- if at all feasible -- may be quite complicated.

The general results $1^\circ$ and $2^\circ$ above notwithstanding, we find it enlightening to present a direct geometric proof of the bound $1^\circ$ when $M=P$. 
The following classical result is due to Chakerian and Stein \cite{[CS 67]}. 

\begin{thm} \label{Chakerian+Stein}
If $Q$ is a parallelotope of the minimum volume containing a convex body $C\subset E^n$, then $C$ contains a homothetic copy of $Q$ of ratio $\frac{1}{n}$.
\end{thm}

On the other hand, we have, in a sense, the following dual theorem \cite{[L 91]}.

\begin{thm} \label{MaxVolContained}
If $P$ is a parallelotope of maximal volume contained in a convex body $C\subset E^n$, then $C$ is contained a homothetic copy of $P$ of ratio $n$.
\end{thm}

\begin{figure}[htbp]

\includegraphics[width=0.60\textwidth]{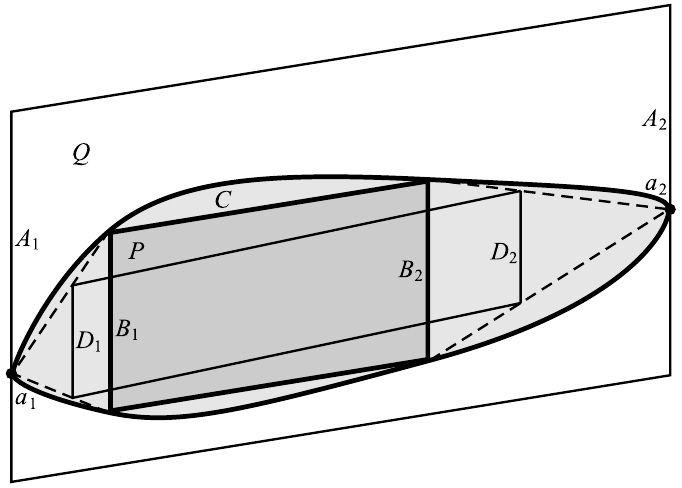} \\   

\vskip 0.1cm
\caption{A parallelotope of the largest volume in a convex body; the case of $n=2$. \label{fig:largestparallelotope} } 

\end{figure}

Here is a sketch of the proof of Theorem \ref{MaxVolContained}. 
Let $Q$ be a homothetic image of $P$ with the smallest positive ratio $r$ such that $C \subset Q$, see Figure \ref{fig:largestparallelotope}.
Suppose that $r>n$.
By the minimality of $Q$ there are parallel facets $A_1$ and $A_2$ of $Q$ whose intersection with $C$ is non-empty and, 
for $i=1, 2$, choose $a_i \in C \cap A_i$. Next, let $B_1, B_2$ be the facets of $P$ parallel to  $A_1, A_2$ numbered so that $B_1$ is closer to $A_1$ than $B_2$ is (and $B_2$ is closer to $A_2$ than $B_2$ is). 
Further, let $D_i$ be the image of $B_i$ under the homothety with ratio $\mu = \frac{(n-1)r}{n(r-1)}$ and center $a_i$. (Note that the assumption $r>n$ implies that $0< \mu <1$.) 
Consider the parallelotope $D= {\rm conv}(D_1 \cup D_2)$ and 
denote by $\delta_i$ the ratio of the distances between the hyperplanes containing $A_i$ and $B_i$, to these containing $B_1$ and $B_2$, then $\delta_1+\delta_2 = r-1$.
A simple calculation shows that the ratio of the distances between the hyperplanes containing $D_1, D_2$ to those containing $B_1, B_2$ equals $\frac{r}{n}$.
Thus $\frac{{\rm Vol}(D)}{{\rm Vol}(P)} = \mu^{n-1}\frac{r}{n} = f(r)$, where $f(x) =\frac{(n-1)^{n-1}}{n^n}\frac{x^n}{(x-1)^{n-1}}$.
Since its derivative $f'(x)=\frac{(n-1)^{n-1}}{n^n}\big(1 - \frac{n}{x}\big)\big(\frac{x}{x-1}\big)^n$ is positive for $x>n$ and since $f$ is continuous on $[n,\infty)$, our hypothesis $r>n$ implies that $f(r)>f(n)=1$. 
Consequently, $\frac{{\rm Vol}(D)}{{\rm Vol}(P)} > 1$, which contradicts the choice of $P$.
This shows that $r \leq n$, as needed.

\vskip0.09cm
Note that most often we get  pairs of different parallelotopes from these two theorems. For instance, consider as $C$ the hexagon with vertices $(\pm 2 , 1), (\pm 2, -1), (\pm 3, 0)$; we leave the details as an exercise for the reader. 

It would be interesting to have a version of the above argument with parallelotopes replaced by cross-polytopes. 

\begin{prob} \label{cross-elementary} 
Produce an elementary proof of an analogue of Theorem \ref{Chakerian+Stein} (or Theorem \ref{MaxVolContained}) with $P$ (resp., $Q$) replaced by a cross-polytope. 
\end{prob}
As noted at the beginning of this section, we know that such analogues do hold as a consequence of Theorem \ref{thm:GLMP}. We do not immediately see an elementary proof, but in all likelihood it can be extracted from  the arguments in \cite {[GLMP]}.  

We conclude this section by recalling a problem that was posed in \cite{[L 93]}. We need the following notation. Given convex body $C\in \mathcal{C}^n$ and a unit vector $u \in S^{n-1}$ we define the {\em width} of $C$ in the direction of $u$ as 
\begin{equation}
  w(C,u)  := \max_{x\in C} \, \langle x,u\rangle \, - \, \min_{x\in C} \, \langle x,u\rangle
\end{equation}
We can now state  
\begin{prob} Let $C\in \mathcal{C}^n$ and let $P\subset C$ be a parallelotope inscribed in $C$ (that is, all vertices of $P$ belong to $\partial C$). Next, let $u_i$, $i = 1, \ldots, n$, be the unit vector perpendicular to
the $i$th pair of parallel facets of $P$. Is it true that we necessarily have  
 \begin{equation}
     \sum_{i=1}^n \tfrac 1{w(C,u_i)} \geq 1\; ? 
 \end{equation}   
\end{prob}
 It was shown in \cite{[L 93]} that the answer is affirmative if $n=2$. Various partial results were presented in \cite{[L 01]}  and \cite{Nev11}.  

\subsection{Approximation by rectangles, rectangular parallelotopes, and similar bodies} \label{sec:rectangles}
 
The present section is sort of an outlier in that it considers concepts that are {\em not necessarily} affinely invariant. 
We start with the following (somewhat vague) definition. 
\begin{ddd} \label{class} 
Let $\cal D$ denote a class of convex bodies. 
We then set 
  \begin{equation} \label{eq:class}
      d_{BM}^{\rm class}(C,{\cal D}) = \inf_{D \in {\cal D}, h_{\lambda}} \{ \lambda \ :\ \ D \subset C \subset h_{\lambda}(D) \},
  \end{equation}  
  where $h_\lambda$ is a homothety with a positive ratio $\lambda$.
\end{ddd}
The quantity defined by \eqref{eq:class} is clearly a generalization of the BM-distance. Indeed, we have $d_{BM}(C,D) = d_{BM}^{\rm class}(C,{\cal D})$, where  ${\cal D}$ consists of bodies that are affinely equivalent to $D$. Other affinely-invariant classes can be legimately considered (see, e.g., Theorem \ref{symmetric-ellp} in Section \ref{sec:classes} and the comments that follow it), but here we will concentrate on a different direction, which is of interest in computational geometry. 

\noindent
\vskip 0.09cm
Let ${\cal K}^n$ be the class of all rectangular parallelotopes in $E^n$.
P\'olya and Szeg\"o \cite{[PoSe]} showed that for every convex body $C \in \mathcal{C}^2$  there is a pair of homothetic rectangles $R, R'=h_{\lambda}R$ with $1 \leq \lambda \leq 3$ such that $R \subset C \subset R'$. An optimal result in this direction (which also supplies additional information) was obtained in \cite{[L 93]}. 
\begin{thm} \label{rectangles}
In any convex body $C \in \mathcal{C}^2$ we can inscribe a rectangle $R$ such that a homothetic copy $R'$ of $R$ is circumscribed about $C$, with 
the (positive) ratio $\lambda$ of the homothety being at most $2$. In other words, for such $C$, 
\begin{equation} \label{eq:rectangles1}
    d_{BM}^{\rm class} (C,{\cal K}^2) \leq  2  
\end{equation} 
and the inequality is saturated when $C=T$. Furthermore, 
\begin{equation} \label{eq:rectangles2}
   \tfrac{1}{2}\cdot{\rm Area}(R') \leq {\rm Area}(C) \leq 2\cdot{\rm Area}(R). 
\end{equation}
\end{thm}

The fact that the bound  in \eqref{eq:rectangles1} cannot be improved follows from  \eqref{eq:measure-simplex}.   
In \cite{[SFRW]} we find a further discussion of the results from \cite{[L 93]}.
In particular, the authors prove that if $C$ is a polygon with $m$ vertices, then this best pair  of rectangles $R, S=h_\lambda (R)$ (i.e., with the smallest possible $\lambda$) can be computed in time $O(\log^2m)$. 
 
\vskip 0.09cm Motivated by Theorem \ref{rectangles} (and by Theorems \ref{Chakerian+Stein}
and \ref{MaxVolContained} from Section \ref{non-sym}), 
we ask the following question.

\begin{prob} \label{prob:rectangles}
Does the bound from Theorem \ref{rectangles} generalize to $d_{BM}^{\rm class} (C, {\cal K}^n) \leq n$, valid for every dimension~$n$ and for every $C\in \mathcal{C}^n$?\;\footnote{A much weaker conjecture was stated by Hadwiger \cite {[Had]}, who asked if there exist rectangular parallelotopes  $P_1 \subset C \subset P_2$ with mutually parallel faces and the ratio of volumes not exceeding $n! \cdot n^n$.}
\end{prob}
Again, we see from \eqref{eq:measure-simplex}  that a better estimate is impossible. A modest evidence in favor of the affirmative answer is the following result due to Makeev \cite{[Mak96]} (who also asked the question from Problem \ref{prob:rectangles}). 
\begin{thm} \label{thm:Makeev} 
For every $C\in {\mathcal C}^n$ there exist two {\bf similar} rectangular parallelotopes $P_1$ and $P_2$,  with similarity ratio not exceeding $n$,  such that $P_1 \subset K \subset P_2$.
\end{thm}

We conclude this section with a selection of results and problems about other similar classes of bodies. 
First, let us denote the class of rhombi by $\cal R$. Since the polar of a rectangle centered at the origin is a rhombus, the classes $\cal R$ and ${\cal K}^2$ are dual to each other, and we have an analogue of Theorem \ref{rectangles}. 
\begin{thm}[\cite{[L 02]}] \label{rhombi}
All statements in Theorem \ref{rectangles} remain true if we replace rectangles by rhombi and the class ${\cal K}^2$ by $\cal R$. 
\end{thm}
In particular, for the homothetic rhombi $R \subset C \subset R'=h_\lambda(R)$ implicit in the assertion, $R$ is inscribed in $C$ and $R'$ is circumscribed about $C$. 
This answered the question of de'Valcourt \cite{de'Val} about such rhombi. 

We point out that while Theorems \ref{rectangles} and \ref{rhombi} are dual to each other, neither of them formally follows from the other.  This is due to the fact that the center of the homothety does not necessarily coincide with the centers of the rectangles or rhombi in question -- a complication that was already noted in Section \ref{sec:intro},  specifically in the paragraphs following \eqref{eq:polar}. Consequently, the answer to the following question is not formally implied by a solution to Problem \ref{prob:rectangles}. 
\begin{prob} \label{prob:rhombi} Let ${\cal R}^n$ be the class of $n$-dimensional rectangular orthoplexes (i.e., the polars of rectangular parallelotopes with center at $0$, and translates of such polars). 
Does the bound from Theorem \ref{rhombi} generalize to $d_{BM}^{\rm class} (C, {\cal R}^n) \leq n$, valid for every dimension~$n$ and for every $C\in \mathcal{C}^n$?
\end{prob}

\vskip 0.09cm
Another result from \cite{[L 02]} is as follows.  For every  $C \in {\cal C}^2$ there exists an isosceles triangle $T \subset C$ such that a homothetic image $T'$ of it with ratio at most $\frac{7}{2}$ contains $C$ (again, see \cite{[L 02]}).
In other words, $d_{BM}^{\rm class} (C, {\cal T}_{I}) \leq \frac{7}{2}$, where ${\cal T}_{I}$ denotes the class of isosceles triangles. 
What is more, we have ${\rm Area} (T') \leq \frac{49}{16} {\rm Area}(C)$.

\begin{prob}
Can these estimates  be improved? If yes, what are the optimal constants? 
\end{prob}

Still another observation from  \cite{[L 02]}, p. 3079, it is that in a similar way we can show that every convex body $C$ contains a rectangular triangle $T$ such that a homothetic copy of it of ratio $\frac{7}{2}$ contains $C$.

\begin{prob} 
Is the bound $\tfrac{7}{2}$ optimal? If not, what is the optimal value? 
What can we say about the areas of rectangular triangles contained in $C$ (respectively, containing $C$)? 
Can these results be generalized to higher dimensional convex bodies and to ``rectangular simplices.''
\end{prob}

Approximation of convex bodies belonging to a given class by bodies from a different class is discussed at length in the survey \cite{Grub} by Gruber. 

\section{Banach-Mazur distances within special classes of normed spaces}  \label{sec:classes} 
In this section we will report on distances within special classes of convex bodies motivated by functional analysis. Accordingly, it will be more convenient to use the language of normed spaces. 

The first classical result, due to 
Gurari\u\i, Kadets and Matsaev \cite{GKM65}, gives precise (or asymptotically precise) bounds on Banach-Mazur distances between $\ell_p^n$-spaces. (See \cite{[Tom 78]} by Tomczak-Jaegermann for analogous results concerning noncommutative $L_p$-spaces, or {\em Schatten-von Neumann classes}.) 
\begin{thm} \label{GKM}
  If $1\leq p \leq q \leq 2\;$ or $\;2\leq p \leq q \leq \infty\;$, then $\;d (\ell_p^n, \ell_q^n) = n^{1/p -1/q}$.  
If $1\leq p \leq 2 \leq q \leq \infty$, then $c_1 n^\alpha \leq d (\ell_p^n, \ell_q^n) \leq  c_2 n^\alpha$, where $\alpha=\max \{1/p - 1/2, \, 1/2 -1/q\}$ and $c_1,c_2 >0$ are constants intependent of $p, q$, and $n$.
Moreover, we have $c_2 = 1$ if $n$ is a Hadamard number. 
\end{thm}
Recall that Hadamard numbers were defined in Section \ref{para-cross}. In the second statement we didn't use the more compact  notation $d (\ell_p^n, \ell_q^n) =  \Theta (n^\alpha)$ to emphasize that the constants $c_1,c_2$ are independent of all parameters involved and not just of the dimension $n$. 

Several proofs of Theorem \ref{GKM} are readily available in the literature (see, e.g., \S 37, \cite{[Tom 89]}), so we will just sketch several crucial tricks that illustrate methods for determining or estimating the BM-distance. 

First, the inequalities $\|\cdot\|_p\leq \|\cdot\|_q\leq n^{1/q-1/p}\|\cdot\|_p$ (for $p\leq q$) imply $\;d (\ell_p^n, \ell_q^n) \leq n^{|1/p -1/q|}$. Next, denote $\lambda=d(\ell_2^n,\ell_\infty^n)$ and let
$u \in \gGL(n)$ be such that $B_2^n \subset u(B_\infty^n) \subset \lambda B_2^n$. Then, for every $x \in \ell_\infty^n$ with $\|x\|_\infty =1$, we have $1\leq \|ux\|_2 \leq \lambda$. Consequently, if $e_1,e_2,\ldots,e_n$ is the standard unit vector basis and $x_k=ue_k$, then $\|x_k\|\geq 1$ and $\|\sum_{k=1}^n \ve_k x_k\|_2 \leq \lambda$ for any choice of  $\ve_k =\pm 1$. Now, by the (generalized) Parallelogram Law, 
\[
\lambda^2 \geq \tfrac 1{2^n} \|\sum_{k=1}^n \ve_k x_k\|_2^2 = 
\sum_{k=1}^n \|x_k\|_2^2 \geq n, 
\]
so $\lambda = d(\ell_2^n,\ell_\infty^n)\geq \sqrt{n}$ (and hence, by the preceding estimate, equal to $\sqrt{n}$). The first statement of of Theorem \ref{GKM} can now be deduced using the triangle inequality \eqref{triangle_ineq} and the duality relations \eqref{eq:duality},  \eqref{ellp-ellq}. Incidentally, analyzing the equality case in the above proof of $d(\ell_2^n,\ell_\infty^n)= \sqrt{n}$ shows that an operator $u \in \gGL(n)$ certifies that equality (i.e., satisfies $B_2^n \subset u(B_\infty^n) \subset \sqrt{n} B_2^n$) if and only if $u \in \gO(n)$ 
(cf. Mishra and Sain \cite{[MiSa]}). 

Concerning the second statement, the key instance is when $1/p+1/q =1$, 
in which case an upper bound $d (\ell_p^n, \ell_q^n) \leq  c_2 n^\alpha$ can be derived by interpolation using Hadamard matrices in a similar way as it was done in Section \ref{para-cross}, with $c_2=(1+\sqrt{2})^{2\alpha}$. 
Next, the lower bound $d(\ell_p^n,\ell_\infty^n)\geq \frac 1{\sqrt{2}} \sqrt{n}$ for $p\in [1,2]$ is obtained in a similar way as it was done above for $p=2$, but using the so-called {\em cotype $2$} property (see \cite{Maurey}) of $L_p$-spaces for $p\in [1,2]$ in place of the Parallelogram Law. The remaining instances can then be deduced using, again, using the triangle inequality \eqref{triangle_ineq} and the duality relations. 

The next special class that we are going to discuss is that of {\em permutationally-symmetric} convex bodies/normed spaces. A basis $(e_k)$ of a normed space $X$ (say, $\dim X=n$) is called permutationally-symmetric (usually abbreviated to ``symmetric'') if, for every sequence of scalars $(x_k)$, every choice of signs $\ve_k=\pm 1$, and every permutation $\pi \in S_n$, we have 
\begin{equation} \label{def:symmetric}
 \Big\|\sum_k\ve_k x_{\pi(k)}e_k\Big\|= \Big\|\sum_k x_ke_k\Big\|.
\end{equation} 
A normed space is {\em symmetric} if it has a symmetric basis, and a norm 
on $\bR^n$ is symmetric if the standard basis is symmetric with respect to that norm. 

Clearly all $\ell_p^n$ spaces are symmetric, but so are Orlicz and Lorentz sequence spaces. (See \cite{Byl} for some estimates concerning the latter.) In geometric terms, any convex body in $\bR^n$ that is invariant under permutations of coordinates and coordinate reflections is a unit ball with respect to a symmetric norm. Thus the class of symmetric spaces seems vastly larger than that of $\ell_p$-spaces, and so it comes as a surprise that their diameters with respect to the BM-distance are asymptotically the same. We have
\begin{thm}  \label{symmetric} 
Let $X, Y$ be symmetric normed spaces with $\dim X = \dim Y = n$. Then 
$d(X,Y)  \leq c\sqrt{n}$, where $c\geq 1$ is a universal constant. 
\end{thm}
The above theorem was proved by Tomczak-Jaegermann \cite{[Tom 83]} and was a culmination of steady progress over the preceding years, with slightly weaker bounds  shown in Gluskin \cite{[Gl 79]} ($O(\sqrt n \log^4 n)$; see also \cite{Gl 83}) and in Tomczak-Jaegermann  \cite{[Tom 79]} ($O(\sqrt n \log^2 n)$). 

A still larger class that of Banach spaces with unconditional basis, defined by replacing \eqref{def:symmetric} by 
\begin{equation} \label{def:unc}
 \Big\|\sum_k\ve_k x_ke_k\Big\|= \Big\|\sum_k x_ke_k\Big\|.
\end{equation} 
That is, in geometric terms, we are talking about convex bodies in $\bR^n$ that are invariant under coordinate reflections. While no conclusive answer concerning the diameter is known in this setting, major progress was achieved by Lindenstrauss and Szankowski in \cite{[LiSz86]}. 
\begin{thm} \label{unc} 
There is $\gamma \in [1/2,2/3]$ such that, for any
$\beta>\gamma$ there exists a constant $c(\beta)>0$ such that if $X,Y$
are $n$-dimensional normed spaces with unconditional bases, then 
$d(X,Y) \leq c(\beta) n^\beta$.
\end{thm}
The value of the constant $\gamma$ from Theorem \ref{unc} appears in \cite{[LiSz86]} as a solution to an optimization problem involving Lipschitz functions on $[0,1]$. The authors show that $\gamma \leq 2/3$ and state that there is numerical evidence (but no rigorous proof) that $\gamma \leq 5/9$. We thus have 
\begin{prob}
What is the optimal value of $\gamma$ in Theorem \ref{unc}? Is it true that, under the hypotheses of Theorem \ref{unc}, $d(X,Y) = O(\sqrt{n})$?
\end{prob}

The unconditional basis property \eqref{def:unc} is a relaxation of the symmetric basis property \eqref{def:symmetric}. Another relaxation is the following.
\begin{ddd} \label{es}
Let $X$ be a finite dimensional normed space and let $G$ be its group of (linear) isometries. We say that $X$ has {\em enough symmetries} if the only operators that commute with all elements of $G$ are multiples of $I_X$, the identity on $X$.
\end{ddd}
We can then ask
\begin{prob}
Determine the order of magnitude (as $n\to \infty$) of the diameter -- with respect to the BM-distance --  of the class of $n$-dimensional normed spaces with enough symmetries. In particular, is the diameter $\gg \sqrt{n}$?
\end{prob}
Clearly, the diameter is $\frac 32$ when $n=2$, at least $\frac 95$  when $n=3$, and at least $\sqrt{n}$ when $n\geq 4$. 

\bigskip 
We complete this section with a ``hybrid'' result due to Davis and Maurey  \cite{[DM 78]}, which involves two classes of spaces. 
\begin{thm} \label{symmetric-ellp} 
Let $X$ be a symmetric normed space with $\dim X= n$. Then 
\begin{equation} \label{eq:symmetric-ellp} 
    \min_{p \in [1,\infty]} d(X,\ell_p^n) = O(n^{1/4}).
\end{equation}
\end{thm}
This improved earlier result of Davis and Enflo \cite{[DE 77]} (see also \cite{Dav}), who established \eqref{eq:symmetric-ellp}  up to a polylogarithmic factor. 
The exponent $1/4$ is optimal, see, e.g., \S 49 in \cite{[Tom 89]}; it is also explained there why no improvement beyond $O(n^{1/2})$ is possible if we replace ``symmetric'' with ``unconditional.'' Similarly, it is not possible to approximate a general $n$-dimensional space within $o(\sqrt{n})$ by a space with unconditional basis (see \cite{FKP}), or even by a space possessing much weaker symmetries (see \cite{Sz83,Sz86}). See also \cite{CDGGJR} for a study of an ``unconditional-like'' property in the context of not-necessarily-symmetric bodies.

\vskip0.9cm 

\section{Other related concepts of distance} 

In this section we will discuss briefly other notions of distance that are affinely invariant, but differ in character from the Banach-Mazur type distances considered earlier.

 \subsection{Volume ratio of convex bodies} \label{sec:vr}

The first concept we consider is related to volume. 
\begin{ddd} \label{def:vr}
Let $K, L \subset E^n$ be convex bodies. 
The volume ratio of $K$ and $L$ is defined as 
\begin{equation} \label{eq:vr}
{\rm vr}(K,L) = \inf\Big(\tfrac{\vol(K)}{\vol(a(L))}\Big)^{1/n},
\end{equation} 
where the infimum is extended over all affine transformations $a$ for which $a(L) \subset K$. 
\end{ddd}
In the case when $L=B=B^n$, or an ellipsoid (in which case we abbreviate  
${\rm vr}(K,B) =: {\rm vr}(K)$, or ${\rm vr}(X)$ if $K=B_X$), the relevance of this concept to functional analysis was noticed in \cite{Sz78} (see also \cite{Kas}), where it was observed that, for all $n\in \bN$,  ${\rm vr}(B_1^n) \leq\sqrt{2\pi/e}$. A formal definition (in the language of normed spaces and for $L=B$) was given in \cite{SzTom}. 

The extremal values of ${\rm vr}(\cdot)$ and ${\rm vr}(B,\cdot)$ are very well understood. 
It was shown by Ball (\cite{Ball89, Ball91}) that 
\begin{thm} \label{vrB} 
We have 
 \begin{equation} \label{eq:vrB1} 
  \max_{M \in {\mathcal M}^n}  {\rm vr}(M) = {\rm vr}(P) \ 
  \hbox{ and } \ \max_{C \in {\mathcal C}^n}  {\rm vr}(C) = {\rm vr}(\Tt)
\end{equation}
and similarly
 \begin{equation} \label{eq:vrB2} 
  \max_{M \in {\mathcal M}^n}  {\rm vr}(B,M) = {\rm vr}(B,B_1^n) \ 
  \hbox{ and } \ \max_{C \in {\mathcal C}^n}  {\rm vr}(B,C) = {\rm vr}(\Tt).
\end{equation} 
In each of these settings, the extremal value of the volume ratio is $\Theta(\sqrt{n})$.  
\end{thm}

Note that, for a specific body $K$, the values of ${\rm vr}(K)$ and ${\rm vr}(B,K)$ can often be determined rather easily since optimality of the  affine transformation $a$ from \eqref{eq:vrB1} or \eqref{eq:vrB2} (the John position, or its dual analogue, the L\"owner position) can be certified via Theorem \ref{decomposition} or its version for the L\"owner ellipsoid. Next, for the concrete bodies $K=T, P, B_1^n$ appearing in \eqref{eq:vrB1} and \eqref{eq:vrB2} (and, more generally, for $B_p^n$, $1\leq p\leq \infty$), closed formulae for $\vol(K)$ (and, consequently, for ${\rm vr}(K)$ and ${\rm vr}(B,K)$) are well-known at least since Dirichlet \cite{Dirichlet}. Specifically, we have \begin{equation} \label{Bpn}
   \vol\big(B_p^n\big)=\frac{2^n\,\Gamma(1+1/p)^n}{\Gamma(n/p+1)}  
\end{equation}
(see, e.g., \cite{Gao} or Exercise 4.30 in \cite{Ver26} for accessible derivations). In particular, the quantities $\Theta(\sqrt{n})$ in Theorem \ref{vrB}  can be replaced by an explicit expression (typically, involving the {\em Gamma} function) or written in the form $\sim \alpha \sqrt{n}$ with an explicit constant $\alpha$.  

However, if neither of the bodies in question is an ellipsoid, determination of the optimal position is no longer straightforward. We will now discuss the values of ${\rm vr}(C,D)$, where one of the bodies $C,D$ is a parallelotope. This is the instance that has been of significant interest in the areas of computational geometry and mathematical optimization, see for example \cite{FS75,Toussaint83, ORurke85,HH23}. The relevant keywords that appear in the literature are the {\em minimum bounding box} (MBB) or the {\em minimum volume enclosing parallelotope} (MVEP). 

Early results in general dimension appeared in Macbeath \cite{Mac} and Hadwiger \cite {[Had]}, but we will start by discussing the case of dimensions $2$ and $3$. We will express the results and questions in terms of areas and volumes since this is how they normally appear in the literature.  

In dimension $2$, the first -- to the best of our knowledge -- results appeared around 1950 (Radziszewski \cite{[Radz1],[Radz2]}, Radziszewski and Bielecki \cite{[BiRa]}, S\"uss \cite{Suss}; see also the relatively recent paper \cite{Bri16} by Brinkhuis for a detailed exposition and comments) and the situation is completely clear. The most comprehensive is the result from \cite{[L 93]}, which was already stated as Theorem \ref{rectangles} in Section \ref{sec:rectangles}: {\em Given $C\in \mathcal{C}^2$, there exist parallelograms $P,P'$ such that $P \subset C \subset P'$ with ${\rm Area}(P) \geq \tfrac{1}{2}\cdot{\rm Area}(C)$ and ${\rm Area}(P')\leq 2\, {\rm Area}(C)$. Moreover, $P$ and $P'$ can be chosen to be homothetic rectangles (with the homothety ratio necessarily not exceeding $2$).} 

In dimension $3$, we have the following results.
\begin{thm} \label{vrP-3} 
Given $C\in \mathcal{C}^3$, there exist parallelepipeds $P,P'$ such that $P \subset C \subset P'$ with 
 \begin{equation} \label{eq:vrP-3i} 
  \vol(P) \geq  \tfrac 29 \vol(C) 
\end{equation}
and
 \begin{equation} \label{eq:vrP-3o} 
\vol(P') \leq 3\sqrt{2} \vol(C) 
\end{equation}
\end{thm}
The  inequality  \eqref{eq:vrP-3i} was established in \cite{[Radz2]},  \cite{[BiRa]} and \cite{Suss}, and is optimal (see Theorem \ref{PinT} below). The  inequality \eqref{eq:vrP-3o} is due to Makeev \cite{[Mak07]}, who (repeatedly) improved an earlier result of Radziszewski \cite{[Radz2]}, in which the upper bound was $6 \vol(C)$. 

To the best of our knowledge,  the optimal constant in \eqref{eq:vrP-3o} is unknown, which is somewhat surprising. It is fairly easy to see that it can not be smaller than $3$ (this is well-known, and a special case of Theorem \ref{TinP} below). We therefore ask 
\begin{prob}
What is the optimal constant in \eqref{eq:vrP-3o}? Is it equal to $3$? Equivalently, do we have $\max_{C\in {\mathcal{C}}^3}{\rm vr}(P,C) = 3^{1/3}$?
\end{prob}

The optimality of the bounds on ${\rm vr}(C,P)$ and ${\rm vr}(P,C)$  in dimension $2$ that are implicit in the results cited above is elementary; it was discussed already in Section \ref{sec:simplex} (see the discussion related to Figure \ref{fig:la25100a} and the comments following Theorem \ref{thm:T-center}).  The bound \eqref{eq:vrP-3i} is a little more delicate, but it is a consequence of the following general result from \cite{[L 99]} (announced in \cite{[L 99a]}; Leng, Zhang and Ma \cite{LenZhaMa} provided later an alternate proof).
\begin{thm} \label{PinT}
    The maximum volume of a parallelotope $P$ contained in an $n$-dimensional simplex $\Tt$ is equal to $\tfrac{n!}{n^n} {\rm}\vol(\Tt)$. Consequently, 
    ${\rm vr}(T,P)<e$ for all $n$ and, as $n\to \infty$, 
    ${\rm vr}(T,P)\sim e$.
\end{thm}
In view of the above results, one is tempted to guess that the simplex $T$ is the maximizer for  ${\rm vr}(\cdot,P)$. This is, however, very far from the truth, even for centrally symmetric sets. 
\begin{thm} \label{khrabrov}
    We have $\max_{M\in \mathcal{M}^n} {\rm vr}(M,P) = \Theta(\sqrt{n}) = \max_{C\in \mathcal{C}^n} {\rm vr}(C,P)$.
\end{thm}
The first equality can be found in Khrabrov \cite{[Kh02]}. Since the English translation of \cite{[Kh02]} may be difficult to navigate, here is a guide to the results. The more difficult $\Omega(\sqrt{n})$ part is an artifact of Gluskin's construction from \cite{[Gl 81]}; the statement can be found in \cite{[Kh02]} in the first Remark following the proof of Theorem 5. The $O(\sqrt{n})$ part is a Corollary to Remark 4 in \cite{[Kh02]}; it asserts in fact that ${\rm vr}(M,B_p^n)\leq \sqrt{en}$ for any $p\in [1,\infty]$  (and, of course, for any $n \in \bN$ and any $M\in \mathcal{M}^n$). The $O(\sqrt{n})$ part of the second equality follows easily from \eqref{eq:vrB1} and an explicit calculation (along the lines sketched in the paragraph following Theorem \ref{vrB}) yielding ${\rm vr}(P,B) \leq \sqrt{\frac{\pi e}2}$; the same argument shows in fact that ${\rm vr}(C,B_p^n) \leq \frac e2 \sqrt{n}$ for all $p\in [2,\infty]$. 

Concerning ${\rm vr}(P,C)$ in general dimension, a fundamental result is  due to Dvoretzky and Rogers \cite{DR50}, who showed that 
\begin{equation} \label{DRvol}
 \max_{M \in \mathcal{M}^n} {\rm vr}(P,M) = O(\sqrt{n}).
\end{equation}
Since,  by the Rogers-Shephard inequality \cite{RoSh1}, $\vol(C-C) \leq \binom{2n}{n} \vol(C)$ for any $C \in \mathcal{C}^n$, the same bound  holds for $\max_{C \in \mathcal{C}^n}$ (with the constant implicit in the $O(\cdot)$ notation multiplied by $4$). Of course, the bound in \eqref{DRvol}
is asymptotically optimal, as seen by considering $M=B$. The estimate from 
from \cite{DR50} was subsequently slightly improved in \cite{PeSz} and \cite{FNZ}, but arguably the most elegant statement in this direction is due to Ball \cite{Ball89}. 
\begin{thm} \label{Ball}
If $M\in \mathcal{M}^n$, then 
\begin{equation} \label{eq:Ball}
    {\rm vr}(P) \leq {\rm vr}(M)\cdot {\rm vr}(P,M) \leq \sqrt{e} \ {\rm vr}(P).
\end{equation} 
\end{thm}
Since, trivially, ${\rm vr}(P)\leq \sqrt{n}$ and ${\rm vr}(M)\geq 1$, \eqref{eq:Ball} implies \eqref{DRvol}. 

The final result belonging to this circle of ideas that we mention concerns ${\rm vr}(P,T)$. It is more convenient to state it in terms of volume. 
\begin{thm} \label{TinP}
If $T$ is an $n$-dimensional simplex contained in a parallelotope $P$, 
then 
\begin{equation} \label{eq:TinP}
    \vol(T)/\vol(P) \leq \tfrac{(n+1)^{(n+1)/2}}{2^n n!}=:\delta(n). 
\end{equation}
This inequality can be saturated if and only if $n+1$ is a Hadamard number. 
\end{thm}

Theorem \ref{TinP} goes back essentially to Hadamard \cite{Hadamard}. The argument is quite simple and well-known, but we will sketch it here since it has bearing on other -- less well-known -- results in this survey. 
The trick is to connect the question at hand to two determinant-optimization problems: 
\begin{equation}\label{det01}
f_n = \max \{ |\det(A)| : A \hbox{ is an } n\times n  \hbox{ matrix with } 0,1  \hbox{ entries}\}, 
\end{equation}
\begin{equation}\label{det1-1}
g_m = \max \{ |\det(H)| : H \hbox{ is an } m\times m  \hbox{ matrix with } \pm 1  \hbox{ entries}\}. 
\end{equation}
Going back to \eqref{eq:TinP}, we can assume that $P$ is the unit cube, $[0,1]^n$. While it is possible that not all of the vertices of a maximum volume simplex contained in $P$ are vertices of the cube, there is always a maximizer that has this property (see Figure \ref{fig:la25100a}). Further, by symmetry, we can assume that one of those vertices is the origin. It follows that 
$\delta(n) = \frac 1{n!} f_n$. 

Next, it is fairly easy to see that 
\begin{equation} \label{fg}
  g_{n+1}= 2^n f_n.   
\end{equation}
Indeed, given $(n+1)\times (n+1)$ matrix $H$ with $\pm 1$ entries, we can successively perform standard row operations (which do not change the absolute value of the determinant) as follows: \\
$1^\circ$ where needed, we multiply columns by $-1$ so that the first row consist exclusively of $1$'s\\
$2^\circ$ where needed, we multiply rows $2$ through $n+1$ by $-1$ so that the first entry of each such row is $-1$\\
$3^\circ$ we add the first row to rows $2$ through $n+1$; every entry of each such row is then $0$ or $2$, with the first entry being $0$.

After these operations, the first column of the obtained matrix $H'$ is $[1,0,\ldots,0]^t$, and the $n\times n$ matrix obtained by deleting the first row and the first column of $H'$ is of the form $2A$, where $A$ is a matrix of the type appearing in \eqref{det01}. In particular, $|\det(H)| = |\det(H')| = 2^n |\det (A)|$. Since this procedure is clearly reversible, \eqref{fg} follows. 

Finally, we note that, by Hadamard's inequality, if $B$ is an $m\times m$ matrix with $\pm 1$ entries, then $|\det(H)| \leq m^{m/2}$, with equality iff $H$ is orthogonal, which -- by definition -- can happen exactly when $m$ is a Hadamard number. Combining the above observations we deduce Theorem \ref{TinP}.  

\begin{prob}
Find an upper bound for ${\rm vr}(P,T)$  when $n+1$ is not a Hadamard number, which nearly matches the $\delta(n)^{-1/n}$, for example in the spirit of Theorem \ref{thm:Paley}. At the very least, show that ${\rm vr}(P,T) \sim \delta(n)^{-1/n}$.
\end{prob}
Note that it follows from Stirling's formula that $\delta(n)^{-1/n} \sim \frac 2e \sqrt{n}$. 

\begin{prob} \label{prob:vrPT}
    Is ${\rm vr}(P,T) = \max_{C\in \mathcal{C}^n} {\rm vr}(P,C)$? 
\end{prob}

The proof of Theorem \ref{TinP} yields additional information, which allows to deduce Theorem \ref{simplex-H}. Namely, it follows from the construction that if $n+1$ is a Hadamard number, $H$ - an $(n+1)\times (n+1)$ Hadamard matrix, and $T \subset P=[0,1]^n$ - an $n$-dimensional simplex (of maximum volume) built from $H$, then $T$ is regular, inscribed into the Euclidean ball $B$ (of radius $\sqrt{n}/2$) circumscribed about $P$, and its centroid $t$ coincides with the centers of $P$ and $B$. Since the homothethy $h$ with center $t$ and ratio $n$ clearly verifies $h(T)\supset B$, it follows that 
$d_{\rm BM}^{ctr} (T,P) \leq n$ and   $d_{\rm BM}^{\rm cen} (T,P) \leq n$. The reverse inequality is trivial since both  $d_{\rm BM}^{ctr}$ and $d_{\rm BM}^{\rm cen}$ {\em always} dominate $d_{\rm BM}$ and $d_{\rm BM} (T,P)=n$ by Theorem \ref{thm:GLMP}. This proves the instance $M=P$ of Theorem \ref{simplex-H}, and the case of the cross-polytope follows by duality, which works with no issues in this configuration. 

It is not immediately clear what the above argument allows to deduce when $n+1$ is not a Hadamard number. Accordingly, we ask

\begin{prob} \label{simplex-nonH}
    If $M$ is either the parallelotope or the cross-polytope, derive reasonable estimates for $d_{BM}^{ctr} (T,M)$ and $d_{\rm BM}^{\rm cen} (T,M)$ for the dimensions $n$ such that $n+1$ is not a Hadamard number, possibly appealing to the relative ubiquity of Hadamard numbers in the spirit of Theorem \ref{thm:Paley}.
\end{prob}

In the case of general bodies, a pretty sharp upper bound on the volume ratios 
was obtained very recently by Galicer, Merzbacher and Pinasco in \cite{GMP26}. 
\begin{thm} \label{thm:vr} 
We have \ $\sup\{{\rm vr}(K,L) :  K, L \in {\mathcal C}^n\} =O\big(\sqrt {n  \log (1+n)}\,\big)$. 
\end{thm}
The above bound improved on an earlier estimate $O(\sqrt n  \log n)$ shown by Giannopoulos and Hartzoulaki \cite{[Gia-Har]} and -- in view of Theorem \ref{vrB} -- is optimal up to the $\sqrt{\log (1+n)}$ factor. 
\begin{prob}
Can we replace the bound in Theorem \ref{thm:vr} by $O(\sqrt{n})$? 
\end{prob}
Further invariants related to the volume ratio have been defined and studied, among others, by Galicer, Merzbacher and Pinasco \cite{GMP21}, Khrabrov \cite{[Kh02]}, and Markesinis \cite{Mar}. 

While ${\rm vr}(\cdot,\cdot)$ does quantify similarity between convex bodies, it doesn't have all features of ``distance,'' in particular it is highly non-symmetric. (For example, as stated above, we have ${\rm vr}(B_1^n,B^n)= \Theta(1)$ while ${\rm vr}(B^n,B_1^n)=\Theta(\sqrt{n})$.)  Here is one natural remedy to this drawback.
\begin{ddd}\label{def:dvol}
Given $K, L \in {\mathcal C}^n$, let us define their {\em volumetric distance} as
\begin{equation} \label{eq:dvol}
d_{\rm vol}(K,L) := {\rm vr}(K,L)\cdot {\rm vr}(L,K) = \inf\Big\{\Big(\tfrac{\vol(b(L)}{\vol(a(L))}\Big)^{1/n} : a(L) \subset K \subset b(L)\Big\},
\end{equation}
where $a, b$ are affine maps. 
\end{ddd}
A version of $d_{\rm vol}$ was defined by Macbeath in \cite{Mac} and subsequently investigated by various authors, e.g., Khrabrov \cite{[Kh01a],[Kh01b],[Kh01c],Kh04}, and implicitly Gluskin \cite{[Gl 81]}.

Since the affine images $a(L)$ and $b(L)$ in \eqref{eq:dvol} are not required to be homothetic, 
this is a relaxation of the BM-distance and  of the Gr\"unbaum distance $d_{BM}^\pm$. Accordingly
\begin{equation}  \label{ineq:dvol}
d_{\rm vol} \leq d_{BM}^\pm \leq d_{BM}.
\end{equation} 
It easily follows that
\begin{thm} \label{thm:dvol-nonsym}
We have $\max_{K, L \in {\mathcal C}^n} d_{\rm vol}(K,L) = n$.
\end{thm}
Indeed, the inequality ``$\leq$'' follows by combining \eqref{ineq:dvol} with \eqref{Gru-n1}, and it is saturated by the pair $(K,L)=(\Tt,B)$. 

The bound from Theorem \ref{thm:dvol-nonsym} is asymptotically correct also for centrally symmetric bodies, but is much more difficult to establish in that context.  It is a corollary of the proof of Theorem \ref{gluskin}. 
\begin{cor} \label{cor:dvol-sym}
We have $\max_{K, L \in {\mathcal M}^n} d_{\rm vol}(K,L) = \Theta(n)$.
\end{cor}
In the setting of the Corollary, the affine maps $a,b$ from \eqref{eq:dvol} can be replaced by linear operators. The reader will notice that in the outline of the proof of Theorem \ref{gluskin}, which we sketched in Section \ref{radius:ngeneral}, the failure (with probability close to $1$) of the inclusion $vM_1 \subset c_1\sqrt{n} M_2$  for {\em all} $v \in \gSL(n)$ 
and similar failure for $M_2 \subset c_1\sqrt{n}\, wM_1$ are considered separately, without requiring that the linear operators $v$ and $w$ be related in any way. 

Based on the above, the reader may suspect that there is a close relationship between $d_{\rm vol}$ and $d_{BM}^\pm$. In general, this is not the case.  Consider, for example $M=B^k \times P^m$, where $k+m=n$ and $m\approx (n/\log n)$. It is easy to show that then $d_{\rm vol}(M,B)= \Theta(1)$, while  $d_{BM}^\pm(M,B)= d_{BM}(M,B)=\sqrt{m}$. On the other hand, if $M$ is the unit ball of a normed space with enough symmetries, then it is easy to see that $d_{\rm vol}(M,B)= d_{BM}(M,B)$. This suggests the following. 
\begin{prob}
If $M_1,M_2$ are unit balls of symmetric normed spaces (or, more generally, spaces with enough symmetries)  and
$d_{\rm vol}(M_1,M_2)$ is small, does it follow that $d_{BM}(M_1,M_2)$ is small?  
\end{prob}
We did not give a lot of thought to this question; it is possible that there is a simple counterexample. However, a modest evidence in support of the affirmative answer is that if $M_1=B_p^n$ and $M_2=B_q^n$, then $d_{BM}(M_1,M_2) = \Theta\big(d_{\rm vol}(M_1,M_2)\big)$,  see \cite{[Kh02]}.

\subsection{The weak distance of normed spaces} \label{sec:weak}
 
The {\em weak distance} of normed spaces, sometimes referred to as the {\em weak Banach-Mazur distance}, was defined by Tomczak-Jaegermann in early 1980's \cite{Tom 81,[Tom 84]}. In a sense, it is a {\em convexification} of the BM-distance. We first reformulate the definition \eqref{def:dBMns} of $d(\cdot,\cdot)$ as follows
\begin{equation}  \label{def:dBMns2}
 d(X,Y) = \inf \big\{ \|u\|\cdot \|v\| : u \circ v = I_X \big\}, 
\end{equation}
where $v: X\to Y$ and $u: Y\to X$. We then define an auxiliary (non-symmetric) quantity 
\begin{equation} \label{def:qweak}
   q(X,Y) := \inf \Big\{ \sum_k \|u_k\|\cdot \|v_k\| : \sum_k u_k \circ v_k = I_X\Big\}, 
\end{equation}
where $v_k: X\to Y$ and $u_k: Y\to X$, $k=1,2,\ldots N$, are finite sequences of operators ({\em a priori} of arbitrary length\footnote{An elementary argument using Carath\'eodory's theorem shows that it is enough to consider $N\leq n^4$, and that the infumum is attained.}), and finally the weak distance between $X$ and $Y$ as
\begin{equation} \label{def:weak}
d_w(X,Y) = \max \{q(X,Y), q(Y,X)\}. 
\end{equation} 
Obviously, 
\begin{equation} \label{weak-BM}
    d_w(X,Y) \leq d(X,Y). 
\end{equation} 
The articles \cite{Tom 81,[Tom 84]} initiated the study of this concept by proving a number of fundamental results and stating a number of open problems. It is also noted in Proposition 3 in \cite{[Tom 84]} that $d_w$  makes sense even if $\dim X \neq \dim Y$ and that it can be alternatively defined using the language of the so-called {\em  ideal norms}. We will not pursue these directions systematically, but we note that appropriately chosen ideal norms provide often -- in combination with \eqref{weak-BM} --  a very effective way of supplying lower bounds for $d(X,Y)$. 

A few years later, Lindenstrauss and Szankowski 
 \cite{[LiSz87], LiSz87a} developed a machinery that -- in principle -- allows to calculate the weak distance between two symmetric normed spaces 
 up to a polylogarithmic factor.  Next, important contributions were made by Rudelson in \cite{[Rud 95]}, and more calculations of the weak distance for (classes of) specific pairs of normed spaces were done in Khrabrov \cite{Kh04,Kh05} and Mankiewicz \cite{Mank12}. 

All the above papers addressed various aspects of two natural questions about the weak distance of $n$-dimensional normed spaces:\\
$1^\circ$ Given two such spaces $X,Y$, how much can $d_w(X,Y)$  differ from $d(X,Y)$?\\
$2^\circ$  What is the maximal weak distance between two such spaces? \\
Here is an overview of some of the progress that was achieved in these directions.

First, it was observed in \cite{Tom 81} that, for any $X$,  $d_w(X,\ell_2^n)=d(X,\ell_2^n)$. Next, it followed from the methods developed in  \cite{[LiSz87]} that, for some symmetric spaces $X,Y$, 
the ratio $d(X,Y)/d_w(X,Y)= \Omega(n^\alpha)$, where $\alpha$ (determined numerically) is approximately $0.02$. For general $X,Y$, it was shown in \cite{[Tom 84]} that for ``the usual suspects,'' the Gluskin-type spaces described in Section \ref{radius:ngeneral}, we    typically have $d_w(X,Y)= O(\sqrt{n})$ while, as determined in the proof of Theorem \ref{gluskin}, the bound $d(X,Y)= \Omega(n)$ holds generically. The same technique was refined in Rudelson \cite{[Rud 95]} to produce examples for which $d_w(X,Y)=O(1)$ (or even, with some caveats,  $1+o(1)$),  while $d(X,Y)=\Omega(n)$. In the same paper, in the direction of the query $2^\circ$ above, it was shown that 
 \begin{equation}
d(X,Y) = O\big(n^{13/14} \log^{15/7}(1+n) \big). 
 \end{equation}
 This leaves as the most important open problem in this area
 \begin{prob} \label{prob:weak1}
What is the order of magnitude of $\max\{d_w(X,Y) : \dim X = \dim Y =n\}$? Is it $O(\sqrt{n})$, perhaps up to a polylogarithmic factor? 
 \end{prob}

Another open problem, versions of which were asked in \cite{[Tom 84]} and \cite{LiSz87a}, asks the following. 
 \begin{prob} \label{prob:weak2}
Is it true that $d(X,\ell_p^n)/d_w(X,\ell_p^n)= O(1)$, perhaps up to a polylogarithmic factor? What if we additionally assume that $X$ has a symmetric basis? 
 \end{prob}

\subsection{Cones and projective transformations} 
Apart from convex bodies, another class of sets that is important in convex analysis is that of (closed) {\em convex cones.} It is easy to see that every nondegenerate (i.e., pointed and full-dimensional) closed convex cone $C \subset \bR^{n+1} =\bR \times \bR^n$ can be, up to affine equivalence, represented as follows
\begin{equation} \label{cone-base}
C = \{ t \hskip0.3mm x \ : \ t\geq 0, \, x\in K\},    
\end{equation}
where $K$ is a convex body in the hyperplane $E^n=\{1\} \times \bR^n$. (It will be convenient to denote the standard basis of $\bR^{n+1}$ by $(e_0,e_1,\ldots,e_n)$, so that $E^n=\{x \in \bR^{n+1} : \langle x,e_0\rangle = 1\}$.) 
In that case, we will refer to $K$ as a {\em base} of $C$. 
Clearly, affinely equivalent bases lead to affinely equivalent cones, but the converse is false. For example, it is easy to see that all cones generated by convex quadrangles are affinely equivalent, while obviously not all quadrangles are. In this context , the correct framework for analyzing convex bodies in $E^n$ is that of {\em projective transformations} (see, e.g., the textbook \cite{Kaplansky} by Kaplansky, which focuses on $n=2$, with the general case left as an exercise).  Ignoring the ``points at infinity'' (which we can if we only care about cones of the form \eqref{cone-base}), projective transformations can be represented as follows. Given $u\in \gGL(n+1)$ with $\langle ue_0,e_0\rangle =1$ and $x \in E^n$, we set
\begin{equation}
    \Phi (x) = \tfrac {ux}{\langle ux,e_0\rangle}, 
\end{equation}
which is defined (and belongs to $E^n$) except for the points being sent by $\Phi$ to infinity. In this language, we define the {\em projective Banach-Mazur distance} of $C, D \in {\mathcal C}^n$ as follows. 
\begin{equation} \label{eq:proj}
d_{\pi} (C,D) = \inf_{\Phi,\Psi, h_{\lambda}} \{ \lambda : \, \Phi(D) \subset \Psi(C) \subset h_{\lambda}(\Phi(D)) \}, 
\end{equation} where $\Phi,\Psi$ are projective transformation with $\Phi(D) \subset E^n$ and $h_{\lambda}$ is a homothety with a {\sl positive}  ratio $\lambda$. 

Obviously, $d_\pi \leq d_{BM}$. The formula \eqref{eq:proj} is more complicated that the definition \eqref{eq:dBM} of $d_{BM}$ : it involves two projective transformations $\Phi,\Psi$ as opposed to one affine map $a$ (this is because $\Phi^{-1} h_\lambda \Phi$ is typically not a homothety). However, while the relevance of projective transformations to the theory of cones was well-known, the invariant $d_{\pi}$ hasn't been -- to our knowledge -- studied in depth, if at all, and so we are not sure to what extent this complication is necessary. In particular, we do not know answers to the following simple-minded questions. (Prove or give a counterexample.)
\begin{prob} 
    If $C, D \in {\mathcal M}^n$, is $d_{\pi} (C,D)= d_{BM} (C,D)$? 
\end{prob}
\begin{prob}
    For $C, D \in {\mathcal C}^n$, is the infimum in \eqref{eq:proj} always attained? 
\end{prob}

Finally, it is of some interest to to determine the diameters and radii of the Bamach-Mazur compacta with respect to $d_{\pi}$, for example
\begin{prob} \label{proj-diam}
 What is $\sup_{C, D \in {\mathcal C}^2}d_{\pi} (C,D)$? What are the extremal sets (assuming the supremum is attained)? What about the same questions for higher dimensions? 
\end{prob}
It appears likely that the asymptotic behavior of the diameters/radii with respect to $d_{\pi}$ is similar to that for $d_{BM}$, so Problem \ref{proj-diam} is probably most interesting in low dimensions. Also, note that since projective images of a simplex or an ellipsoid are again a simplex and an ellipsoid, some of the questions about $d_{\pi}$ have simple answers. 

Even more exotic analogues of the Banach-Mazur distance appear in areas such as symplectic geometry (for example, symplectic Banach-Mazur distance or contact Banach-Mazur distance, see, e.g., \cite{PRSZ,Melistas, SZ,DGRZ}). However, the context in which they normally arise includes non-convex sets, and so we will not focus on that direction of research in this survey.

\vskip0.4cm

\section{\bf BM-distance between sections of convex bodies} \label{sec:sections}

In this section we present a sample of results addressing the relationship between BM-distance of convex bodies and that of their sections (resp., normed spaces and their subspaces). Note that in view of the duality between sections and projections (the polar of the section $K\cap E$ of a body $K$ is the projection of $K^\circ$ onto $E$;  see, e.g., section 1.1.4 in \cite{AuSz17}), statements  about subspaces (or sections) typically have dual companion formulations involving quotient spaces (or projections, or affine images), with some important caveats in the not-necessarily-symmetric case that were already pointed out in Section \ref{sec:intro} (see the comments following \eqref{eq:polar}). We will concentrate on the first point of view and most results will be stated in the language of normed spaces. In the setting of bodies from $\mathcal{C}^n$, we will restrict ourselves to a few positive results and examples that indicate possible difficulties.

\vskip 0.09cm
We start with the following result shown in  \cite{[L 08]}.

\begin{thm} \label{$(2c+1)^2$}
Let $X,Y$ be $n$-dimensional normed spaces and let $X_0 \subset X, Y_0\subset Y$ be subspaces of codimension $k$. Then \\
{\rm (i)} $d(X_0,Y_0) \leq (2k+1)^2 d(X,Y)$\\
{\rm (ii)} $d(X,Y) \leq (2k+1)^2 d(X_0,Y_0)$. 
\end{thm}
Technically,  \cite{[L 08]} contains only a statement of part (i) with $X=Y$. However, the proof is based on a result of type (ii) (note that if $X_0,Y_0 \subset X$ are both subspaces of codimension $k$, then they contain a common subspace $Z$, whose codimension in both $X_0$ and $Y_0$ is at most $k$, and so part (ii) with $X_0=Y_0=Z$ applies). Next, the generalization to $X\neq Y$ (resp., $X_0\neq Y_0$) is straightforward. 

In view of Corollary \ref{john2}, the bounds in Theorem \ref{$(2c+1)^2$} are nontrivial only if $k$ is much smaller than $n$ (at most $O(\sqrt{n})$, for starters).  It is also unlikely that the factor $(2k+1)^2$ in Theorem \ref{$(2c+1)^2$} is optimal. Accordingly we ask
\begin{prob}
 What is the optimal factor in   Theorem \ref{$(2c+1)^2$}? In particular, is it $O(k)$? 
\end{prob}
Clearly, the optimal factor can not be $o(k)$. Indeed, if $n=2k$ and $X=Y=X_0\oplus Y_0$, where $X_0, Y_0$ are arbitrary $k$-dimensional spaces, the best bound we can have for $d(X_0,Y_0)$ is $\Theta(n)$ (by Theorem \ref{gluskin}). The author of \cite{[L 08]} muses that the optimal value could be $2k+1$. This is mildly supported by the fact that, in the infinite-dimensional setting, it is known that there is a space $X$ and two subspaces $X_0, Y_0$ of codimension $1$ such that $d(X_0,Y_0) =3$, see Cambern \cite{Cambern}.  Specifically, this happens if $X=c$, the space of convergent sequences with the supremum norm, and $X_0$ and $Y_0$ are subspaces isometric to, respectively, $c$ and $c_0$, the subspace consisting of sequences converging to $0$. 

Of course, Theorem \ref{$(2c+1)^2$} can be formulated equivalently in the language of centrally symmetric convex bodies and their central sections. The symmetry hypothesis is important, but not entirely indispensable. More detailed analysis of related issues can be found in \cite{Szar14}; here we restrict ourselves to presenting one simple positive result and an example hinting at the intricacies of the not-necessarily-symmetric setting. 

Let $K\subset E^n$ be a convex body and $H\subset E^n$ a hyperplane. We will call the section $H\cap K$ of $K$ {\em balanced} if $H$ is located exactly halfway between the two hyperplanes parallel to it and supporting to $K$. 

\begin{pro} \label{section-nonsym}
Let $K, L \in \mathcal{C}^3$ and let $K_0, L_0$ be their balanced sections. Then
\begin{equation} \label{eq:section-nonsym}
d_{BM}(K,L) \leq 9 d_{BM}(K_0,L_0).     
\end{equation} 
\end{pro}
Here is an outline of the proof, which is similar to, but simpler than,  the argument from \cite{[L 08]}. By replacing $K$ with an affine image, we may assume that $K\subset \bR^n$, $H =\{x \in \bR^n : x_n=0\}$, $\pm e_n \in \partial K$, and $K \subset \{x \in \bR^n : |x_n|\leq 1\}$. It is then an elementary exercise to show that if $K'= \frac 23 K_0 \times [-\frac 13,\frac 13] \subset \bR^{n-1} \times \bR \sim \bR^n$, then $K' \subset K \subset 3K'$. Since, for a similarly constructed $L'$, we clearly have $d_{BM}(K',L') \leq d_{BM}(K_0,L_0)$, the conclusion \eqref{eq:section-nonsym} follows from the triangle inequality \eqref{triangle_ineq}.  
Proposition \ref{section-nonsym} is an analogue of assertion (ii) of Theorem \ref{$(2c+1)^2$} for $k=1$. 
\begin{prob}
 Do we have an analogue of assertion (i) from Theorem \ref{$(2c+1)^2$} in the spirit of Proposition \ref{section-nonsym}? Is there a reasonable generalization to $k>1$?
\end{prob}
Note that, in the non-symmetric setting, a statement of type (ii) doesn't formally imply (i) since, even if $K_i = K\cap H_i$ for $i=1,2$ are balanced sections of $K$, it doesn't follow that $K \cap (H_1\cap H_2)$ is a balanced section of $K_1$ or $K_2$. 

While Proposition \ref{section-nonsym} can be generalized to the situation where the sections are ``reasonably'' balanced, it is quite clear that without any restriction on the location of the sections, no non-trivial estimates in the spirit of  Theorem \ref{$(2c+1)^2$} or Proposition  \ref{section-nonsym} are possible. For starters, two hyperlane sections of the same body can be of a different dimension. 
For a slightly less trivial example, consider $n>2$ and  two arbitrary bodies $K_0, K_1 \subset \bR^{n-1}$, which we may assume to verify $B \subset K_s \subset (n-1)B$. Next, let $K_s'$ be an affine image of $K_s$ contained in the hyperplane $H_s :=\{x \in \bR^n : x_n = (-1)^s\}$ and let $K= \conv(K_0'\cup K_1'\cup\{\pm te_n\})$, where $t$ is just slightly larger than $1$ (e.g, $t=\frac n{n-2}$ will work). Then $K\cap H_s = K_s'$ for $s=0,1$; in other words, $(n-1)$-dimensional sections of $K$ can be as different as it is possible. 

There are reasons to think that considering projections of convex bodies (or, equivalently, affine images of smaller dimension) may be more appropriate that considering sections. For example, if $K\in \mathcal{C}^n$ and $K'$ is an affine image of dimension $k$, then $d_{BM}(K', B^k)\leq d_{BM}(K,B^n)$, while a similar statement for sections is clearly false in view of the example from the preceding paragraph.  It would be of interest to investigate whether this approach is useful, but -- to the best of our knowledge --  this direction of research did not seems to attract much attention. Here is a vague question to encourage such exploration.
\begin{prob} \label{prob:nonsym-proj} 
Are there reasonable generalizations of statements in the spirit of Theorem \ref{$(2c+1)^2$} to bodies $K\in \mathcal{C}^n$ if we replace subspaces/sections of codimension $k$ by affine images of dimension $n-k$?
\end{prob}

\medskip 
The next few theorems demonstrate further limits for what can be inferred about the body from the knowledge of its sections, even in the centrally symmetric/normed space setting. The first result is due to Rudelson \cite{[Rud 04]}. He estimates extremal distances between $k$-dimensional sections of two centrally symmetric convex bodies in $E^n$.
For a positive integer $k < n$ the symbol $d_k(K,D)$ denotes the smallest Banach-Mazur distance between $k$-dimensional sections of $K$ and $D$.
By $\Delta(k,n)$ the author means the supremum of $d_k(K,D)$ over all centrally symmetric convex bodies $K$ and $D$ in $E^n$.
He proves the following theorem

\begin{thm} \label{Rudelson}
For every $k < n$ we have $\Delta(k,n) \sim_{\log n} \sqrt k$ if $k \leq n^\frac{2}{3}$ and $\Delta(k,n) \sim_{\log n} \frac{k^2}{n}$ if $k> n^\frac{2}{3}$.
Here $A \sim_{\log n} B$ means that $c/\log^a n \cdot A \leq B \leq c^{-1}\log^a n \cdot A\;$ for some absolute positive constants $C$ and $a$.
\end{thm}

\vskip 0.09cm 
The preeminent result concerning sections of convex bodies is arguably Dvoretzky's theorem (\cite{Dvoretzky}). 
\begin{thm} \label{dvoretzky}
    Given $\ve >0$, there exists a sequence $k_n \to \infty$ such that every $K\in \mathcal{C}^n$ contains a $k_n$-dimensional section $K'$ verifying $d_{BM}(K',B) \leq 1+\ve$. Moreover, if $K$ is in the John position, then for most $k_n$-dimensional vector subspaces $E\subset \bR^n$ we have 
    $r(B\cap E) \subset K\cap E \subset (1+\ve)r(B\cap E)$, for some appropriate scaling constant $r$. 
\end{thm}

Dvoretzky's theorem is usually stated in the language of normed spaces and their subspaces, but most proofs yield the more general version stated above. The ``John position'' hypothesis can likewise be substantially relaxed. Simpler proofs were subsequently given in Milman \cite{[Mi71]}, Szankowski \cite{Szan}, and Figiel \cite{Figiel75}. A very influential work addressing related circle of ideas was the paper \cite{FLM} by Figiel, Lindenstrauss and Milman.

The original proof of Theorem \ref{dvoretzky} gave $k_n =\Omega(\sqrt{\log n})$. This was improved to the asymptotically optimal $\Omega(\log n)$ in the seminal paper \cite{[Mi71]}. Additionally, \cite{[Mi71]} supplies a ``formula'' for the radius $r$ of the ``almost spherical'' sections, and consequently it is sometimes referred to as a ``tangible'' version of Dvoretzky's theorem. The optimal dependence on $\ve$ is not known; the method of \cite{[Mi71]} gave $\Omega(\ve^2)$, improved to $\Omega\big(\ve/(\log (1+1/\ve))^2\big)$ by Schechtman \cite{Sch13}. In this context, we would like to advertise the following sample question, an affirmative answer to which would go a long way towards addressing the ``dependence on $\ve$'' issue (cf. Milman \cite{Mi88}). 
\begin{prob} \label{knaster} 
    Given continuous function $f : S^{n-1} \to \bR$, $m \leq \sqrt{n}$, and a collection of points $x_1, x_2,\ldots, x_m \in S^{n-1}$, does there exist an isometry $U \in \gO(n)$ such that 
    $f(Ux_1)=f(Ux_2)=\ldots =f(Ux_m)$?
\end{prob}
The above question is a weakened variant of the well-known Knaster problem \cite{Knaster47}, which asked whether the answer is positive if $m=n$, and which was solved in the negative (for sufficiently large dimensions) in \cite{KaSz03} (see also \cite{HiRi05}). 

As an example of a quantitative version of Theorem \ref{dvoretzky} (i.e., potentially large $\ve$) we mention the  following result of Milman and Schechtman \cite{MS95, MS99}, and  Gordon, Guedon and Meyer \cite{GGM}.

\begin{thm} \label{GeedonMeyerPajor}
Given $n\in \bN$ and $\log n \leq k \leq n$, every $n$-dimensional normed space has a $k$-dimensional subspace $E$ such that $d(E, \ell_2^k)  =  O (\sqrt{k / (1+ \log (n /k)})$. Similar estimate holds for sections of not-necessarily symmetric convex bodies provided $k\leq cn/(1+\log n)$ for appropriate $c>0$. 
\end{thm}

We conclude this section with two spectacular results showing presence of Euclidean-like structure of very high dimension under surprisingly weak assumptions.  The first result appeals to the concept of volume ratio introduced in Section \ref{sec:vr}. 
\begin{thm} \label{vr-kashin}
    Let $X$ be an dimensional normed space verifying ${\rm vr}(X) \leq v$ and let $\theta\in (0,1)$. Then there exists a subspace $Y\subset X$ with $\dim Y = k \geq \theta n$ and $d(Y, \ell_2^k \leq \lambda$, where $\lambda$ depends only on $v$ and $\theta$. If $X$ is in the John position, the above holds for most $k$-dimensional vector subspaces of $X$. 
\end{thm}
Theorem \ref{vr-kashin} as stated was shown in \cite{Sz78, SzTom}, but it has roots in \cite{Kas}, which was focused on the instance $X=\ell_1^n$ and the following phenomenon implied by the Theorem, often referred to as {\em Kashin decomposition}. 
\begin{cor}\label{kashin}
    Given $k\in \bN$, there exist two orthogonal $k$-dimensional subspaces $E_1,E_2$ of $\ell_2^{2k}$ such that 
    \begin{equation}\label{eq:kashin}
     c \|x\|_2  \leq \tfrac{\|x\|_1}{\sqrt{n}} \leq \|x\|_2 
    \end{equation}
  for all $\ x\in E_i, \ i=1,2$,  where $c>0$ is a universal constant. 
\end{cor}
Of course, similar assertion holds (for example) for three orthogonal $k$-dimensional subspaces $E_1,E_2,E_3$ with \eqref{eq:kashin} holding for all 
$x \in E_i +E_j$, for any combination of  $i,j \in \{1,2,3\}$. This particular fact was exploited, e.g.,  in \cite{Pis86} or \cite{KaSz03}. 

The second result is even more spectacular as it holds without any additional hypotheses. 

\begin{thm} [\cite{Mi85}] \label{MilmanQS}
  Given  $\theta \in (0,1)$,  there exists a constant $c=c(\theta) \geq 1$ such that, for any $n$-dimensional normed space $X$, some subspace $Z$ of a quotient of $X$, say $Z$, verifies $\dim Z =:k \geq \theta n$ and $d(Z,\ell_2^k)\leq c$. Equivalently, for any $M \in \mathcal{M}^n$ there is an affine subspace $F$ passing through the center of $M$ and an affine map $a : F\to E^k$ such that $d_{BM}(a(M\cap F),B^k \leq c$. 
    \end{thm}
    
The literature addressing the circle of ideas related to Dvoretzky's theorem, volume ratio, and Euclidean structures in normed spaces is enormous, here is (in a chronological order) a sample of classical  references \cite{DilSzar, BoMi87, Mi87, Bou88, Gia-Mil, ManTom, FJ19, GLMP}.   

\vskip0.7cm

\section{Bodies with extremal or near-extremal distance to the Euclidean ball} \label{sec:extreme}

In this section we will collect and cross-reference results concerning bodies, which have extremal or near-extremal distance to the Euclidean ball. In the not-necessarily-symmetric case, the picture is reasonably clear. First, as we pointed out in the comments following John's Theorem \ref{john}, the only body $C \in \mathcal{C}^n$ for which $D_{BM}(C,B) = n$ is the simplex. Stable version of this fact was proved by Kobos in \cite{Kobos19} and subsequently improved by Brandenberg, Gonz\'alez Merino and Grundbacher in \cite{BGG}, who showed the following. 
\begin{thm} \label{thm:BGGstable} 
 If $K \in\mathcal{C}^n$  verifies $d_{BM}(K,B)\geq n-\ve$  for some $\ve \in [0,\frac 12)$, then $d_{BM}(K,T)\leq  \frac{1-2\ve^2}{1-2\ve}$.
\end{thm}
Theorem \ref{thm:BGGstable} is a consequence of the following elegant result from \cite{BGG}, which unifies the statements about symmetric and not-necessarily-symmetric bodies from John's Theorem \ref{john}. 
\begin{thm} \label{thm:BGGas}
 If $K \in\mathcal{C}^n$, then  
 \begin{equation} \label{eq:BGGas}
  d_{BM}(K,B)\leq \sqrt{n \cdot{\rm as}(K)}. 
 \end{equation}
\end{thm}
Indeed, if $d_{BM}(K,B)\geq n-\ve$, then Theorem \ref{thm:BGGas} implies that ${\rm as}(K) \geq \frac{(n-\ve)^2}n > n-2\ve$ and it remains to appeal to Theorem \ref{thm:BGG}. 

Surprisingly, the situation is more complicated in the setting of $\mathcal{M}^n$. While, as pointed out already in Section \ref{sec:john}, $d_{BM} (M,B) = \sqrt n$ when $M$ is a parallelotope or a cross-polytope, complete characterization of {\em all} extremal bodies is feasible only in low dimensions.    For $n=2$, the only such bodies are parallelograms. This can be seen by reversely engineering the derivation of Theorem \ref{john} from Theorem \ref{decomposition}, which we presented in Section \ref{sec:john}, but was known at least since 1930's. When $n=3$, parallelepipeds and cross-polytopes are no longer affinely equivalent; the fact that all maximizers are of that form appears to have been folklore since early 1990's, but hasn't seemed to be clearly documented in the literature. We refer to the 2025 paper \cite{GruKob25} by Grundbacher and Kobos for a detailed chronological record, references, and a rigorous argument based on Theorem \ref{distance-position}. 

When $n\geq 4$, it follows (for example) from Theorem \ref{decomposition} that there are infinitely many maximizers,  which can not be easily characterized.  In fact, we have the following. 
\begin{thm} \label{thm:extreme}
Let $n\geq 4$ and let $B=B_2^n$ and $Q=B_\infty^n$ be the Euclidean ball and the cube. Then there exists a set $S$ of vertices of $Q$ of cardinality strictly smaller than $2^{n-1}$ such that any $0$-symmetric convex body $M$ satisfying 
$B\subset M \subset Q$ and $S\subset M$ verifies 
\begin{equation} \label{eq:extreme}
   d_{BM}(M,B) = d_{\rm vol}(M,B) = \sqrt{n}.
\end{equation}
\end{thm}
Since the polytope $M_0 :=\conv\big(B\cup S \cup (-S)\big)$ is a proper subset of $Q$ (in many cases much smaller than $Q$), any ($0$-symmetric convex) body that is intermediate between $M_0$ and $Q$ will be a maximizer for $d_{\rm vol}(\,\cdot\,,B)$ (and hence for $d_{BM}(\,\cdot\,,B)$). 

The argument leading to \eqref{eq:extreme} is well known and versions off it appeared in the literature, see for example \cite{Tcaciuc}.  We note first that $B$ is the ellipsoid of largest volume contained in $Q$ (this follows from the uniqueness of such ellipsoid and symmetries, or from the fact that the standard basis $e_1,e_2,\ldots,e_n$ yields a collection of ``inner'' contact points verifying \eqref{eq:john} with all $\lambda_k=1$). {\em A fortiori}, $B$ is the ellipsoid of largest volume contained in $M$. Next, by the analogue of \eqref{eq:john} for L\"owner ellipsoid, there is a set $S$ of ``outer'' contact points of $Q$ and $\sqrt{n}\hskip0.4mm B$ of cardinality at most $n(n+1)/2$ such that every $0$-symmetric ellipsoid containing $S$ has volume at least as large as that of $\sqrt{n}\hskip0.4mm B$. In other words (cf. Definitions \ref{def:vr} and \ref{def:dvol}), the hypotheses of Theorem \ref{thm:extreme} force  ${\rm vr}(B,M) \cdot {\rm vr}(M,B) \geq \Big(\vol(\sqrt{n}B)/\vol(B)\Big)^{1/n} = \sqrt{n}$, which exactly means that $d_{\rm vol}(B,M) \geq \sqrt{n}$, as needed (the converse inequality is trivial). It remains to verify that the cardinality of $S$ is smaller than $2^{n-1}$. For $n\geq 5$, this follows from $n(n+1) < 2^n$; for $n=4$ we can directly see that there is an ``eligible'' set $S$ of four outer contact points given by rows of any $4 \times 4$ Hadamard matrix (see Section \ref{para-cross}). 

The instance $n=4$, or any other Hadamard number,  is most spectacular. In that case, the sufficient condition for $M \in \mathcal{M}^n$ to be extremal (i.e., $d_{BM}(M,B)=\sqrt{n}$) reduces to 
\begin{equation} \label{eq:extreme2}
    M_0 \subset M \subset Q=B_\infty^n,
\end{equation}
where $M_0$ is inscribed in $Q$ and is similar to $\sqrt{n}\hskip0.4mm B_1^n$. 

For $n=2$, there is also a stability result in the spirit of Theorem \ref{thm:BGGstable} : {\it If $K \in \mathcal{M}^2$ and $d_{BM}(K,B) \geq \sqrt{2}- \ve$, then $d_{BM}(K,P) -1 = O(\ve)$} (see \cite{GruKob25}). In view of Theorem \ref{thm:extreme}, no result of similar form seems possible when $n\geq 4$.  This leaves open the case $n=3$ (we are not aware of it being addressed in the literature). 
\begin{prob} \label{stable3}
 If $K \in \mathcal{M}^3$ and $d_{BM}(K,B) \geq \sqrt{3}- \ve$, does it follow that 
 $$
 \min\{d_{BM}(K,B_1^3),d_{BM}(K,B_\infty^3)\} -1 = O(\ve), 
 $$
or at least  $O(\ve^\alpha)$ for some (universal) $\alpha \in (0,1)$?   
\end{prob}

Stability results mentioned above work only if $d_{BM}(K,B)$ is very close to the extremal value. However, {\em some} information can be obtained under much weaker ``isomorphic'' hypotheses. We will mention here two such results, describing them in the language of normed spaces, in which they usually appear in the literature. 
\begin{thm} \label{thm:MW}
Given $\alpha \in (0,1]$ there exist constants $c_1(\alpha), c_2(\alpha) >0$, depending only on $\alpha$, such that if 
an $n$-dimensional normed space satisfies $d(X,\ell_2^n) \geq \alpha\sqrt{n}$, then $X$ contains a $k$-dimensional subspace $F$ with $k \geq c_1(\alpha) \log n$ and  verifying $d(F,\ell_1^k) \leq c_2(\alpha)$. 
\end{thm}
Theorem \ref{thm:MW} originates in the work \cite{MW} of Milman and Wolfson, who showed a weaker version with (unspecified) sequence $k_n \to \infty$ in place of $c_1(\alpha) \log n$. The version stated above can be found in the monograph \cite{[Tom 89]}; it is credited in \cite{Gia-Mil} to the ``work of Kashin, Bourgain and Tomczak-Jaegermann.''

As unsatisfying as it may seem, the order  $\Omega(\log n)$ in Theorem \ref{thm:MW} is optimal, as can be seen by considering $X =\ell_\infty^n$. This apparent obstruction suggest the following question.
\begin{prob}
Is it possible to improve the dependence of $\dim F$ on $n$ in Theorem \ref{thm:MW}  if instead of 
$d(F,\ell_1^k) \leq c_2(\alpha)$ we require that 
$\min\{d(F,\ell_1^k),d(F,\ell_\infty^k)\} \leq c_2(\alpha)$? 
\end{prob}

The second result of ``isomorphic'' type that we want to mention is due to Anisca, Tcaciuc and Tomczak-Jaegermann \cite{ATT} (see also \cite{Tcaciuc}). Since the precise statement is rather technical, we will only describe it here in general terms. In a nutshell, it says that  $d(X,\ell_2^n) = \Omega(\sqrt{n})$ if an only if, for some $k = \Omega(n)$ and some $k$-dimensional subspace $F\subset X$, the configuration from \eqref{eq:extreme2} can be approximately reproduced in $F$, where ``approximately'' means that $M_0$ and $Q$ need to be replaced by homotheties (with ratios depending only on the constants implicit in the $\Omega(\cdot)$ notation) of their ``well-conditioned'' affine images. 
Given that $d(F,\ell_2^k)=\Omega(\sqrt{n})$, which is feasible if $k := \dim F=\Omega(n)$, obviously forces  $d(X,\ell_2^n) \geq \Omega(\sqrt{n})$, the most that we may expect to derive from the latter condition is information about some subspace of $X$ of dimension $\Omega(n)$. Accordingly, the results of \cite{ATT} are the best ``characterization'' of spaces with ``nearly extremal'' distance to the Euclidean space that we can hope for. 

What is missing in this circle of ideas, are analogues of Theorem \ref{thm:MW} and the above-mentioned results from \cite{ATT} in the not-necessarily symmetric case. Here is a sample question that can be asked.
\begin{prob}
    Does there exist a sequence $k_n\to \infty$ such that if $C \in \mathcal{C}^n$ verifies $d_{BM}(C,B)=\Omega(n)$, then there is $C'$, a $k_n$-dimensional affine image of $C$, such that $d_{BM}(C', T^{k_n}) = O(1)$?  Is it possible to have $k_n= \Omega(n)$?
\end{prob}
   
\vskip 0.2cm
\section{Miscellanea and further reading} \label{sec:misc} 

The research investigating the Banach-Mazur distance and the Banach-Mazur compactum, and/or meaningfully applying these concepts is enormous. Out of necessity, this survey covers only a limited selection of topics. 
Here is a selection of papers addressing some of the other aspects of the subject in question, organized by themes. 

When considering the Banach-Mazur compactum as a geometric object, we mostly concentrated on global geometric parameters such as the diameter or radii with various centers. Other points of view include studying it as a manifold, a topological space, or investigating its complexity/metric entropy, see   \cite{SB}, \cite{SBF}, \cite{ABR03}, \cite{Osa03}, \cite{ABR04},  \cite{Szar14}, \cite{Pis15}, \cite{AriKov}, \cite{CheLyuGuoWu}, \cite{KobSwa}. 

We barely touched the very interesting topic of approximation of convex bodies by polytopes (with not too many faces or vertices), which has implications, in particular, to the question of complexity of the Banach-Mazur compactum and for the reduction of convex optimization problems to linear programming. This approximation can be quantified in the language of the BM-distance, but for some purposes approximation in terms of volume may be more relevant. A very partial bibliography of the topic includes \cite{FLM}, \cite{BF88}, \cite{CP88}, \cite{Bar97}, \cite{SW03}, \cite{Bar14}, \cite{Szar14}, \cite{AuSz17a}, \cite{BCH17}, \cite{Huang19}, \cite{Nas19}, \cite{NNR20}, \cite{PSW22}, \cite{HR26}. 

Except for a few isolated comments, we focused on the Banach-Mazur distance for finite-dimensional normed spaces or convex bodies. However, Banach-Mazur distance perfectly makes sense in the infinite-dimensional setting. Here is a sample of the literature on the topic: \cite{Lus76}, \cite{ApSe}, \cite{BeCa}, \cite{Jar},  \cite{Pia},  \cite{Ger}, \cite{Dirk},  \cite{Ron}, \cite{RonSom}, \cite{GerPia}, \cite{PiaVil}, \cite{KorPle}, \cite{MalPia}, \cite{Villada}, \cite{KorPle2}, \cite{Malec26}. In this context, related concepts are those of Kadets or Gromov-Hausdorff distances, see \cite{Kadets75}, \cite{Ost84}, \cite{KO99}, \cite{CDK22}.

We discussed BM-distancs  within various special classes of normed spaces/convex bodies primarily in Section \ref{sec:classes} and approximation by bodies from a certain class (for example) in Section \ref{sec:rectangles}. Other topics that fit into this theme are the relationship between the 
BM-distances and the smoothness/convexity properties of the related spaces \cite{[Mi]},  \cite{DMT}, \cite{[Tom 81]}, \cite{BoMi} or BM-distances for subspaces of $L_p$ :  \cite{JS}, \cite{Le84}, \cite{SchZva}, \cite{[GGMP]}.  Another notable fact involving ``$d_{BM}$ under constraints'' is the result from \cite{BoMi} (see also the proof of Theorem 46.3 in \cite{[Tom 89]}) stating that, for an $K \in \mathcal{M}^n$, 
\begin{equation}\label{K-Kpolar}
d_{BM}(K, K^\circ) = O\big(\log^6(1+n)\big)d_{BM}^{5/3}(K,B) = O\big(n^{5/6}\log^6(1+n)\big), 
\end{equation} 
which, among others, found applications in foundations of physics, see \cite{ALPSW}. It would be interesting to determine whether the bound in \eqref{K-Kpolar} is precise on the power scale. 
\begin{prob} \label{prob:K-Kpolar}
Is the bound in  \eqref{K-Kpolar} optimal (up to a polylogarithmic factor)?  
For $K\in \mathcal{C}^n$, can we find an upper bound for  $d_{BM}(K,K^\circ)$ that  goes beyond Theorem \ref{bizeul-klartag}? For the second question, the more appropriate quantity to consider may be $\min_{z\in {\rm int} K}d_{BM}\big(K,(K-z)^\circ\big)$.
\end{prob}

Finally, let us mention miscellaneous applications in geometry and convex geometry, for example \cite{KLS}, \cite{Dek}, \cite{Bor10}, \cite{BorHug}, \cite{DiaJim}, \cite{MitSai},  \cite{Usher}, \cite{ZhaHe}, \cite{RZ26}.

\baselineskip 13pt

\vskip0.25cm
\noindent
Marek Lassak

\noindent
{\it University of Science and Technology, Bydgoszcz, Poland}

\vskip0.1cm
\noindent
Email: {lassak@pbs.edu.pl}

\vskip0.3cm
\noindent
Stanis{\l}aw J. Szarek

\noindent
{\it Case Western Reserve University, Cleveland, Ohio 44106-7058, U.S.A.}  and

\noindent
{\it Institut de Math\'ematiques de Jussieu-PRG, Sorbonne Universit\'e, 4 place Jussieu,\\ 75005 Paris, France}

\vskip0.1cm
\noindent
Email: {szarek@cwru.edu}


\begin{thebibliography}{999}
\addcontentsline{toc}{section}{References}


\baselineskip 7.5pt \small


\bibitem{Ader}
O. B. Ader, An affine invariant of convex regions, {\it Duke Math. J.} {\bf 4} (1938), no. 2, 291--299.

\bibitem{ABR03}
S. M. Ageev, S. A. Bogatyi and  D. D. Repo$\check{\rm v}$s, The complement $Q_E(n)$ of the point Eucl of Euclidean space in the Banach-Mazur compactum $Q(n)$ is a $Q$ manifold, {\it Russ. Math. Surv.} {\bf 58} (2003),  607--609. 

\bibitem{ABR04}
S. M. Ageev, S. A. Bogatyi and  D. D. Repo$\check{\rm v}$s, Banach-Mazur compacta are Aleksandrov compactifications of $Q$-manifolds, {\it Mathematical Notes} {\bf 76}, no. 1 (2004), 3--9.


\bibitem{AMP26}
J. Alonso, P. Mart\'in and P. L. Papini, 
Wheeling around Chebyshev centers and Jung constant in normed planes, 
{\it Rend. Circ. Mat. Palermo} (2) {\bf 75} (2026), no. 5, Paper No. 131, 17 pp.

\bibitem {ATT}
R. Anisca, A. Tcaciuc and N. Tomczak-Jaegermann, Structure of normed spaces with extremal distance to the Euclidean space, {\it Houston J. Math.} 31 (2005), No. 1, 267--283 (electronic). 					  	
	
\bibitem {ApSe} 
J. Appel and J. M. Semenov, Estimating the Banach-Mazur distance of some pairs of symmetric spaces, {\it Math. Nachr.}, {\bf 132} (1987), 7--14. 

\bibitem{AriKov}
A. Arias and V. Kovalchuk, A Remark on geodesics in the Banach-Mazur distance, {\it Proc. Amer. Math. Soc.} {\bf 55} (2023), No. 9., 3959--3965.

\bibitem{ArtGiaMi-1} 
S. Artstein-Avidan, A. Giannopoulos and V. D. Milman, Asymptotic geometric analysis. Part I, {\it Mathematical Surveys and Monographs} {\bf 261} (2015), Providence, RI: American Mathematical Society, 451 p. 

\bibitem{ArtGiaMi-2}  
S. Artstein-Avidan, A. Giannopoulos and V. D. Milman, Asymptotic geometric analysis. Part II, {\it Mathematical Surveys and Monographs} {\bf 261} (2021), Providence, RI: American Mathematical Society, 645 p.

\bibitem {Aspl}  
E. Asplund, Comparison between plane symmetric convex bodies and parallelograms, {\it Math. Scand.} {\bf 8} (1960), 171--180. 

\bibitem {AuSz17} 
G. Aubrun and S. J. Szarek, {\it Alice and Bob Meet Banach:  The Interface of Asymptotic Geometric Analysis and Quantum Information Theory}, Mathematical Surveys and Monographs, Vol. {\bf 223}, Amer. Math. Soc. 2017.  

\bibitem {AuSz17a} 
G. Aubrun and S. J. Szarek, 
Dvoretzky's theorem and the complexity of entanglement detection,  {\it Discrete Anal.} {\bf 2017}:1, 20 pp.

\bibitem {ALPSW}
G. Aubrun, L. Lami, C. Palazuelos, S. J. Szarek and A. Winter, 
Universal Gaps for XOR Games from Estimates on Tensor
Norm Ratios, {\it Commun. Math. Phys.} {\bf 375} (2020),  679--724.

\bibitem {BY26} E. Bakaev and A. Yehudayoff, A simplex-based measure of symmetry, arXiv:2607.03815 (2026).

\bibitem {BHP01} R. C.  Baker,  G. Harman, and J. Pintz,  The difference between consecutive primes, II, {\it Proc. Lond. Math. Soc.}  {\bf 83} (3) (2001),  532--562. 

\bibitem{BMT}
V. Balestro, H. Martini and R. Teixera, Convexity from the geometric point of view, Springer, Cham (2024), 1184 p + xix.  

\bibitem {Ball89} 
K. Ball, Volumes of sections of cubes and related problems, {\it Geometric aspects of functional analysis (1987–88)}
Lecture Notes in Math. {\bf 1376} (1989), 251--260,
Springer-Verlag, Berlin.  

\bibitem {Ball91} 
K. M. Ball, Volume ratios and a reverse isoperimetric inequality, {\it J. London Math. Soc.} {\bf 44} (1991), no. 2, 351--359.

\bibitem {Ball92} 
K. Ball, Ellipsoids of maximal volume in convex bodies, {\it Geom. Dedicata} {\bf 41} (1992), No. 2, 241--250. 

\bibitem {Ball95} 
K. Ball, Mahler's conjecture and wavelets, 
{\it Discrete Comput. Geom.} {\bf 13} , no. 3-4 (1995), 271--277.

\bibitem {Ball97} 
K. Ball, {\it An elementary introduction to modern convex geometry, Flavors of geometry,} Math. Sci. Res. Inst. Publ., vol. {\bf 31}, Cambridge Univ. Press, Cambridge, 1997, pp. 1--58.

\bibitem {Ban32} 
S. Banach, {\it Th\'eorie des op\'erations lin\'eaires,} Monogr. Mat. {\bf 1}. Warszawa (1932).  [English translation: 
Theory of linear operations. Translated from the French by F. Jellett. With comments by A. Pe{\l}czy\'nski and Cz. Bessaga. North-Holland Mathematical Library, 38. North-Holland Publishing Co., Amsterdam, 1987.] 

\bibitem {BLPS99}
W. Banaszczyk, A. E. Litvak, A. Pajor and S. J. Szarek,  The flatness theorem for nonsymmetric convex bodies via the local theory of Banach spaces, {\it Math. Oper. Res.} {\bf 24} (1999), No. 3, 728--750. 

\bibitem{BF88} 
I. B\'ar\'any and Z. F\"uredi, Approximation of the sphere by polytopes having few vertices, 
{\it Proc. Amer. Math. Soc.} {\bf 102} (1988), no. 3,  651--659.

\bibitem{Bar97} 
I. B\'ar\'any, 
Approximation by random polytopes is almost optimal, 
II International Conference in ``Stochastic Geometry, Convex Bodies and Empirical Measures'' (Agrigento 1996). {\it Rend. Circ. Mat. Palermo} (2)
Suppl. No. {\bf 50} (1997), 43--50.

\bibitem {Bar02} 
A. Barvinok, {\it A course in convexity,} Graduate Studies in Mathematics, vol. {\bf 54}, American Mathematical Society, Providence, RI, 2002.

\bibitem{Bar14}  
A. Barvinok, Thrifty Approximations of Convex Bodies by Polytopes, {\it Int. Math. Res. Notices} {\bf 2014}, no. 16  (2013), 4341--4356.

\bibitem {BarVit}
Y. M. Baryshnikov and R. A. Vitale,  Regular simplices and Gaussian samples, {\it Discrete Comput. Geom.} {\bf 11} (1994), No. 2, 141--148.


\bibitem {BasRom} 
J. Bastero and M. Romance, John's decomposition of the identity in the non-convex case, {\it Positivity} {\bf 6} (2002), 1--16. 

\bibitem{BeCa}
E. Behrends, M. Cambern, An isomorphic Banach-Stone theorem, {\it Stud. Math.} {\bf 90} (1988), No. 1, 15--26.

\bibitem{BeGo}  
Y. Benyamini and Y. Gordon, Random factorization of operators between Banach spaces, {\it J. d'Anal Math.} {\bf 39} (1981), 45--74.

\bibitem{Bes}
A. S. Besicovitch, Measures of asymmetry of convex curves, {\it J. London Math. Soc.} {\bf 23} (1945), 237--240. 

\bibitem {[BiRa]} 
A. Bielecki and K. Radziszewski, Sur les parall\' el\' epip\`edes inscrits dans les corps convexes, {\it Ann. Univ. Mariae Curie-Sk\l odowska, Sect A} {\bf 7} (1954), 97--100.


\bibitem{Bizeul}
P. Bizeul, {\it Optimal $MM^*$ bounds for convex bodies}, arXiv:2607.29458 (2026), 20pp.

\bibitem{BK}
P. Bizeul and B. Klartag,  Distances between non-symmetric convex bodies: optimal bounds up to polylog,  arXiv:2510.20511 (2025), 30pp.

\bibitem {Bor05} 
K. B\"or\"oczky, Jr., The stability of the Rogers-Shephard inequality and of some related inequalities, {\it Adv. Math.} {\bf 190} (2005), No. 1, 47--76. 

\bibitem {Bor10} 
K. B\"or\"oczky, Jr., Stability of the Blaschke-Santal\'o and the affine isoperimetric inequality, {\it Adv. Math.} {\bf 225} (2010), No. 4, 1914--1928. 


\bibitem{BorHug}
K. B\"or\"oczky, Jr. and D. Hug, Stability of the reverse Blaschke-Santal\'o inequality for zonoids and applications, {\it Adv. in Appl. Math.} {\bf 44} (2010), No. 4. 309--328.

\bibitem {Born}  
W. B\"orner, An extremal problem on hexagons for the computation of special Banach-Mazur distance, {\it Beitr. Algebra Geom.} {\bf 30} (1990), 169--180. 

\bibitem{Bou84} 
J. Bourgain, On martingale transforms in finite dimensional lattices with an appendix on the $K$-convexity constant, {\it Math. Nachr.}  {\bf 119}  (1984), 41--53.

\bibitem{Bou88}  
J. Bourgain, On finite-dimensional homogeneous Banach spaces, {\it Lecture Notes in Maths.} {\bf 1317} (1988), 232--238.

\bibitem {BoMi}  
J. Bourgain and V. D. Milman, Distances between normed spaces, their subspaces and quotient spaces, {\it Eq. Operator Th.} {\bf 9} (1986), 31--46.

\bibitem {BoMi87} 
J. Bourgain and V. D. Milman, New volume ratio properties for convex symmetric bodies in $\bR^n$, {\it Invent. Math.} {\bf 88} (1987), no. 2, 319--340.

\bibitem {[BoSz88]} 
J. Bourgain and S. J. Szarek, The Banach-Mazur distance to the cube and the Dvoretzky-Rogers factorization, {\it Israel J. Math.} {\bf 62} (1988), No. 2, 169--180. 

\bibitem {BGG}
R. Brandenberg, B. Gonz\'alez Merino and F. Grundbacher, Tight stability estimates near the simplex and applications for the Banach-Mazur distance and Rogers-Shephard-type inequalities, arXiv:2607.27041 (2026), 31pp.

\bibitem {BG25}
R. Brandenberg and F. Grundbacher, Tightening inequalities on volume extremal $k$-ellipsoid using asymmetry measures,\; {\it Mathematika} {\bf 71} (2025),  no. 4, Article ID e70051 (43p.).

\bibitem {BK15} R. Brandenberg, S. K\"onig, Sharpening geometric inequalities using computable symmetry measures, {\it Mathematika}, {\bf 61} (2015), no. 3, 559--580.

\bibitem {[BL]}
 P. Brass and M. Lassak, Problems on approximation by triangles, {\it Geombinatorics} {\bf 10} (2001), No. 3, 103--115. 

\bibitem {[BMP]} 
P. Brass, W. Moser and J. Pach, {\it Research problems in discrete geometry,} Springer, 2005. 

\bibitem{BCH17} 
S. Brazitikos, G. Chasapis and L. Hioni, 
Random approximation and the vertex index of convex bodies. 
{\it Arch. Math.} (Basel) {\bf 108} (2017), no. 2, 209--221.

\bibitem {Bri16} 
J. Brinkhuis, Inner and outer approximation of convex sets using alignment, {\it Optim. Lett.} {\bf 10}, No. 7 (2016), 1403--1416.

\bibitem{BPP18}  
S. Brodiuk, N. Palko and A. Prymak, On Banach-Mazur distance between planar convex bodies, {\it Aequat. Math.} {\bf 92} (2018), 993--1000.

\bibitem{Byl}
O. P. Bylinkina, Some estimates of the Banach-Mazur distance between finite-dimensional $L\sb{p,q}$ spaces. (Russian) 
In {\it Theory of operators in function spaces},  Voronezh. Gos. Univ., 1983, 14--28.  

\bibitem{Cambern}
M. Cambern, On mappings of sequence spaces, {\it Studia Math.} {\bf 30} (1968), 73--77. 

\bibitem{CDGGJR}
A. Caragea, K. von Dichter, K. K. Gottwald, F. Grundbacher, T. Jahn and M. Runge, 
Minkowski chirality: a measure of reflectional asymmetry of convex bodies. 
{\it J. Math. Anal. Appl.} {\bf 557} No. 1 (2026) 130287, 34 pp. 

\bibitem{CP88} 
B. Carl and A. Pajor, Gelfand numbers of operators with values in a Hilbert space, {\it Invent. Math.} {\bf 94} (1988), 479--504. 

\bibitem{CEKP07}
P. G. Casazza, D. Edidin, D. Kalra and V. I. Paulsen. Projections and the Kadison-Singer problem. Operators and Matrices {\bf 1} (3):391408, 2007.

\bibitem {[CS 67]} 
G. D. Chakerian and A. K. Stein, Some intersection properties of convex bodies, {\it Proc. Amer. Math. Soc.} {\bf 18} (1967), 109--112.

\bibitem{CheLyuGuoWu}  
Y. Cheng, Y. Lyu , S. Gao and S. Wu, Construction of $\ve$-nets for the space of planar convex bodies endowed with the Banach-Mazur metric, {\it Mathematics} {\bf 13} (2025), 1358--1370. 

\bibitem{CC}
C. Chu, H. B. Cohen, Isomorphisms of spaces of continuous functions, {\it Pacific J. Math.} {\bf 155} (1992), No. 1, 71--86.

\bibitem{CD-handbook} C. J. Colbourn and J. H. Dinitz, Editors, {\it Handbook of Combinatorial Designs}, 2nd Ed. Chapman \& Hall, Boca Raton, London, New York, 2007.

\bibitem{CDK22}
M. C\'uth, M. Doucha and O. Kurka, 
Complexity of distances: reductions of distances between metric and Banach spaces, {\it Israel J. Math.} {\bf 248} (1) (2022), 383--439.

\bibitem{DalDasLauStr}
H. H. Dales, F. K. Dashiell, A. Lau and D. Strauss, Banach spaces functions as dual spaces, {\it CMB Books in Mathematics}, Cham: Springer (2016), 277p.

\bibitem{DGRZ}  
J. Dardnnes, J. Gutt, V. G. B. Vinitius and J. Zhang, Coarse distance from dynamically convex to convex, {\it Int. Math. Res. Not.} {\bf 2025}, no. 12, 17p. 

\bibitem{Dav}
W. Davis, The distance of symmetric spaces from $\ell_p^{(n)}$, {\it Lecture Notes in Maths.} {\bf 604} (1977).

\bibitem {[DE 77]} 
W. Davis and P. Enflo,  The distance of symmetric spaces from $\ell_p^n$, {\it Banach spaces of analytic functions (Proc. Pelczynski Conf., Kent State Univ., Kent, Ohio, 1976)}, pp. 25--28. Lecture Notes in Math., Vol. {\bf 604}, Springer, Berlin, 1977. 

\bibitem{DMT}  
W. J. Davis, V. D. Milman and N. Tomczak-Jaegermann, The distance between certain $n$-dimensional Banach spaces, {\it Israel J. Math.} {\bf 39} (1981), No. 1-2, 1--15. 

\bibitem {[DM 78]}
 W. Davis and B. Maurey, {\it The distance of a symmetric space from some $\ell_p$ space,}\; Proceedings of the International Conference on Operator Algebras, Ideals, and their Applications in Theoretical Physics\; (Leipzig, 1977), pp. 58--68, Teubner, Leipzig, 1978. 

\bibitem{Dek}
B. V. Dekster, Simplexes with prescribed edge lengths in Minkowski and Banach spaces, {\it Acta Math. Hung.} {\bf 86} (2000), No. 4, 343--358.

\bibitem{de'Val} 
B. A. de'Valcourt, Axially symmetric polygons inscribed in and circumscribed about convex sets, {\it Elem. Math.} {\bf 22} (1967), 121--133.

\bibitem{DiaJim}  
R. Diaz-Garcia and A. Jim\' enez-Melado, A norm inequality in James' space and stability of fixed point property, {\it J. Fixed Point Theory Appl.} {\bf 21} (2019), No. 1. Paper No. 2, 9p. 

\bibitem{DilSzar}
S. Dilworth and S. Szarek, The cotype constant and an almost Euclidean decomposition for finite-dimensional normed spaces. {\it Isr. J. Math.} {\bf 52} (1985), 82--96.


 \bibitem{Dirichlet}
 P. G. L. Dirichlet, Sur une nouvelle m\'ethode pour la d\'etermination des int\'egrales multiples. (French) {\it J. Math. Pures Appl.} {\bf 4} (1839), 164--168.

\bibitem{Dirk}
W. Dirk, On Wolfgang Lusky's paper ``The Gurarij spaces are unique", {\it Arch. Math.}, {\bf 121} (2023), No. 5--6, 615--624.

\bibitem{Dok10} D. Dokovi\'c,  Hadamard matrices of small order and Yang conjecture. {\it J. Combin. Designs} {\bf 18} (2010), 254--259.

\bibitem{Dvoretzky} A. Dvoretzky,  Some results on convex bodies and Banach spaces, {\it Proc. Internat. Sympos. Linear Spaces (Jerusalem, 1960)}, Jerusalem Academic Press, Jerusalem;
Pergamon, Oxford (1961), 123--160.

\bibitem{DR50} A. Dvoretzky and C. A. Rogers, Absolute and unconditional convergence in normed linear spaces, {\it Proc. Nat. Acad. Sci. U.S.A.} {\bf 36} (1950), 192--197.

\bibitem {FMRWY92} 
R. Fleischer, K. Mehlhorn, G. Rote, E. Welzl and C. Yap, Simultaneous inner and outer approximation of shapes, 
{\it Algorithmica} {\bf 8} (1992), 365--389. 

\bibitem{Figiel75}
T. Figiel,  A short proof of Dvoretzky's theorem, 
{\it S\'eminaire {M}aurey-{S}chwartz 1974--1975}, Exp. No. XXIII, 6 pp. (erratum, p. 3),
\'Ecole Polytech., Paris, 1975.

\bibitem{FJ19}
T. Figiel and W. Johnson, Quotients of essentially Euclidean spaces, {\it Can. Math. Bull.} {\bf 62} (2019), No. 1, 71--74.

\bibitem{FKP} 
T. Figiel, S. Kwapie\'n and A. Pe\l czy\'nski, 
Sharp estimates for the constants of local unconditional structure of Minkowski spaces.  {\it Bull. Acad. Polon. Sci. Sér. Sci. Math. Astronom. Phys.} {\bf 25}  (1977), no. 12, 1221--1226. 

\bibitem{FLM}
T. Figiel, J. Lindenstrauss, and V. D. Milman, The dimension of almost spherical sections of convex bodies,  {\em  Acta Math.} {\bf 139} (1977), no. 1-2, 53--94.

\bibitem{FNZ}
F. Fodor,  M. Nasz\'odi and T. Zarn\'ocz, 
On the volume bound in the Dvoretzky-Rogers lemma. 
{\it Pacific J. Math.} {\bf 301} (2019), no. 1, 89--99.

\bibitem{FS75}
H. Freeman and R. Shapira, Determining the minimum-area encasing rectangle for an arbitrary closed curve. {\it Comm. ACM} {\bf 18} (1975), 409--413.

\bibitem {Fri1} 
O. Friedland, An improved lover bound on the Banach-Mazur distance to the cross-polytope,  arXiv:2602.10665 (2026), 33pp.

\bibitem {Fri2} 
O. Friedland, A $5/8$ lower bound on the Banach-Mazur distance to the cross-polytope, arXiv:2608.17743 (2026), 12pp.

\bibitem {Fri3} 
O. Friedland, The sharp exponent for the Banach-Mazur distance to the cube, arXiv:2610.36917 (2026), 16pp.

\bibitem {FNY} 
O. Friedland, A. Naor and P. Youssef, 
On the Banach-Mazur radius of $\ell_\infty^n$, 
arXiv:2608.21112 (2026), 7pp.


\bibitem{GMP21} 
D Galicer, M Merzbacher and D Pinasco, Asymptotic estimates for the largest volume ratio of a convex body, {\it Rev. Mat. Iberoam.} {\bf 37} (2021), no. 6, 2347--2372. 

\bibitem{GMP26} 
D Galicer, M Merzbacher and D Pinasco, 
An Improved Volume-Ratio Bound via Isotropic Positions, 
arXiv:2608.02825 (2026), 8 p.

\bibitem{Gao}
F. Gao, 
Volumes of Generalized Balls, {\it Amer. Math. Monthly} {\bf 120}(2) (2013), 130.


\bibitem{Ger}  
A. Gergont, A note on the Banach-Mazur distance between $c_0$ and other $\ell_1$-preduals, {\it Ann. Univ. Mariae 
Curie-Sk\l odowska, Sect. A} {\bf 76} (2022), No. 1, 25--30. 

\bibitem{GerPia}  
A. Gergont and \L. Piasecki, The Banach-Mazur distance between isomorphic spaces of continuous functions is not always an integer,   {\it J. Math. Anal. Appl.} {\bf 537} (2024), No. 2, 12 p. (Article ID 128305). 

\bibitem {[Gi 95]} 
A. A. Giannopoulos, A note on the Banach-Mazur distance to the cube.  In: Lindenstrauss, J., Milman, V. (eds), Geometric Aspects of Functional Analysis, {\it Oper. Theory Adv. Appl.} {\bf 77} (1995), 67--73, Birk\"hauser, Basel.  

\bibitem {[Gi 96]} 
A. A. Giannopoulos, A proportional Dvoretzky-Rogers factorization result, {\it Proc. Amer. Math. Soc.}, {\bf 124} (1996), 233--241. 

\bibitem {[Gi-EoM]} 
A. A. Giannopoulos, {\it Banach-Mazur compactum,} Encyclopedia of Mathematics,  \url{https://encyclopediaofmath.org/} (October 2017).

\bibitem{[Gia-Har]}
A. Giannopoulos and M. Hartzoulaki, On the volume ratio of two convex bodies, {\it Bull. London Math. Soc.} {\bf 34} (2002), No. 6, 703--707. 

\bibitem {Gi-Hm} 
A. A. Giannopoulos and A. Hmadi, Banach-Mazur distances and basis constants of isotropic log-concave random spaces,  arXiv:2604.12692 (2026).

\bibitem {Gia-Mil} 
A. Giannopoulos and V.D. Milman, Euclidean structure in finite dimensional normed spaces. In {\em Handbook of the Geometry of Banach Spaces}, Vol. {\bf 1}, pp. 709--777,  North-Holland, 2001.

\bibitem {[GianPT]} 
A. Giannopoulos, I. Perissinaki and A. Tsolomitis, John's theorem for an arbitrary pair of convex bodies, {\it Geom. Dedicata} {\bf 84} (2001), No. 1-3, 63--79.

\bibitem {[Gl 81]}
E. D. Gluskin, Diameter of the Minkowski compactum is approximately equal to $n$ (in Russian), {\it Funktsional. Anal. i Prilozhen.,} {\bf 15}, No. 1 (1981), 72--73. English translation: {\it Funct. Anal. Appl.} {\bf 15}, No. 1 (1981), 57--58.  

\bibitem {Gl 83}  
E. D. Gluskin, On distances between some symmetric spaces (in Russian), {\it J. Soviet Math.} {\bf 22} (1983), 1841--1846.

\bibitem {[Gl 87]}
E. D. Gluskin, {\it Probability in the geometry of Banach spaces} (in Russian), Proceedings of the International Congress of Mathematics (Berkeley, Calif., 1986), 924--938, Amer. Math. Soc., Providence, RI, 1987. 

\bibitem {[Gl 79]}
E. D. Gluskin, An estimate of distances between finite-dimensional symmetric spaces (in Russian), Investigations on linear operators and the theory of functions, IX. Zap. Nauchn. Sem. Lenigrad. Otdel. Mat. Inst. Steklov. {\bf 92} (1979), 268--273. 

\bibitem {GluLit} E. D. Gluskin and  A. E. Litvak, On the asymmetry constant of a body with few vertices, {\it Geom. Dedicata} {\bf 90} (2002), 45--48.

\bibitem {GGM}
Y. Gordon, O. Gu\'edon and M. Meyer, An isomorphic Dvoretzky's theorem for convex bodies, {\it Stud. Math.} {\bf 127} (1998), No. 2, 191--200.  

\bibitem {[GGMP]}
Y. Gordon, O. Gu\'edon, M. Meyer and A. Pajor, Random Euclidean sections of some classical Banach spaces, {\it Math. Scand.} {\bf 91} (2002), No. 2, 247--268.  

\bibitem{GLMP}
D. Golicer, A. Litvak, M. Merzbacher and D. Pinasco, On the volume ratio of projection of convex bodies, {\it J. Funct. Anal.} {\bf 286} (2024), No. 3., Paper No. 110242, 27pp.

\bibitem{GomRol}
C. Gomez-Navarro and E. Rold\'an-Pensado, Transversals to colorful intersecting convex sets, {\it Discrete Comput. Geom.} {\bf 73} (2025), No. 4, 1079--1096. 

\bibitem {[GLMP]} 
Y. Gordon, A. R. Litvak, A. Meyer and A. Pajor, John's decomposition in the general case and its applications, {\it J. Differential Geom.} {\bf 68} (2004), No. 1, 99--119.

\bibitem {Grub}   
P. M. Gruber,  Aspects of approximation of convex bodies, {\it Handbook of Convex Geometry}, Vol. A, 319--345, North-Holland, Amsterdam, 1993.

\bibitem{GruKob25}  
F. Grundbacher, T. Kobos, On certain extremal Banach-Mazur distances and Ader's characterization of distance ellipsoids, 
{\it Mathematica}, {\bf 72} (2026), No.1, Paper No. e70062, 36 pp.

\bibitem{GruKob26}  
F. Grundbacher, T. Kobos, John-type decompositions for affinely-optimal positions of convex bodies,  arXiv:2602.14847 (2026).

\bibitem {[Grun60]} B. Gr\" unbaum, Projection constants. {\it Trans. Am. Math. Soc.} {\bf 7} (1960), 451--465. 

\bibitem {[Grun63]} 
B. Gr\"unbaum,  Measures of symmetry of convex sets. In {\it Convexity, Proc. Sympos. Pure Math.}, vol. {\bf 7} (1963), Amer. Math. Soc., Providence, R.I.,  233--270.

\bibitem {[Gu 05]} 
Q. Guo, Stability of the Minkowski measure of asymmetry for convex bodies, {\it Discrete Comput. Geom.} {\bf 34} (2005), No. 2, 351--362.

\bibitem {[Gu 12]} 
Q. Guo,  A simple approach to estimating the Banach-Mazur distance between simplex and convex body, \; {\it Chinese Quartely Journal of Mathematics,}\; {\bf 27} (2012), No. 4,  504--508. 

\bibitem {[GK 99]} 
Q. Guo and S. Kaijser, On the distance between convex bodies, {\it Northeast. Math. J.} {\bf 15} (1999), No. 3, 323--331.

\bibitem {[GK 02]} 
Q. Guo and S. Kaijser, On asymmetry of convex bodies, {\it Discrete Comput. Geom.} {\bf 27} (2002), 239--247.

\bibitem {[GK 03]} 
Q. Guo and S. Kaijser, Approximation of convex bodies by convex bodies, {\it Northeast. Math. J.} {\bf 19} (2003), No. 4, 323--332. 

\bibitem {GKM65} 
V. I. Gurari\u\i, M. I. Kadets and V. I. Matsaev, On Banach-Mazur distance between certain Minkowski spaces, {\it Bull. Acad. Polon. Sci. S\'er. Sci. Math. Astronom. Phys.} {\bf 13} (1965), 719--722.

\bibitem {GKM66} 
V. I. Gurari\u\i, M. I. Kadets and V. I. Matsaev, Distances between finite dimensional analogs of the $L_p$ spaces. (Russian) {\it Mat. Sb.} {\bf 70} (1966), 481--489.  

\bibitem {Haa82} U. Haagerup, The best constants in the Khintchine inequality, {\it Stud. Math.}  {\bf 70} (1982), 231--283.

\bibitem {Hadamard}  J. Hadamard,  R\'esolution d'une question relative aux d\'eterminants, {\it Bull. Sci. Math.} {\bf 17} (1893), 240--246.

\bibitem {[Had]}
H. Hadwiger, Volumsch\"{a}tzung f\"{u}r die einem Eik\"{o}rper \"{u}berdeckenden and unterdeckenden Parallelotope, {\it Elem. Math.} {\bf 10} (1955), 122--124.  

\bibitem {Henk}
 M. Henk, L\"owner-John ellipsoids,
{\it Doc. Math.} (Bielefeld) Extra Vol. Optimization Stories (2012), 95--106. 

\bibitem {HiRi05}
A. Hinrichs and  C. Richter, 
The Knaster problem: more counterexamples, 
{\it Israel J. Math.} {\bf 145} (2005), 311--324. 

\bibitem {Hmadi26} 
A. Hmadi, Gaussian polytopes with large Banach–Mazur distance to the cross-polytope, arXiv:2609.27796 (2026), 13pp.

\bibitem {HH23} 
J. Hu and K. Huang, Identifiable Bounded Component Analysis Via Minimum Volume Enclosing Parallelotope, {\it ICASSP 2023 - 2023 IEEE International Conference on Acoustics, Speech and Signal Processing (ICASSP)}, Rhodes Island, Greece  (2023),  1--5. 

\bibitem{Huang19} 
H. Huang, 
John’s position is not good for approximation, 
{\it Israel J. Math.}  {\bf 234} (2019), 37--53. 

\bibitem{HR26} 
H. Huang and  M. Rudelson,
Hardness of approximation of centered convex bodies by polytopes, 
arXiv:2602.23034 (2016),  20pp. 

\bibitem{Hua18} X. Huang, 
Notes on the Minkowski measure, the Minkowski symmetral, and the Banach-Mazur distance, {\it J. Korean Math. Soc.}  {\bf 55} (3) (2018), 695--704.

\bibitem{Jar}
K. Jarosz, Small isomorphisms on $C(X,E)$ spaces, {\it Pacific J. Math.} {\bf 138} (1989), No. 2, 295--315.

\bibitem{JN11} 
C. H. Jim\'enez and M. Nasz\'odi, On the extremal distance between two convex bodies, {\it Israel J. Math.} {\bf 183} (2011), 103--115. 

\bibitem {[John]}
F. John,  Extremum problems with inequalities as subsidiary conditions, {\it Courant Anniversary Volume}, 1948, pp. 187--204.

\bibitem {[JL]}
W. B. Johnson and J. Lindenstrauss, {Basic concepts in the geometry of Banach spaces}, Handbook of the geometry of Banach spaces, Vol.  {\bf 1}, 1--84, North-Holland, Amsterdam, 2001. 

\bibitem {JS}
W. B. Johnson and G. Schechtman,  On subspaces of $L\sb{1}$ with maximal distances to Euclidean space. Proceedings of research workshop on Banach space theory (Iowa City, Iowa, 1981), 83--96, Univ. Iowa, Iowa City, IA, 1982.  

\bibitem {Kadets75}
M.I. Kadets, Note on the gap between subspaces, 
{\it Funkts. Anal. Prilozhen.} {\bf 9} (1975), no. 2, 73--74 (Russian); English transl. in {\it Funct. Anal. Appl.} {\bf 9} (1975), 156--157.

\bibitem {KO99}
N. J. Kalton and M. I. Ostrovskii, Distances between Banach-spaces, {\it Forum Math.} {\bf 11} (1999), 17--49. 

\bibitem{KLS}  
R. Kannan, L. L\'ovasz and H. E. Scarf, The shapes of polyhedra, Math. Oper. Res. {\bf 15} (1990), 364--380. 

\bibitem{Kaplansky}
I. Kaplansky,  {\it Linear algebra and geometry. A second course.} 
Allyn and Bacon, Inc., Boston, MA, 1969. 

\bibitem{Kas}
B. S. Kashin, Sections of some finite dimensional sets and classes of smooth functions (in Russian), {\it Izv. Akad. Nauk. SSSR Ser. Mat.} {\bf 41} (1977), 334--351. 

\bibitem{KaSz03}
B. S. Kashin and S. J. Szarek, The Knaster problem and the geometry of high-dimensional cubes,
{\it C. R. Math. Acad. Sci. Paris} {\bf 336} (2003), no. 11, 931--936.

\bibitem{Kha79}
L. G. Khachiyan, 
A polynomial algorithm in linear programming. (Russian)
{\it Dokl. Akad. Nauk SSSR} {\bf 244} (1979), no. 5, 1093--1096. 
English translation in: {\it Soviet Math. Dokl.} {\bf 20} (1979), no. 1, 191--194.

\bibitem{Kha80}
L. G. Khachiyan, 
Polynomial algorithms in linear programming.(Russian)
{\it Zh. Vychisl. Mat. i Mat. Fiz.} {\bf 20} (1980), no. 1, 51--68, 260.

\bibitem{[Kh01a]}  
A. I. Khrabrov, Estimates of distances between the sums of spaces $l\sp p\sb n$, {\it Vestnik St. Petersburg Univ. Math.} {\bf 33} (2001), No. 3, 45--50.

\bibitem{[Kh01b]}   
A. I. Khrabrov, Extremal volume ratios for the sums of normed spaces, In: Function theory and partial differential equations, {\it J. Math. Sci. (New York)} {\bf 105} (2001), No. 5, 2485--2494.  

\bibitem{[Kh01c]}    
A. I. Khrabrov, Distances between spaces with unconditional bases, {\it Function theory and partial differential equations, J. Math. Sci. (New York)} {\bf 107} (2001), No. 3, 3952--3962. 
 
\bibitem{[Kh02]}   	
A. I. Khrabrov, Generalized volume ratios and the Banach-Mazur distance (in Russian), {\it Mat. Zametki} {\bf 70} (2001), No. 6, 918--926. Translation in {\it Math. Notes} {\bf 70} (2001), No. 5--6, 838--846.

\bibitem{Kh04}  
A. I. Khrabrov, 
Comparison of some distances between sums of the spaces $\ell_p^n$, English translation in
{\it Funct. Anal. Appl.} {\bf 38} (2004), no. 1, 75--77.

\bibitem{Kh05}  
A. I. Khrabrov, Estimates of distances between the sums of spaces $l\sp p\sb n$, {\it J. Math. Sci. New York} {\bf 129} (2005), No. 4, 4040--4048. 
		
\bibitem{[KN]} 
S-J. Kim and K. Nakparsit, Coloring the complements of intersection graphs of geometry figures, {\it Discr. Math.} {\bf 308} (2008), 4589--4594.

\bibitem{Knaster47} 
B. Knaster, Problem 4, {\it Colloq. Math.} {\bf 30} (1947), 30–31.


\bibitem{Kobos19}
T. Kobos, Stability result for the extremal Gr\"unbaum distance between convex bodies, {\it J.
Convex Anal.} {\bf 26} (4) (2019), 1277--1296.

\bibitem{Kobos20}
T. Kobos, Extremal Banach-Mazur distance between a symmetric convex body and an arbitrary convex body in the plane, {\it  Mathematika} {\bf 66} (2020), 161--177.

\bibitem{KobSwa}  
T. Kobos and K. Swanepoel, Equilateral dimension of the planar Banach-Mazur compactum,  arXiv:2503.21478 (2025). 

\bibitem{[KV 24]} 
T. Kobos and M. Varivoda, On the Banach-Mazur distance in small dimensions, {\it Discrete Comput. Geom.} {\bf 74} (2024), No. 2, 399--427. 

\bibitem{[KTS99]} H. K\"onig, C. Sch\"utt, and N. Tomczak-Jaegermann, Projection constants of symmetric spaces and variants of Khintchine's inequality, {\it J. Reine Angew. Math.}  {\bf 511} (1999),  1--42. 

\bibitem{Kominers}
S. D. Kominers, private communication, July 2026. 

\bibitem{KorPle}
M. Korpalski and G. Plebanek, Bounds for Banach-Mazur distance between some $C(K)$-spaces,  arXiv:2511.03435 (2025). 

\bibitem{KorPle2}
M. Korpalski and G. Plebanek, How many miles from $L_\infty$ to $l_\infty$,  {\it Proc. Amer. Mat. Soc.} {\bf 154}, No. 8 (2026), 3403--3413.

\bibitem {[L 89]}  
M. Lassak, Approximation of plane convex bodies by centrally symmetric bodies, {\it J. London Math. Soc.} {\bf 40} (1989), 369--377.                   
\bibitem {[L 91]}  
M. Lassak, Approximation of convex bodies by parallelotopes, {\it Bull. Pol. Ac. Math.} {\bf 39} (1991), 219--223. 

\bibitem {[L 92]} 
M. Lassak, Approximation of convex bodies by triangles, {\it Proc. Amer. Math. Soc.} {\bf 115} (1992), 207--210.

\bibitem{[L 92a]} 
M. Lassak, On the Banach-Mazur distance between convex bodies, {\it J. Geom.} {\bf 44} (1992), 11--12. 

\bibitem {[L 93]} 
M. Lassak, Approximation of convex bodies by rectangles, {\it Geom. Dedicata} {\bf 47} (1993), 111--117. 

\bibitem {[L 98]}
M. Lassak, Approximation of convex bodies by centrally symmetric bodies, {\it Geom Dedicata}, {\bf 72} (1998), 63--68.

\bibitem {[L 99]} 
M. Lassak, Parallelotopes of maximum volume in a simplex, {\it Discrete Comput. Geom.}, {\bf 21} (1999), 449--462.

\bibitem {[L 99a]} 
M. Lassak, On large parallelotopes in a simplex and on generalized Banach-Mazur distance, {\it Extracta Math.} {\bf 14} (1999), p. 85 (Report on the meeting ``Functional Analysis and its Applications", Gargano, Italy, 1998).

\bibitem {[L 01]} 
M. Lassak, Relationships between widths of a convex body and of an inscribed parallelotope, {sl Bull. Austral. Math. Soc.} {\bf 63}  (2001), no. 1, 133--140.

\bibitem {[L 02]}
M. Lassak, Approximation of convex bodies by axially symmetric bodies, {\it Proc. Amer. Math. Soc.} {\bf 130} (2002), 3075--3084.

\bibitem {[L 06]} 
M. Lassak, Packing an $n$-dimensional convex body by $n+1$ homothetical copies, {\it Rev. Romaine Mat. Pure Appl.} {\bf 51} (2006), 43--47.

\bibitem {[L 07]} 
M. Lassak, Banach-Mazur distance of planar bodies, {\it Aequat. Math.} {\bf 74} (2007), 282--286.

\bibitem{[L 08]}
M. Lassak, Banach-Mazur distance of central sections of a centrally symmetric convex body, {\it Beitr. Algebra Geom.} {\bf 49} (2008) No. 1, 243--246. 

\bibitem {[L 10]} 
M. Lassak, Simplices of maximum volume contained in the unit ball of a normed space, {\it Publ. Math. Debrecen} {\bf 77} (2010), 79--86.

\bibitem {[L 11]} 
M. Lassak, Approximation of convex bodies by inscribed simplices of maximum volume, {\it Beitr. Algebra Geom.} {\bf 52} (2011), No. 2, 389--394.

\bibitem {[L 14]} 
M. Lassak, Banach-Mazur distance between convex quadrangles, {\it Demonstratio Math.} {\bf 47} (2014), No. 4, 889--993.

\bibitem{[L 21]}
M. Lassak, Banach-Mazur distance from the parallelogram to the affine-regular hexagon an other affine-regular even-gons, 
{\it Results Math.} {\bf 76} (2021), art. 62, 12pp.

\bibitem{[L 24a]}
M. Lassak, The centroid Banach-Mazur distance between the parallelogram and the triangle, {\it J. Convex Anal.} {\bf 31} (2024), No. 1, 51--58. 

\bibitem{[L 24b]}
M. Lassak,  Estimation of the centroid Banach-Mazur distance between planar convex bodies, {\it Ukr. Math. J.} {\bf 76} (2024), No. 5, 872--879.

\bibitem{[L 24c]}
M. Lassak, Position of the centroid of a planar convex body, {\it Aequationes Math.} {\bf 98} (2024), No. 3, 687--695.

\bibitem {Lasse15}     
J. Lasserre, A generalization of L\"owner-John's ellipsoid theorem, {\it Math. Program., Ser. A} {\bf 152} (2015), No. 1-2, 559--591.

\bibitem {LMOT}   
R. Lata{\l}a, P. Mankiewicz, K. Oleszkiewicz and N. Tomczak-Jaegermann, Banach-Mazur distances and projections on random subgaussian polytopes, {\it Discrete Comput. Geom.} {\bf 38} (2007), 29--50.

\bibitem {Lei} 
K. Leichtweiss, \"Uber die affine Exzentrizit\"at konvexer K\"orper. (German) {\it Arch. Math.} {\bf 10} (1959), 187--199.

\bibitem{LenZhaMa}
G. Leng, Y. Zhang and B. Ma, Largest parallelotopes contained in simplices, {\it Discrete Math.} {\bf 211} (2000). No. 1--3, 111--123.  

\bibitem{Le79}   
D. R. Lewis, Ellipsoids defined by Banach ideal norms,\; {\it Mathematika}\; {\bf 26} (1979), 18--29.

\bibitem{Le84}   
D. R. Lewis, The distance from $\ell_r^n$ to a subspace of $L_p$, {\it Match. Nachr.} {\bf 119} (1984), 189--193. 

\bibitem {LiMi}   
J. Lindenstrauss and V. D. Milman,  The local theory of normed spaces and its applications to convexity. In  {\it Handbook of convex geometry}, 1149--1220, North-Holland, Amsterdam, 1993. 

\bibitem{[LiSz86]}
J. Lindenstrauss and A. Szankowski,  On the Banach-Mazur distance between spaces having an unconditional basis,\; {\it Aspects of positivity in functional analysis (T\"ubingen, 1985)}, 119--136,\; North--Holland Math. Stud.\; {\bf 122},\; Amsterdam, 1986.

\bibitem{[LiSz87]} J. Lindenstrauss and A. Szankowski, The weak distance between Banach spaces with a symmetric basis, {\it J. Reine Angew. Math.}  {\bf 373} (1987), 108--147.

\bibitem{LiSz87a} J. Lindenstrauss and A. Szankowski,  The relation between the distance and the weak distance for spaces with symmetric basis, Proc of GAFA Seminar, Tel·Aviv 1985/86, {\it Lecture Notes in Math.} {\bf 1267} (1987), 21--38, Springer Verlag, Berlin.

\bibitem{Lus76} 
W. Lusky, The Gurarij spaces are unique,
{\it Arch. Math. (Basel)} {\bf 27} (1976), no. 6, 627--635.

\bibitem{LutYanZha05}
E. Lutwak, D. Yang and G. Zhang, John ellipsoids, {\it Proc. London. Math. Soc.} III. Ser. {\bf 90} (2005), No. 2, 497--520.

\bibitem{Mac}
A. M. Macbeath, A compactness theorem for affine-eqivalence classes of convex regions, {\it Canad. J. Math.} {\bf 3} (1951), 54--61. 

\bibitem {[Mak96]} 
V. V. Makeev, Applications of topology to some problems in combinatorial geometry, {\it Amer. Math. Soc. Transl.} {\bf 174} (1996), 223--228.

\bibitem {[Mak07]} 
V. V. Makeev, On the volume of parallelepipeds circumscribed about a three-dimensional convex body. (Russian) {\it J. Math. Sci. (N.Y.)} {\bf 140} (2007), no. 4, 542--547. Translation of  : {\it Zap. Nauchn. Sem. S.-Peterburg. Otdel. Mat. Inst. Steklov. (POMI)} {\bf 329} (2005), 79--87, 196. 

\bibitem {MalPia} 
M. Malec and \L. Piasecki, The Banach-Mazur distance between $C([1, \omega^n])$ and $C([1, \omega])$, {\it Studia Math.} {\bf 285} (2025), No. 1, 91--104.

\bibitem {Malec26} 
M. Malec, The Banach–Mazur distance between $L_1$-preduals, 
Ph. D. thesis, Uniwersytet Marii Curie-Sklodowskiej, Lublin, 2026.

\bibitem {Mank12}  P. Mankiewicz, 
Weak distances between random subproportional quotients of $\ell_1^m$,  
{\it Bull. Pol. Acad. Sci. Math.} {\bf 60} (2012), no. 3, 285--294.

\bibitem{ManSza} P. Mankiewicz and S. Szarek, Random Banach spaces. The limitations of the method, {\it Mathematica} {\bf 41} (1994), 239--250,  Corrigenda {\bf 42} (1995), 220--221. 

\bibitem{ManTom} 
P. Mankiewicz and N. Tomczak-Jaegermann, 
Quotients of finite-dimensional Banach spaces; random phenomena. In {\it Handbook of the geometry of Banach spaces}, Vol. {\bf 2}, 1201--1246, 
North-Holland, Amsterdam, 2003. 

\bibitem{Mar}
E. Markesinis, On the mean width ratio of convex bodies, {\it Bull. Lond. Math. Soc.} {\bf 55} (2023), No. 3, 1254--1263.

\bibitem{Maurey}
B. Maurey, Type, Cotype and $K$-Convexity, 
In {\em Handbook of the Geometry of Banach Spaces}, Vol. {\bf 2}, pp. 1299--1332,  North-Holland, 2003.

\bibitem{MT}
T. Medina and G. Tresch, Transportation cost spaces and stochastic trees,  arXiv:2501.08656 (2025). 

\bibitem{Melistas}
T. Melistas, 
The contact Banach-Mazur distance and large scale geometry of overtwisted contact forms,  arXiv:2009.07369 (2020).

\bibitem{[Mi71]} 
V. D. Milman, A new proof of A. Dvoretzky’s theorem on cross-sections of convex bodies. (Russian) {\it Funkcional. Anal. i Prilo\v zen.} {\bf 5} (1971), no. 4, 28--37.

\bibitem{[Mi]}  
V. D. Milman, The Banach-Mazur distance, and the geometric properties of a $B$-space (in Russian) {\it Teor. Funkci{\u\i}  Funkcional. Anal. i Prilo\v zen.} {\bf 19} (1974), 23--32, 158. 

\bibitem{Mi85}  
V. D. Milman, Almost Euclidean quotient spaces of subspaces of finite dimenional normed spaces, {\it Proc. Amer. Math. Soc.} {\bf 94} (1985), 445--449. 

\bibitem{Mi87}  
V. D. Milman, 
Some remarks on Urysohn's inequality and volume ratio of cotype 2-spaces. In: Lindenstrauss, J., Milman, V.D. (eds) {\it Geometrical Aspects of Functional Analysis.}  Lecture Notes in Math. {\bf 1267},  Springer, Berlin, Heidelberg. 

\bibitem{Mi88} 
V. D. Milman, A few observations on the connection between local theory and some other fields,  {\it Lecture Notes in Math.} {\bf 1317} (1988), 283--289.

\bibitem{MS86} 
V. D. Milman and G. Schechtman, 
Asymptotic theory of finite-dimensional normed spaces.
With an appendix by M. Gromov.
{\it Lecture Notes in Math.} {\bf 1200} (1986), 
Springer-Verlag, Berlin. 

\bibitem{MS95} 
V. D. Milman and G. Schechtman, 
An ``isomorphic'' version of Dvoretzky's theorem. 
{\it C. R. Acad. Sci. Paris Sér. I Math.} {\bf 321} (1995), no. 5, 541--544. 

\bibitem{MS99} 
V. D. Milman and G. Schechtman, 
An ``isomorphic'' version of Dvoretzky's theorem. II. In: Convex geometric analysis (Berkeley, CA, 1996), 
{\it Math. Sci. Res. Inst. Publ.} {\bf 34} (1999) 159--164, 
Cambridge University Press, Cambridge. 

\bibitem{MW}
V. D. Milman and H. Wolfson, Minkowski spaces with extremal distance from the Euclidean space, {\it Israel J. Math.} {\bf 29} (1978), No. 2--3, 113--131.	

\bibitem{Mink}
H. Minkowski, Allegemeine Lehrs\"atze \"uber konvexe Polyeder, {\it Nachr. Ges. Wiss. Gottingen} (1897), 198--219.

\bibitem {[MiSa]}
A. Mishra and D. Sain, {On the Banach-Mazur distance from $\ell_1^n$}, {\it Coll. Math.} {\bf 177} (2024), 161--173.

\bibitem{MitSai}
K. Mitani and K. Saito, On von Neumann-Jordan constant of generalized Bana\'s-Fr\c aczek spaces. II., {\it Linear Nonlinear Anal.}, {\bf 8} (2022), No. 2, 217--223. 

\bibitem{Nachbin} L. Nachbin, A theorem of the Hahn-Banach type for linear transformations, {\it Trans. Amer. Math. Soc.} {\bf 68} (1950) 28--46.  

\bibitem{Nas19} 
M. Nasz\'odi,
Approximating a convex body by a polytope using the epsilon-net theorem.
{\it Discrete Comput. Geom.} {\bf 61} (2019), no. 3, 686--693.

\bibitem{NNR20} 
M.  Nasz\'odi, F.  Nazarov and D. Ryabogin, 
Fine approximation of convex bodies by polytopes, 
{\it Amer. J. Math.} {\bf 142} (2020), no. 3, 809--820. 

\bibitem {NY83} 
A. S. Nemirovsky and D. B. Yudin, 
Problem complexity and method efficiency in optimization, 
John Wiley \& Sons, Inc., New York, 1983.

\bibitem {Nev11}
V. M.  Nevski\u i, 
On the Lassak Conjecture for a Convex Body. (Russian)
{\it Model. Anal. Inf. Sist.} {\bf 18} (2011), no. 3, 5--11.

\bibitem{[Nov]}
P. Novotny, Approximation of convex bodies by simplices, {\it Geom. Dedicata} {\bf 50} (1994), 35--55. 

\bibitem{ORurke85} 
J. O'Rourke, Finding minimal enclosing boxes, {\it Internat. J. Comput. Inform. Sci.} {\bf 14} (1985), no. 3, 183--199.

\bibitem {Osa03} 
D. Osajda, The topology of the punctured Banach-Mazur compacta, 
Ph. D. Thesis (2003), Mathematical Institute,  Polish Academy of Sciences, \url{http://www.math.uni.wroc.pl/~dosaj/trav/tpbmc-web.pdf} 

\bibitem {Ost84}
M. I. Ostrovskii, On properties of the opening and related closeness characterizations of Banach spaces, {\it Amer. Math. Soc. Transl.} {\bf 136} (2) (1987), 109--119. Translation of  : {\it Teor. Funkts., Funkts. Anal. Prilozh.} {\bf 42} (1984), 97--107 (in Russian). 

\bibitem{Paley} R. E. A. C.  Paley,  On orthogonal matrices,  {\it  J. Math. and Phys.}  {\bf 12} (1-4) (1933),  311--320.  

\bibitem {[Pal 92]} 
O. Palmon, The only convex body with extremal distance from the ball is the simplex, {\it Israel J. Math.} {\bf 80} (1992), No. 3, 337--349.

\bibitem {Pelc84} 
A. Pe\l czy\'nski,\; {\it Structural theory of Banach spaces and its interplay with anaysis and probability,}\;\; Proceedings of the ICM\; 1983,\;\; PWN -- North Holland, 1984,\; pp. 237--269.

\bibitem {Pelc90} 
A. Pe\l czy\'nski, Remarks on John's theorem on the ellipsoid of maximal volume inscribed into a convex symmetric body in $R^n$, {\it Note Math.} {\bf 10} (1990), Suppl. 2,\; 395--410.

\bibitem {PeSz}   
A. Pe\l czy\'nski and S. Szarek, On parallelepipeds of minimal volume containing a convex symmetric body in $R^n$, {\it Math. Proc. Camb. Phil. Soc.} {\bf 109} (1991), 125--148.

\bibitem{Pia}  
\L. Piasecki, On Banach space properties that are not invariant under the  Banach-Mazur distance $1$, {\it J. Math. Anal. Appl.} {\bf 467} (2018), No. 2, 1129--1147. 

\bibitem{PiaVil}  
\L. Piasecki and J. Villada, The Banach-Mazur distance between $C(\Delta)$ and $C_0(\Delta)$ equals $2$, {\it Topol. Methods Nonlinear Anal.} {\bf 63} (2024), No. 1, 227--232.

\bibitem{Pis86}  
G. Pisier, Factorization of Linear Operators and Geometry of Banach Spaces, {\it CBMS Regional Conference Series in Mathematics} {\bf 60}, Amer. Math. Soc., Providence, RI, 1986.

\bibitem{Pis15}  
G. Pisier, On the metric entropy of the Banach-Mazur compactum, {\it Mathematika} {\bf 61} (2015), 179--198. 

\bibitem{PRSZ} 
L. Polterovich, D. Rosen, K. Samvelyan and J. Zhang,   
Topological Persistence in Geometry and Analysis, 
{\it Univ. Lecture Ser.}, {\bf 74}
Amer. Math. Soc., Providence, RI, 2020. See also \url{https://sites.google.com/site/polterov/topological-persistence-in-geometry-and-analysis-revision} for errata and updates. 

\bibitem {[PoSe]} 
G. P\'olya and G. Szeg\' o, Isoperimetric inequalities in mathematical physics, {\it Ann. Math. Stud.} {\bf 27}, 1951. Princeton:  Princeton University Press, XVI, 279p.

\bibitem{PSW22}
J. Prochno, C. Sch\"utt and E. M. Werner, 
Best and random approximation of a convex body by a polytope,
{\it J. Complexity} {\bf 71} (2022) 101652, 19pp.

\bibitem {[Radz1]}
C. Radziszewski, 
Sur un probl\`eme extr\'emal relatif aux figures inscrites dans les figures convexes. (French)
{\it C. R. Acad. Sci. Paris} {\bf 235} (1952), 771--773.

\bibitem {[Radz2]}
K. Radziszewski, Sur un probl\' eme extr\' emal r\'elatif aux figures inscrites et circonsrites aux figures convexes,\; {\it Ann. Univ. Mariae Curie--Sk\l odowska, Section A}, {\bf 6} (1953), 5--18.

\bibitem {RoSh1}
C. A. Rogers and G. C. Shephard, The difference body of a convex body, {\it Arch. Math.} {\bf 8} (1957), 220--233.

\bibitem {RoSh2}
C. A. Rogers and G. C. Shephard, Convex bodies associated with a given convex body, {\it J. London Math. Soc.} {\bf 33} (1958), 270--281.

\bibitem {Ron}  
J. Rond\'os, On the Banach-Mazur distance between continuous function saves with scattered boundaries, {\it Czech. Math. J.} {\bf 73} (2023), 367--393. 

\bibitem {RonSom}  
J. Rond\'os and J. Somaglia, Isomorphic of $\mathcal C(K,E)$ spaces and height of $K$, {\it Mediterr. J. Math.} {\bf 55} (2023), No. 6, Paper No. 194, 13pp.

\bibitem {[Rud 95]} 
M. Rudelson,  Estimates of the weak distance between finite-dimensional Banach spaces, {\it Israel J. Math.} {\bf 89} (1995), 189–-204.

\bibitem {[Rud 00]} 
M. Rudelson, Distances between non-symmetric convex bodies and the $MM*$-estimate, {\it Positivity} {\bf 4} (2000), 161--178.

\bibitem {[Rud 04]} 
M. Rudelson, Extremal distances between sections of convex bodies, {\it Geom. Funct. Anal.} {\bf 14} (2004), No. 5, 1063--1088.	

\bibitem{RZ26}
D. Ryabogin and A. Zvavitch, 
Zonoids whose polars are zonoids: the Banach--Mazur distance need not tend to one, 
arXiv:2609.10852 (2026), 30 pp. 

\bibitem{Sch13}
G. Schechtman, Euclidean Sections of Convex Bodies. In: Ludwig, M., Milman, V., Pestov, V., Tomczak-Jaegermann, N. (eds) Asymptotic Geometric Analysis. {\it Fields Institute Communications} {\bf 68}. Springer, New York, NY. 

\bibitem{SchZva}
G. Schechtman and A. Zvavitch, Embedding subspaces of $L_p$ into $\ell_p^N, 0<p<1$, {\it Math. Nachr.} {\bf 227} (2001), 133--142. 

\bibitem {Sch}  
C. Sch\"utt, On the Banach-Mazur distance of finite dimensional symmetric Banach spaces and hypergeometric distribution, {\it Stud. Math.} {\bf 72} (1982), 109--129. 

\bibitem{SW03} 
C. Sch\"utt and E. M. Werner, 
Polytopes with Vertices Chosen Randomly from the Boundary of a Convex Body, 
{\it Lecture Notes in Math.}  {\bf 1807},  Springer, Berlin, Heidelberg 2003, 241--422. 

\bibitem{SB}  
A. Sergei and S. A. Bogatyi, The Banach-Mazur compactum is not homeomorphic to the Hilbert cube, {\it Russ. Math. Surv.} {\bf 53} (1998), 205--207. 

\bibitem{SBF}  
A. Sergei, S. A. Bogatyi and P. Fabel, The Banach-Mazur compactum $Q(n)$ is an absolute retract, {\it Mosc. Univ. Math. Bull.} {\bf 53} (1998), 10--12. 

\bibitem{Sch09} R. Schneider, Stability for some extremal properties of the simplex. {\it J. Geom.} {\bf 96} (2009), 135--148. 

\bibitem {[SFRW]}
O. Schwarzkopf, U. Fuchs, G. Rote and E. Welzl, Approximation of convex figures by pairs of rectangles, {\it Comput. Geom.} {\bf 10} (1998), No. 2, 77--87. 

\bibitem {[SG]}
Y. Shao, Q. Guo, The equivalence and estimate of several affine invariant distances between convex bodies, {\it J. of Mathematics (PRC)} {\bf 35} (2015), No. 2, 287--292.


\bibitem {SS12}
D. A. Spielman and N. Srivastava, 
An elementary proof of the restricted invertibility theorem. 
{\it Israel J. Math.} {\bf 190} (2012), 83--91.

\bibitem {Sri10}
N. Srivastava, Spectral sparsification and restricted invertibility. Ph.D. Thesis, Yale University, March 2010.

\bibitem{SZ}
V. Stojisavljevi'c and J. Zhang, 
Persistence modules, symplectic Banach-Mazur distance and Riemannian metrics.(English summary)
{\it Internat. J. Math.} {\bf 32} (2021), no. 7, Paper No. 2150040, 76 pp.

\bibitem {[St 81]}
W. Stromquist, The maximum distance between two-dimensional Banach spaces, {\it Math. Scand.} {\bf 48} (1981), 205--225.	
	
\bibitem {Suss}  
W. S\" uss, \" Uber Parallelograme und Rechtecke, die sich ebenen Eibereichen einbeschrieben lassen, {\it Rend. Mat. e Appl.} {\bf 14} (1955), 338--341.

\bibitem{SwaVil}  
K. J. Swanepoel and R. Villa, Maximal equilateral sets, {\it Discrete Comput. Geom.} {\bf 50} (2013), No. 2, 354--373. 

\bibitem{Sylvester}   J. J. Sylvester, Thoughts on inverse orthogonal matrices, simultaneous sign successions, and tessellated pavements in two or more colours, with applications to Newton's rule, ornamental tile-work, and the theory of numbers, {\it Lond. Edinb. Dublin Philos. Mag. J. Sci.}, {\bf 34} (1867), 461--475.  

\bibitem{Szan}  
A. Szankowski, On Dvoretzky's theorem on almost spherical sectiona of convex bodies, {\it Israel J. Math.} {\bf 17} (1974), 325--338. 

\bibitem{Sz76}  
S. J. Szarek, On the best constants in the Khinchin inequality,  {\it Stud. Math.} {\bf 58} (1976), 197--208.  

\bibitem{Sz78}  
S. J. Szarek, On Kashin's almost Euclidean orthogonal decomposition of $\ell^1_n$, {\it Bull. Acad. Polon. Sci.  S\'er. Sci. Math. Astronom. Phys.} {\bf 26} (1978), 691--694.

\bibitem{Sz80}  
S. J. Szarek, Volume estimates and nearly Euclidean decomposition for normed spaces, {\it Semin. Anal. Fonct. 1979-1980}, Exp. No. 25 (1980), 8pp.

\bibitem{Sz83}  
S. J. Szarek, The finite dimensional basis problem, with an appendix on nets of Grassman manifold, {\it Acta Math.} {\bf 159} (1983), 153--179. 

\bibitem{Sz86}  
S. J. Szarek, On the existence and uniquness of complex structure and spaces with ``few" opertors, {\it Trans. Am. Math. Soc.} {\bf 239} (1986), 339--353.

\bibitem {Sza90}
S. J. Szarek, Spaces with large distances to $\ell_{\infty}^n$ and random matrices, {\it Amer. J. Math.} {\bf 112} (1990), 899--942.

\bibitem {Szar91} 
S. J. Szarek, On the geometry of the Banach-Mazur compactum, {\it Lecture Notes in Math.} {\bf 1470} (1991), 48--59. 

\bibitem {[Szar06]} 
S. J. Szarek,  Convexity, complexity and high dimensions, {\it Proceedings of the International Congress of Mathematicians (Madrid, 2006)}, Vol. II, 1599--1621. 
European Mathematical Society (EMS), Z\"urich, 2006.

\bibitem {Szar14} 
S. J. Szarek,  Coarse approximation of convex bodies by polytopes and the complexity of Banach--Mazur compacta, preprint, June 2014.

\bibitem {[SzTa1]}
S. J. Szarek and M. Talagrand, An ``isomorphic" version of the Sauer-Shelah lemma and the Banach-Mazur distance to the cube, GAFA - Israel Funct. Anal. Sem., {\it Lecture Notes in Math.} {\bf 1376}, Springer Verlag, 1989, 105--112.

\bibitem {[SzTa2]}  
S. J. Szarek and M. Talagrand, On the convexified Sauer-Shelah theorem, {\it J. Comb. Theory}, Ser. B, {\bf 69} (1997), No. 2, 183--192.

\bibitem {SzTom}  
S. J. Szarek and N. Tomczak-Jaegermann, On nearly Euclidean decomposition for some cleaa]ses of Banach spaces, {\it Compos. Math.} {\bf 40} (1980), 367--385.

\bibitem {Taschuk} 
B. Taschuk, The Banach-Mazur distance to the cube in low dimensions, {\it Discrete Comput. Geom.} {\bf 46} (2011), No. 1, 175--183.

\bibitem {Tcaciuc} 
A. Tcaciuc, Structural properties of finite and infinite dimensional Banach spaces, Ph. D. Thesis, University of Alberta, 2005. 

\bibitem {[Thompson]}
A. C. Thompson, {\it Minkowski Geometry.} Encyclopedia of Mathematics and its Applications {\bf 63}, Cambridge University Press, Cambridge, 1996. 

\bibitem {Tik19} 
K. Tikhomirov, On the Banach-Mazur distance to cross-polytope, {\it Adv. Math.} {\bf 345} (2019), 598--617.

\bibitem {[Tom 78]}  
N. Tomczak-Jaegermann, The Banach-Mazur distance between the trace classes $c\sp{n}\sb{p}$, {\it Proc. Amer. Math. Soc.} {\bf 72 } (1978), No. 2, 305--308.   

\bibitem {[Tom 79]} 
N. Tomczak-Jaegermann, On the Banach-Mazur distance between symmetric spaces, {\it Bull. Akad. Polon. Sci. S\'er. Sci. Math.} {\bf 27}. No. 3--4, 273--276.

\bibitem {Tom 81} 
N. Tomczak-Jaegermann, Dimensional gradations and cogradations of operator ideals. The weak distance between Banach spaces, {\it Semin. Anal. Fonct. 1980-81} (1981), Exp. No. 9,  Ecole Polytechnique, Palaiseau.

\bibitem{[Tom 81]}  
N. Tomczak-Jaegermann, The distance between a symmetric space and a $2$-convex or $2$-concave space, {\it J. London Math. Soc. (2)}, {\bf 24} (1981), No. 2, 272--282.

\bibitem {[Tom 83]} 
N. Tomczak-Jaegermann, The Banach-Mazur distance between symmetric spaces, {\it Israel J. Math.} {\bf 46} (1983), No. 1-2, 40--66.
	
\bibitem{[Tom 84]}  
N. Tomczak-Jaegermann, The weak distance between finite-dimensional Banach spaces, {\it Math. Nachr.} {\bf 119} (1984), 291--307.
	 
\bibitem {[Tom 89]} 
N. Tomczak-Jaegermann, {\it Banach-Mazur Distances and Finite-Dimensional Operator Ideals},  Pitman Monographs and Surveys in Pure and Applied Mathematics {\bf 38}, Longman Scientific \& Technical, New York, 1989.

\bibitem{To15}
G. Toth, {\it Measures of Symmetry for Convex Sets and Stability}, Springer International Publishing, Switzerland, 2015.

\bibitem{Toussaint83}
G. T. Toussaint, Solving Geometric Problems with the ``Rotating Calipers,'' {\it Proc. of IEEE MELECON 83}, Athens, Greece (May 1983).

\bibitem{Usher}  
M. Usher, Symplectic Banach-Mazur distances between subsets of $\mathcal C^n$, {\it J. Topol. Anal.} {\bf 14} (2022), No. 1, 231--286. 

\bibitem{Ver01}
R. Vershynin, John's decompositions: Selecting a large part, {\it Israel J. Math.} {\bf 122} (2001), 253--277.

\bibitem{Ver11}
R. Vershynin, Lectures in Geometric Functional Analysis (2011),   \url{https://www.math.uci.edu/~rvershyn/papers/GFA-book.pdf}.
 
\bibitem{Ver26}
R. Vershynin, High-Dimensional Probability: An Introduction with Applications in Data Science,  2nd Ed., Cambridge University Press, 2026. 

\bibitem{Villada}
J. Villada Bedoya, An estimate of the Banach-Mazur distance between the space of convergent sequences and a concrete model of a space of affine continuous functions, {\it Topol. Methods Nonlinear Anal.} {\bf 66} (2025), No. 2, 783--792.


\bibitem {VK}  
V. V. Voskyanan and R. Khaldi, The Banach-Mazur distance between some two-dimensional normed spaces (in Russian), {\it Izv. Akad. Armyan. SSR, Ser. Mat.} {\bf 94} (1993), No. 1, 207--220. 

\bibitem{wendel} 
J. G. Wendel, A problem in geometric probability, {\it Math. Scand.} {\bf 11} (1962), 109--111.

\bibitem {Xue 18} 
F. Xue, On the Banach--Mazur distance between the cube and crosspolytope, {\it Math. Ineq. Appl.} {\bf 21} (2018), No. 4, 931--943. 

\bibitem {Xue 19} 
F. Xue, Convex geometry of numbers: covering, succesive minima and Banach-Mazur distance, Ph. D. Thesis, Technische Universit\"at Berlin, 2019.

\bibitem {Youssef14} 
P. Youssef, Restricted invertibility and the Banach-Mazur distance to the cube, {\it Mathematika} {\bf 60}, No. 1, 201--218.

\bibitem {YN76} 
D. B. Yudin and A. S. Nemirovski\v{i}, Informational complexity and efficient methods for the solution of convex extremal problems. Matekon. 1976; {\bf 13} (2):22–45. 

\bibitem{Zha}
C. Zhang, Completly bounded Banach-Mazur distance, {\it Proc. Edinburgh Mat. Soc.} {\bf 40} (1977), No. 2, 247--260.   

\bibitem{ZhaHe}
X. Zhang and C. He, Partitioning bounding sets in symmetric spaces into subsets with reduced diameter, {\it Math. Inequal. Appl.} {\bf 27} (2024), 703--713.

\bibitem{ZhaMenWu}  
L. Zhang, L. Meng and S. Wu, Banach-Mazur distance from $\ell_p^n$ to $\ell_p^\infty$, {\it Math. Notes} {\bf 114} (2023), 1045--1051. 

\bibitem{Zip81} M. Zippin, The finite-dimensional $P_\lambda$ spaces with small $\lambda$,  {\it Israel J. Math.} {\bf 39} (1981), no. 4, 359--364.

\bibitem{Zip84} M. Zippin, Correction to ``The finite-dimensional $P_\lambda$ spaces with small $\lambda$,''  {\it Israel J. Math.} {\bf 48} (1984), no. 2-3, 255--256.




\end{thebibliography}
\end{document}